\documentclass[11pt]{amsart}

\usepackage[T1]{fontenc}
\usepackage{amsmath,amssymb,amsthm,mathtools}
\mathtoolsset{showonlyrefs}
\usepackage[margin=1.15in]{geometry}
\usepackage{mathrsfs}
\usepackage{enumitem}
\usepackage{microtype}
\usepackage[colorlinks=true,citecolor=blue,linkcolor=blue,urlcolor=blue]{hyperref}

\numberwithin{equation}{section}
\allowdisplaybreaks[2]

\newtheorem{theorem}{Theorem}[section]
\newtheorem{proposition}[theorem]{Proposition}
\newtheorem{lemma}[theorem]{Lemma}
\newtheorem{corollary}[theorem]{Corollary}
\theoremstyle{definition}
\newtheorem{definition}[theorem]{Definition}
\newtheorem{example}[theorem]{Example}

\newcommand{\C}{\mathbb C}
\newcommand{\R}{\mathbb R}
\newcommand{\Sph}{\mathbb S}
\newcommand{\Pcal}{\mathscr P}
\newcommand{\Ric}{\operatorname{Ric}}
\newcommand{\Rm}{\operatorname{Rm}}
\newcommand{\spt}{\operatorname{spt}}
\newcommand{\dist}{\operatorname{dist}}
\newcommand{\tr}{\operatorname{tr}}
\newcommand{\ImmReg}{\operatorname{Reg}_{\mathrm{imm}}}
\newcommand{\ImmSing}{\operatorname{Sing}_{\mathrm{imm}}}
\newcommand{\eps}{\varepsilon}
\providecommand{\llbracket}{[\![}
\providecommand{\rrbracket}{]\!]}

\title[Rigidity, regularity and compactness under Simons stability]{Rigidity, Regularity and Compactness of Special Lagrangian Submanifolds under Simons Stability}
\author{Caiyan Li}
\address{(C.L.) School of Mathematical Sciences, Xiamen University, 361005, Xiamen, P.~R.~China}
\email{caiyanli@xmu.edu.cn} 
\thanks{C.L. is supported by NSFC (Grant No. 12501274).}

\date{}

\begin{document}

\begin{abstract}
In this paper, we prove that, for suitable values of $\delta$, every
complete connected $\delta$-Simons stable special Lagrangian immersion
in $\C^n$, $n\ge3$, is flat. We also prove that a
$\delta$-Simons stable special Lagrangian cone $C\subset\C^n$ with
$C\cap\Sph^{2n-1}$ smooth and compact is a finite union of planes if
$\delta>(n-2)^2/(4n)$, and that this coefficient is optimal.
At equality, we characterize the $\delta$-Simons stable cones. We further establish regularity and compactness
for sequences of proper $\delta$-Simons stable special Lagrangian immersions
when $n\ge3$ and $\delta>(n-2)/(n+1)$, under suitable restrictions on
the initial singular sets, assuming integral current extensions $T_j$
with $\partial T_j=0$ and uniformly bounded mass on compact subsets. In dimension three, we allow initial
singular sets $S_j$ satisfying $\mathcal H^1(S_j)=0$. For suitable values
of $\delta$, we also prove properness, Euclidean volume growth, and finite
total curvature for complete connected immersions of finite $\delta$-index.
\end{abstract}

\maketitle

\section{Introduction}

Harvey--Lawson \cite{HarveyLawson1982} introduced special Lagrangian
submanifolds of $\C^n$ in the theory of calibrations.
A real submanifold $M^n\subset\C^n\simeq\R^{2n}$ is Lagrangian if
$J(T_pM)=(T_pM)^\perp$ for all $p\in M$, where $J$ is the standard
complex structure on $\C^n$.
An oriented Lagrangian submanifold $M$ is called special Lagrangian
of phase $\theta$ if it satisfies the calibration condition
\[
\operatorname{Re}\bigl(e^{-i\theta}dz_1\wedge\cdots\wedge dz_n\bigr)
\big|_M=\operatorname{vol}_M
\qquad\text{for a constant }\theta\in\R.
\]
Here $z_1,\ldots,z_n$ are the standard complex coordinates and
$\operatorname{vol}_M$ is the induced volume form.
Thus $M$ is calibrated and locally volume minimizing, and hence minimal.

For entire special Lagrangian gradient graphs, Yuan proved that the
potential is quadratic under convexity \cite{Yuan2002} or supercritical
phase \cite{Yuan2006}. Tsui--Wang \cite{TsuiWang2002} obtained Bernstein
theorems for complete minimal Lagrangian submanifolds under restrictions
on the Gauss image.

Harvey--Lawson also constructed the nonflat Clifford torus cones
\cite[\S III.3.A]{HarveyLawson1982}. In $\C^3$, Haskins
\cite{Haskins2004} constructed infinitely many special Lagrangian cones
over embedded tori and proved that every special Lagrangian cone over
a smoothly immersed two-sphere is flat
\cite[Theorem~B]{Haskins2004}.
For special Lagrangian cones, Li
\cite[Theorems~1 and~2]{Li2021} proved rigidity under a flat normal
bundle assumption when $2\le n\le7$, provided the cone is graphical,
and under the condition that the Gauss image lies in an open geodesic
ball of radius $\sqrt2\pi/4$ in the oriented Grassmannian when
$2\le n\le8$. Ogden \cite{Ogden2026} proved in arbitrary dimension that
a connected Euclidean submanifold calibrated by a form with constant coefficients and having flat normal bundle is flat.
Lee--Wang \cite[Theorem~1.4]{LeeWang2026} proved that a complete connected
minimal Lagrangian submanifold satisfying the almost isoclinic condition
is flat if the canonical function on the Gauss image is uniformly positive.

An oriented immersed minimal hypersurface $M^n\subset\R^{n+1}$ is called
stable if
\[
\int_M\bigl(|\nabla\phi|^2-|A_M|^2\phi^2\bigr)\,d\mu\ge0,
\qquad \forall\,\phi\in C_c^\infty(M),
\]
where $A_M$ denotes the second fundamental form. Simons
\cite[\S6.1]{Simons1968} proved that a stable minimal hypercone
$C^n\subset\R^{n+1}$ with $C\cap\Sph^n$ a smooth closed immersed
hypersurface is flat for $n\le6$. For $3\le n\le5$, Schoen--Simon--Yau
\cite[Theorems~2 and~3]{SchoenSimonYau1975} obtained pointwise curvature
estimates depending on the volume ratios $R^{-n}\mu(M\cap B_R(x))$. In particular, complete stable
minimal hypersurfaces in $\R^{n+1}$ with Euclidean volume growth are flat.

More generally, for $\delta>0$, $M$ is called $\delta$-stable if
\[
\int_M\bigl(|\nabla\phi|^2-\delta|A_M|^2\phi^2\bigr)\,d\mu\ge0,
\qquad \forall\,\phi\in C_c^\infty(M).
\]
For complete $\delta$-stable minimal hypersurfaces
with $n\ge3$ and $\delta>((n-1)/n)^2$, Fu
\cite[Theorem~1.1 and Corollary~1.2]{Fu2011} proved the
vanishing of $L^2$ harmonic one forms and compactly supported first
cohomology, as well as the existence of exactly one end. Li
\cite{Li2020,Li2021} proved cone rigidity using scalar stability,
curvature gradient estimates, and radial Hardy inequalities.
Hong--Li--Wang \cite[Theorem~1.1]{HongLiWang2024} established regularity
and compactness for sequences of oriented immersed $\delta$-stable
minimal hypersurfaces with $n\ge3$ and $\delta>(n-2)/n$, assuming
uniformly bounded local mass and
$\dim_{\mathcal H}S_j\le\sigma<n-4+4/n$ for the initial singular sets $S_j$.
Minter--Xiao \cite[Theorem~A]{MinterXiao2026} proved regularity and
compactness for oriented stable minimal hypersurfaces properly immersed
outside closed sets $S_j$ with $\mathcal H^{n-2}(S_j)=0$, with uniformly
bounded local mass.

In this paper, we prove Bernstein theorems, cone rigidity, and regularity
and compactness for special Lagrangian immersions under a scalar stability
condition, and study the cones at the critical coefficient and complete
immersions of finite scalar index.

For a smooth special Lagrangian immersion $F:M^n\to\C^n$, $n\ge2$, let
$g$ be the induced metric and $A_M$ the second fundamental form.
The quartic term in the Simons identity \eqref{eq:exact-simons} is
$|\Ric_g|^2+|\Rm_g|^2$. We use the potential
\begin{equation}\label{eq:simons-potential}
\Pcal\coloneqq 
\begin{cases}
\dfrac{|\Ric_g|^2+|\Rm_g|^2}{|A_M|^2},& |A_M|\ne0,\\[2mm]
0,& |A_M|=0.
\end{cases}
\end{equation}

\begin{definition}\label{def:simons-stability}
The immersion $F$ is called \emph{$\delta$-Simons stable} if
\begin{equation}\label{eq:delta-stability}
\delta\int_M\Pcal\phi^2\,d\mu
\le \int_M|\nabla\phi|^2\,d\mu,
\qquad \forall\,\phi\in C_c^\infty(M).
\end{equation}
For the quadratic form
\[
Q_{\delta,M}(\phi)
\coloneqq \int_M|\nabla\phi|^2\,d\mu
-\delta\int_M\Pcal\phi^2\,d\mu,
\quad \phi\in C_c^\infty(M),
\]
we define the $\delta$-index by
\[
\operatorname{ind}_{\delta}(M)
\coloneqq \sup\bigl\{\dim V:
Q_{\delta,M}(\phi)<0\quad\forall\,\phi\in V\setminus\{0\}\bigr\},
\]
where $V$ ranges over the finite dimensional subspaces of $C_c^\infty(M)$.
\end{definition}

Condition~\eqref{eq:delta-stability} differs from spectral stability of
special Lagrangian cones and from normal Jacobi stability. We also recall
that, for a closed special Lagrangian cone $C\subset\C^n$, $n>2$, with an
isolated singularity at the vertex, Joyce \cite[Definition~3.6]{Joyce2004} defined
a stability index using the Laplace spectrum of $C\cap\Sph^{2n-1}$.
Among special Lagrangian torus cones in $\C^3$, Haskins proved that the
Clifford cone is the only one with zero stability index, up to special
unitary equivalence \cite[Theorem~A]{HaskinsComplexity2004}.
Ohnita \cite{Ohnita2007} computed the stability index and the Legendrian
index for cones over certain compact minimal Legendrian orbits.
Dimler--Lee \cite{DimlerLee2025} studied strict normal Jacobi stability
of calibrated cones.

By Proposition~\ref{prop:curvature}, $\operatorname{Scal}_g=-|A_M|^2$ and
\[
\frac{n+1}{n(n-1)}|A_M|^2\le\Pcal\le\frac32|A_M|^2.
\]
Thus stability with potential $\kappa|A_M|^2$ implies
$\delta$-Simons stability for $\delta\le2\kappa/3$.
Conversely, $\delta$-Simons stability implies stability with potential
$\delta(n+1)|A_M|^2/(n(n-1))$.

We first state the Bernstein theorem for complete immersions.  Let
\[
\delta_{\rm K}(n) \coloneqq \frac{n-2}{n+1}, \quad
\delta_{\rm cone}(n) \coloneqq \frac{(n-2)^2}{4n}.
\]
For $n\ge3$, let $D_n$ be the constant in
\eqref{eq:all-dimensional-constant}, with
\[
\max\{\delta_{\rm K}(n),\delta_{\rm cone}(n)\}\le D_n<14n^3.
\]
In Appendix~\ref{app:coefficients}, we construct $D_n$ and prove the local
curvature estimate for $\delta>D_n$
(Proposition~\ref{prop:stable-intrinsic-ball}).
Here $H_c^1(M;\R)$ denotes cohomology with compact support.

\begin{theorem}\label{thm:intro-bernstein}
Let $F:M^n\to\C^n$, $n\ge3$, be a complete connected special
Lagrangian immersion without boundary. Assume that $F$ is
$\delta_+$-Simons stable and that one of the following conditions holds:
\begin{enumerate}[label=\textup{(\roman*)}]
\item
\[
\begin{cases}
\delta_+\ge(5+3\sqrt5)/40,&n=3,\\
\delta_+\ge7/12,&n=4,\\
\delta_+\ge3/2,&n=5,\\
\delta_+\ge7/2,&n=6,\\
\delta_+\ge10,&n=7.
\end{cases}
\]
\item $n=3$, $\delta_+>1/4$, and $\dim H_c^1(M;\R)<\infty$.
\item $\delta_+>D_n$.
\end{enumerate}
Then $F$ is flat.
\end{theorem}

For an immersed special Lagrangian cone $C^m\subset\C^m$, $A_C$ denotes
the second fundamental form and $\Sigma$ denotes the domain of the
immersion into $C\cap\Sph^{2m-1}$.

\begin{theorem}\label{thm:intro-cone}
Let $C^m\subset\C^m$, $m\ge3$, be a special Lagrangian cone with
$\Sigma\coloneqq C\cap\Sph^{2m-1}$ smooth and compact. If the immersion domain of
$C\setminus\{0\}$ is $\delta$-Simons stable for
$\delta>\delta_{\rm cone}(m)$, then $C$ is a finite union of special
Lagrangian $m$-planes with the same phase and positive integer
multiplicities.
\end{theorem}

The coefficient $\delta_{\rm cone}(m)$ in
Theorem~\ref{thm:intro-cone} is sharp. Indeed, the Clifford torus cone
$\mathrm{Cliff}^{m}$ is nonflat and $\delta_{\rm cone}(m)$-Simons stable
(Example~\ref{ex:clifford-cone}).
For a four dimensional cone as in Theorem~\ref{thm:intro-cone}, suppose
\[
\alpha\int_C|A_C|^2\phi^2\,d\mu_C
\le\int_C|\nabla\phi|^2\,d\mu_C,
\qquad \forall\,\phi\in C_c^1(C\setminus\{0\}).
\]
By $\Pcal\le3|A_C|^2/2$, $C$ is $2\alpha/3$-Simons stable.
Thus $C$ is flat for $\alpha>3/8$ by Theorem~\ref{thm:intro-cone};
the sufficient condition in \cite[Theorem~1.2]{Li2020} is $\alpha>1$.

We next consider cones at $\delta=\delta_{\rm cone}(n)$.
For a cone $C\subset\C^n$ with smooth, compact, connected
$\Sigma\coloneqq C\cap\Sph^{2n-1}$, define at $r=1$
\begin{equation}\label{eq:cone-tensor-definition}
(h_\Sigma)_p(X,Y,Z)\coloneqq (h_C)_{(1,p)}(X,Y,Z)
=\langle A_\Sigma^{\mathbb S}(X,Y),JZ\rangle,
\end{equation}
for $p\in\Sigma$ and $X,Y,Z\in T_p\Sigma$. Here $h_C$ is the tensor in \eqref{eq:lagrangian-cubic} for $C$, and $A_\Sigma^{\mathbb S}$ is the second
fundamental form of $\Sigma\subset\Sph^{2n-1}$. We compute
$\nabla^\Sigma h_\Sigma$ using the Levi--Civita connection of $g_\Sigma$
in all three arguments. By Theorem~\ref{thm:critical-cone-non4},
$C$ is $\delta_{\rm cone}(n)$-Simons stable if and only if
$\nabla^\Sigma h_\Sigma=0$ when $n\ne4$; for $n=4$, stability is
equivalent to $h_\Sigma\equiv0$ or
\eqref{eq:critical-four-characterization}.

We also prove that a global special Lagrangian gradient graph cone in
$\C^n$, $n\ge3$, that is smooth away from the vertex and has
finite $\delta_{\rm cone}(n)$-index is flat
(Corollary~\ref{cor:critical-graphical-cone}).

We next consider complete immersions of finite scalar index.
For complete minimal hypersurfaces of dimensions three through six,
Tysk \cite{Tysk1989} proved that finite Morse index is equivalent to
finite total curvature under Euclidean volume growth. Chodosh--Li
\cite{ChodoshLi2024} and Chodosh--Li--Minter--Stryker
\cite{ChodoshLiMinterStryker2025} removed the growth assumption in
dimensions three and four, respectively. For $n=3,4$, Ranieri
\cite[Theorem~1.2]{Ranieri2026} established
the equivalence with finite $\delta$-index in place of finite Morse index
when $\delta>((n-1)/n)^2$.

In higher codimension, Wang \cite{Wang2003} introduced super stability
and the super index with potential $|A_M|^2$. For a complete connected
minimal submanifold of dimension at least three, Wang proved that super
stability and finite total scalar curvature imply that the submanifold
is flat, and that finite super index implies finitely many
ends. With the potential $\Pcal$, we obtain the following result.

\begin{theorem}\label{thm:intro-finite-index}
Let $F:M^n\to\C^n$, $n\ge3$, be a complete connected special
Lagrangian immersion without boundary. Assume
$\operatorname{ind}_{\delta_+}(M)<\infty$ and that one of the following holds:
\begin{enumerate}[label=\textup{(\roman*)}]
\item
\begin{equation}\label{eq:low-dimensional-index-range}
\begin{cases}
\delta_+>(5+3\sqrt5)/30,&n=3,\\
\delta_+>(9+6\sqrt6)/40,&n=4,\\
\delta_+\ge3/2,&n=5,\\
\delta_+\ge7/2,&n=6,\\
\delta_+\ge10,&n=7,
\end{cases}
\end{equation}
\item $\dim H_c^1(M;\R)<\infty$ and either $n=3$ with $\delta_+>1/4$,
or $n=4$ with $\delta_+\ge7/12$.
\item $\delta_+>D_n$.
\end{enumerate}
Then $F$ is proper,
$\int_M|A_M|^n\,d\mu<\infty$, and there exist constants $\Lambda,C,R_0>0$ such that
\[
\begin{aligned}
\mu\bigl(F^{-1}(B_R(0))\bigr)&\le\Lambda R^n &&\forall\,R\ge1,\\
|A_M|(x)&\le C|F(x)|^{-n-1} &&\forall\,x\in\{|F|\ge R_0\}.
\end{aligned}
\]
\end{theorem}

We finally state the regularity and compactness theorem. For $\delta>0$,
define
\begin{equation}\label{eq:compactness-critical-dimension}
n_\delta\coloneqq
\min\bigl\{m\in\mathbb N:m\ge3,\ \delta\le\delta_{\rm cone}(m)\bigr\}.
\end{equation}
A point of $\spt T$ is immersed regular if $T$ is locally a
finite sum of smooth proper special Lagrangian immersions with positive
integer multiplicities (Definition~\ref{def:imm-reg}). We denote
$\spt T\setminus\ImmReg T$ by $\ImmSing T$.

\begin{theorem}\label{thm:intro-regularity-compactness}
Let $n\ge3$, $\delta>\delta_{\rm K}(n)$, and let $S_j\subset B_2$
be relatively closed sets satisfying
\begin{equation}\label{eq:initial-singular-dimension}
\dim_{\mathcal H}S_j\le\sigma<n-4+\frac6{n+1}\qquad\forall\,j\ge1.
\end{equation}
Suppose $F_j:M_j^n\to B_2\setminus S_j\subset\C^n$ are smooth proper
special Lagrangian immersions without boundary over $B_2\setminus S_j$,
with a common phase $\Theta_j$ on all components.
Assume that each $F_j$ is $\delta$-Simons stable, the current of $F_j$ extends to
an integral current $T_j$ with $\partial T_j=0$ in $B_2$, and
$\sup_j\|T_j\|(B_2)\le\Lambda<\infty$.
Then there exists a subsequence such that $\Theta_j\to\Theta$ modulo $2\pi$ and
\[
T_j\rightharpoonup T\quad\text{locally as currents},\qquad
|T_j|\longrightarrow|T|\quad\text{locally as varifolds as }j\to\infty,
\]
where $T$ is a special Lagrangian integral current of phase $\Theta$.
For $S\coloneqq\ImmSing T\cap B_1$,
\[
\begin{cases}
S=\varnothing,&n<n_\delta,\\
S\text{ is locally finite in }B_1,&n=n_\delta,\\
\dim_{\mathcal H}S\le n-n_\delta,&n>n_\delta.
\end{cases}
\]
Moreover, along the same subsequence, the extensions of $F_j$ across
the removable singularities in $S_j$ converge smoothly as immersions on
compact subsets of $B_1\setminus S$ to a proper $\delta$-Simons stable
special Lagrangian immersion of phase $\Theta$ representing
$T\llcorner(B_1\setminus S)$.
\end{theorem}

In dimension three, we replace \eqref{eq:initial-singular-dimension}
by $\mathcal H^1(S_j)=0$ and prove local smooth immersed convergence
for $\delta>1/4$
(Theorem~\ref{thm:three-dimensional-null-set-compactness}).

For $\delta=1$, we have $n_\delta=8$, so the current $T$ in
Theorem~\ref{thm:intro-regularity-compactness} satisfies
$\ImmSing T\cap B_1=\varnothing$ for $n\le7$ and
$\ImmSing T\cap B_1$ is locally finite for $n=8$. Since
$\delta_{\rm K}(n)>\delta_{\rm cone}(n)$ for $3\le n\le5$, $S=\varnothing$ throughout
the local regularity range. In particular, for
$\delta>\max\{\delta_{\rm K}(n),\delta_{\rm cone}(n)\}$,
we obtain the following curvature estimate and Bernstein theorem under
Euclidean volume growth.

\begin{corollary}\label{cor:curvature-and-growth-bernstein}
Let $n\ge3$ and
$\delta>\max\{\delta_{\rm K}(n),\delta_{\rm cone}(n)\}$.
For each $0<\Lambda<\infty$, there exists $C=C(n,\delta,\Lambda)$ such
that any smooth proper $\delta$-Simons stable special Lagrangian
immersion $F:M^n\to B_2\subset\C^n$, without boundary over $B_2$,
with a common phase on all components and $\mu(M)\le\Lambda$, satisfies
\begin{equation}\label{eq:curvature-mass-bound}
\sup_{F^{-1}(B_1)}|A_M|\le C.
\end{equation}
Moreover, any complete connected $\delta$-Simons stable special
Lagrangian immersion $F:M^n\to\C^n$ without boundary satisfying
\[
\mu\bigl(F^{-1}(B_R(0))\bigr)\le\Lambda R^n\qquad\forall\,R\ge1
\]
is flat.
\end{corollary}

We prove the local curvature estimates in
Section~\ref{sec:local-estimates}. In Section~\ref{sec:cones}, we prove
cone rigidity and characterize stability at the critical coefficient.
Theorem~\ref{thm:intro-regularity-compactness} and
Corollary~\ref{cor:curvature-and-growth-bernstein} are proved in
Section~\ref{sec:immersed-compactness}. In
Section~\ref{sec:small-singular-sets}, we extend the compactness theorem
to $\mathcal H^1$-null initial singular sets in dimension three.
Sections~\ref{sec:complete} and~\ref{sec:compactness} contain the finite
total curvature, finite index, and Bernstein arguments.
We prove Proposition~\ref{prop:stable-intrinsic-ball} in
Appendix~\ref{app:coefficients}.

\section{Curvature estimates}\label{sec:local-estimates}

Let
$M^n\subset\C^n$ be a smooth special Lagrangian immersion with constant
phase on each connected component. Let $g$ be the induced metric,
$\overline\nabla$ the Euclidean connection, and $\nabla$ the Levi-Civita
connection and the induced connections on tensor bundles. We use the
Laplace--Beltrami operator $\Delta\coloneqq \operatorname{div}_g\nabla$. For smooth tangent
vector fields $X,Y,Z,U$ on $M$, our conventions for the second fundamental
form $A_M$ and the Riemann curvature tensor $\Rm_g$ are
\begin{equation}\label{eq:curvature-conventions}
	\begin{aligned}
		A_M(X,Y)&\coloneqq (\overline\nabla_XY)^\perp,\\
		\Rm_g(X,Y,Z,U)
		&\coloneqq \bigl\langle\nabla_X\nabla_YU-\nabla_Y\nabla_XU
		-\nabla_{[X,Y]}U,Z\bigr\rangle.
	\end{aligned}
\end{equation}
Let
\begin{equation}\label{eq:lagrangian-cubic}
	h_M(X,Y,Z)\coloneqq \langle A_M(X,Y),JZ\rangle.
\end{equation}
Fix a local orthonormal tangent frame $e_1,\dots,e_n$ and set
$(h_M)_{ijk}\coloneqq h_M(e_i,e_j,e_k)$ and
$(\Rm_g)_{ijk\ell}\coloneqq \Rm_g(e_i,e_j,e_k,e_\ell)$. The Ricci tensor and sectional
curvature are $(\Ric_g)_{ij}\coloneqq \sum_k(\Rm_g)_{ikjk}$ and
$K_g(e_i,e_j)\coloneqq (\Rm_g)_{ijij}$ for $i\ne j$, respectively. The scalar curvature is
$\operatorname{Scal}_g\coloneqq \sum_i(\Ric_g)_{ii}$, and
$|\Rm_g|^2\coloneqq \sum_{i,j,k,\ell}(\Rm_g)_{ijk\ell}^2$.

\begin{proposition}\label{prop:curvature}
Let $M^n\subset\C^n$, $n\ge2$, be a smooth special Lagrangian
immersion. Then
\begin{gather}
\frac12\Delta|A_M|^2 =|\nabla A_M|^2-|\Ric_g|^2-|\Rm_g|^2,\label{eq:exact-simons}\\
\frac{n+1}{n(n-1)}|A_M|^2 \le\Pcal\le\frac32|A_M|^2.\label{eq:P-estimates}
\end{gather}
\end{proposition}

\begin{proof}
Using $(h_M)_{ijk}=(h_M)_{jik}=(h_M)_{ikj}$ and
$\sum_a (h_M)_{aap}=0$ in the Gauss equation with
\eqref{eq:curvature-conventions}, and then contracting
$(\Rm_g)_{ikjk}$ over $k$, we obtain
\begin{equation}\label{eq:ric-matrices}
\begin{aligned}
(\Rm_g)_{ijk\ell}
&=\sum_p\bigl((h_M)_{ikp}(h_M)_{j\ell p}-(h_M)_{i\ell p}(h_M)_{jkp}\bigr)
 =[\mathsf S_i,\mathsf S_j]_{k\ell},\\
(\Ric_g)_{ij}&=-\sum_{a,p}(h_M)_{iap}(h_M)_{jap}
 =-\tr(\mathsf S_i\mathsf S_j),
\end{aligned}
\end{equation}
where $\mathsf S_k\coloneqq ((h_M)_{ijk})_{i,j=1}^n$.
Consequently,
\begin{equation}\label{eq:comm-Rm}
\sum_{i,j}\bigl(\tr(\mathsf S_i\mathsf S_j)\bigr)^2=|\Ric_g|^2,
\quad
\sum_{i,j,k,\ell}[\mathsf S_i,\mathsf S_j]_{k\ell}^2=|\Rm_g|^2.
\end{equation}
In the normal frame $Je_1,\ldots,Je_n$, by
\cite[p.~63, Eq.~(2.22)]{ChernDoCarmoKobayashi1970} with zero ambient curvature,
\[
\Delta\mathsf S_i
=-\sum_j\tr(\mathsf S_i\mathsf S_j)\mathsf S_j
 -\sum_j[\mathsf S_j,[\mathsf S_j,\mathsf S_i]].
\]
Taking the inner product with $\mathsf S_i$,
summing over $i$, and using \eqref{eq:comm-Rm} and
$\langle\mathsf S_i,[\mathsf S_j,[\mathsf S_j,\mathsf S_i]]\rangle
=\sum_{k,\ell}[\mathsf S_i,\mathsf S_j]_{k\ell}^2$, we obtain
\eqref{eq:exact-simons}.

Tracing \eqref{eq:ric-matrices}, we obtain
$\operatorname{Scal}_g=2\sum_{i<j}(\Rm_g)_{ijij}=-|A_M|^2$.
By Cauchy--Schwarz and the skew symmetry of $\Rm_g$ in each pair of
indices,
\[
\begin{aligned}
|\Ric_g|^2\ge&\frac1n\Bigl(\sum_i(\Ric_g)_{ii}\Bigr)^2
=\frac1n|A_M|^4,\\
|\Rm_g|^2\ge&4\sum_{i<j}(\Rm_g)_{ijij}^2 
\ge \frac8{n(n-1)}\Bigl(\sum_{i<j}(\Rm_g)_{ijij}\Bigr)^2
=\frac2{n(n-1)}|A_M|^4.
\end{aligned}
\]
Together with \eqref{eq:comm-Rm} and \cite[Eq.~(5)]{LiLi1992},
\[
\frac{n+1}{n(n-1)}|A_M|^4
\le |\Ric_g|^2+|\Rm_g|^2
\le\frac32|A_M|^4.
\]
Thus \eqref{eq:P-estimates} follows after division by $|A_M|^2$ on
$\{|A_M|>0\}$, with $\Pcal=0$ on $\{A_M=0\}$.

\end{proof}

By \eqref{eq:exact-simons} and the refined Kato inequality
\cite[Lemma~5.17, (5.39)]{Fox2021} with $k=3$ for $h_M$ in
\eqref{eq:lagrangian-cubic},
\begin{equation}\label{eq:scalar-simons}
\begin{aligned}
\frac12\Delta|A_M|^2+\Pcal|A_M|^2
=&|\nabla A_M|^2=|\nabla h_M|^2 
\ge \frac{n+4}{n+1}|\nabla|h_M||^2
=\frac{n+4}{n+1}\bigl|\nabla|A_M|\bigr|^2.
\end{aligned}
\end{equation}
\begin{lemma}\label{lem:caccioppoli}
Let $M^n\subset\C^n$ be a smooth special Lagrangian immersion.
If $M$ is $\delta$-Simons stable and
$\delta_{\rm K}(n)<\alpha<\min\{\delta,1\}$, then
\begin{align}
&\int_{M\cap\{|A_M|^\alpha>k\}}\Bigl(1-\frac{k}{|A_M|^\alpha}\Bigr)
|\nabla(|A_M|^\alpha)|^2\eta^2\,d\mu\notag\\
\le &  C\int_{M\cap\{|A_M|^\alpha>k\}}(|A_M|^\alpha-k)^2|\nabla\eta|^2\,d\mu
+Ck^{2+2/\alpha}
\int_{M\cap\{k<|A_M|^\alpha<\delta k/(\delta-\alpha)\}}\eta^2\,d\mu,
\label{eq:level-cacc}
\end{align}
for $k\ge0$ and $\eta\in\operatorname{Lip}_c(M)$, where
$C=C(n,\delta,\alpha)$.
\end{lemma}

\begin{proof}
On $\{|A_M|>0\}$, by \eqref{eq:scalar-simons},
\begin{equation}\label{eq:power-ineq}
\Delta(|A_M|^\alpha)+\alpha\Pcal |A_M|^\alpha
\ge\frac{\alpha-\delta_{\rm K}(n)}{\alpha|A_M|^\alpha}
|\nabla(|A_M|^\alpha)|^2.
\end{equation}
For $\phi\in\operatorname{Lip}_c(M)$, choose $\phi_j\in C_c^\infty(M)$
with supports in a fixed compact subset of $M$ and
$\phi_j\to\phi$ in $W^{1,2}(M)$ as $j\to\infty$. By
\eqref{eq:delta-stability},
\begin{equation}\label{eq:lipschitz-stability}
\begin{aligned}
\delta\int_M\Pcal\phi^2\,d\mu
=&\lim_{j\to\infty}\delta\int_M\Pcal\phi_j^2\,d\mu 
\le \lim_{j\to\infty}\int_M|\nabla\phi_j|^2\,d\mu
=\int_M|\nabla\phi|^2\,d\mu.
\end{aligned}
\end{equation}
For $k>0$, multiplying \eqref{eq:power-ineq} by
$(|A_M|^\alpha-k)^+\eta^2$ and integrating by parts over $M$, we obtain,
by \eqref{eq:lipschitz-stability} with $\phi=(|A_M|^\alpha-k)^+\eta$,
\begin{equation}\label{eq:truncation-integral}
\begin{aligned}
&\frac{\alpha-\delta_{\rm K}(n)}{\alpha}\int_{\{|A_M|^\alpha>k\}}\Bigl(1-\frac{k}{|A_M|^\alpha}\Bigr)
|\nabla(|A_M|^\alpha)|^2\eta^2\,d\mu\\
\le&\int_{\{|A_M|^\alpha>k\}} 
 \alpha\Pcal|A_M|^\alpha(|A_M|^\alpha-k)\eta^2 -|\nabla((|A_M|^\alpha-k)\eta)|^2
 +(|A_M|^\alpha-k)^2|\nabla\eta|^2 \,d\mu\\
\le&\int_{\{|A_M|^\alpha>k\}}(|A_M|^\alpha-k)^2|\nabla\eta|^2  + \Pcal(|A_M|^\alpha-k)
\bigl(\delta k-(\delta-\alpha)|A_M|^\alpha\bigr)\eta^2\,d\mu.
\end{aligned}
\end{equation}
By Young's inequality and \eqref{eq:P-estimates},
\begin{equation}\label{eq:truncation-positive-part}
\begin{aligned}
&\Pcal(|A_M|^\alpha-k)^+\bigl(\delta k-(\delta-\alpha)|A_M|^\alpha\bigr)\\
\le&\frac{3\alpha^2}{8(\delta-\alpha)}
 \Bigl(\frac{\delta}{\delta-\alpha}\Bigr)^{2/\alpha}k^{2+2/\alpha}
 \mathbf1_{\{k<|A_M|^\alpha<\delta k/(\delta-\alpha)\}},
\end{aligned}
\end{equation}
Thus, for $k>0$, we obtain \eqref{eq:level-cacc} from
\eqref{eq:truncation-integral} and \eqref{eq:truncation-positive-part}.

For $k=0$, set $f_\eps\coloneqq (|A_M|^2+\eps^2)^{1/2}$, $\eps>0$.
By \eqref{eq:scalar-simons} and the chain rule,
\begin{equation}\label{eq:regularized-power}
\begin{aligned}
\Delta(f_\eps^\alpha)+\alpha\Pcal f_\eps^\alpha
&=\alpha f_\eps^{\alpha-2}
 \bigl(|\nabla A_M|^2+(\alpha-2)|\nabla f_\eps|^2+\eps^2\Pcal\bigr)\\
&\ge\frac{\alpha-\delta_{\rm K}(n)}{\alpha f_\eps^\alpha}
|\nabla(f_\eps^\alpha)|^2
 +\alpha\eps^2\Pcal f_\eps^{\alpha-2},
\end{aligned}
\end{equation}
Multiplying \eqref{eq:regularized-power} by $f_\eps^\alpha\eta^2$ and
integrating by parts over $M$, we obtain
\begin{equation}\label{eq:regularized-energy}
\begin{aligned}
&\alpha\int_M\Pcal f_\eps^{2\alpha}\eta^2\,d\mu\\
\ge&\frac{2\alpha-\delta_{\rm K}(n)}{\alpha}
\int_M|\nabla(f_\eps^\alpha)|^2\eta^2\,d\mu
+2\int_M f_\eps^\alpha\eta
\langle\nabla(f_\eps^\alpha),\nabla\eta\rangle\,d\mu\\
\ge&\frac{\alpha-\delta_{\rm K}(n)}{\alpha}
\int_M|\nabla(f_\eps^\alpha)|^2\eta^2\,d\mu
+\delta\int_M\Pcal f_\eps^{2\alpha}\eta^2\,d\mu
-\int_Mf_\eps^{2\alpha}|\nabla\eta|^2\,d\mu.
\end{aligned}
\end{equation}
Here the second inequality follows from \eqref{eq:lipschitz-stability}
with $\phi=f_\eps^\alpha\eta$. Since
$\delta>\alpha>\delta_{\rm K}(n)$ and
$f_\eps^\alpha\to|A_M|^\alpha$ locally uniformly,
\eqref{eq:regularized-energy} implies
$f_\eps^\alpha\rightharpoonup|A_M|^\alpha$ in $W^{1,2}_{\rm loc}(M)$
as $\eps\downarrow0$. By weak lower semicontinuity,
\begin{equation}\label{eq:zero-level-energy}
\begin{aligned}
&\frac{\alpha-\delta_{\rm K}(n)}{\alpha}\int_M|\nabla(|A_M|^\alpha)|^2\eta^2\,d\mu
 +(\delta-\alpha)\int_M\Pcal|A_M|^{2\alpha}\eta^2\,d\mu 
\le \int_M|A_M|^{2\alpha}|\nabla\eta|^2\,d\mu.
\end{aligned}
\end{equation}
Thus \eqref{eq:level-cacc} holds for $k=0$.
\end{proof}

\begin{proposition}\label{thm:eps-reg}
Let $S\subset B_{2r}(x)$ be relatively closed, and let
$F:M^n\to B_{2r}(x)\setminus S\subset\C^n$, $n\ge3$, be a smooth
proper special Lagrangian immersion without boundary over $B_{2r}(x)\setminus S$.
Assume that $|F|$ extends to a stationary integral varifold on $B_{2r}(x)$,
$F$ is $\delta$-Simons stable, and
$\delta_{\rm K}(n)<\alpha<\min\{\delta,1\}$. If $S\ne\varnothing$,
assume also that $\dim_{\mathcal H}S<n-2-2\alpha$.
There exist $\eps_0=\eps_0(n,\delta,\alpha)>0$ and
$C=C(n,\delta,\alpha)>0$ such that, if
\[
r^{2\alpha-n}\int_{M\cap B_{2r}(x)}|A_M|^{2\alpha}\,d\mu\le\eps_0,
\]
then
\begin{equation}\label{eq:eps-sup}
r\sup_{M\cap B_{r/2}(x)}|A_M|
\le C\Bigl(
r^{2\alpha-n}\int_{M\cap B_{2r}(x)}|A_M|^{2\alpha}\,d\mu
\Bigr)^{1/(2\alpha)}.
\end{equation}
\end{proposition}

\begin{proof}
Applying $z\mapsto(z-x)/r$, we may assume $x=0$ and $r=1$.
When $S\ne\varnothing$, we first extend the integral inequalities across
$S$, following the cutoff argument in \cite[Lemma~3.3 and the proof of
Lemma~3.2]{HongLiWang2024}.
For a stationary integral $n$-varifold $V$ on an open set $U$, the
monotonicity identity \cite[Chapter~4, \textup{(3.6)}]{Simon2014GMT} is
\begin{equation}\label{eq:stationary-monotonicity}
\begin{aligned}
\frac{\|V\|(B_t(y))}{t^n}-\frac{\|V\|(B_s(y))}{s^n}
=&\int_{B_t(y)\setminus B_s(y)}
\frac{|\pi_{(T_zV)^\perp}(z-y)|^2}{|z-y|^{n+2}}\,d\|V\|(z)\ge0,
\end{aligned}
\end{equation}
where $0<s<t$, $\overline B_t(y)\Subset U$, and $T_zV$ is the
approximate tangent plane. Applied to the stationary extension of $|F|$,
\eqref{eq:stationary-monotonicity} implies
$\mu(F^{-1}(B_\rho(y)))\le C\rho^n$ on each fixed smaller ball.
For $\eta\in\operatorname{Lip}_c(B_2)$ and each $\ell$, cover
$S\cap\spt\eta$ by finitely many balls $B_{r_i}(x_i)$, with
$\max_i r_i\to0$ and $\sum_i r_i^{n-2-2\alpha}\to0$ as $\ell\to\infty$.
Choose $\chi_\ell\in\operatorname{Lip}(M)$ equal to zero on
$F^{-1}(\bigcup_iB_{r_i}(x_i))$ and one outside
$F^{-1}(\bigcup_iB_{2r_i}(x_i))$, using the minimum of functions whose
gradients have norm at most $2/r_i$. Then
\begin{equation}\label{eq:singular-cutoff-gradient}
\int_M|\nabla\chi_\ell|^2\,d\mu
\le C\sum_i r_i^{n-2}\longrightarrow0\quad\text{as }\ell\to\infty.
\end{equation}
Substituting $(\eta\circ F)\chi_\ell$ into \eqref{eq:delta-stability}
and using \eqref{eq:P-estimates}, we obtain
\[
\frac\delta n\int_M|A_M|^2(\eta\circ F)^2\chi_\ell^2\,d\mu
\le2\int_M\chi_\ell^2|\nabla(\eta\circ F)|^2\,d\mu
 +2\|\eta\|_\infty^2\int_M|\nabla\chi_\ell|^2\,d\mu.
\]
Letting $\ell\to\infty$ by \eqref{eq:singular-cutoff-gradient} and
then choosing $\eta$ on concentric balls,
we have
\begin{equation}\label{eq:singular-curvature-L2}
\int_{F^{-1}(B_\rho(y))}|A_M|^2\,d\mu\le C\rho^{n-2}
\end{equation}
on each fixed smaller ball. By \eqref{eq:singular-curvature-L2}
and H\"older's inequality,
\begin{equation}\label{eq:singular-weighted-cutoff}
\begin{aligned}
\int_M|A_M|^{2\alpha}|\nabla\chi_\ell|^2\,d\mu
\le&C\sum_i r_i^{-2}
 \Bigl(\int_{F^{-1}(B_{2r_i}(x_i))}|A_M|^2\,d\mu\Bigr)^\alpha
 \mu(F^{-1}(B_{2r_i}(x_i)))^{1-\alpha}\\
\le&C\sum_i r_i^{n-2-2\alpha}\longrightarrow0\quad\text{as }\ell\to\infty.
\end{aligned}
\end{equation}
By the Michael--Simon inequality
\cite[Theorem~2.1]{MichaelSimon1973} and H\"older's inequality, we have
\begin{equation}\label{eq:MS}
\Bigl(\int_M|w|^{\frac{2n}{n-2}}\,d\mu\Bigr)^{\frac{n-2}n}
\le C\int_M|\nabla w|^2\,d\mu,\quad \forall\,w\in W_c^{1,2}(M),
\end{equation}
where $C=C(n)>0$.

Using $(\eta\circ F)\chi_\ell$ and letting $\ell\to\infty$, by
\eqref{eq:singular-weighted-cutoff} and \eqref{eq:zero-level-energy}
we extend \eqref{eq:level-cacc} and \eqref{eq:zero-level-energy} to
$\eta\circ F$, and \eqref{eq:MS} to $(|A_M|^\alpha-k)^+(\eta\circ F)$.
For $S=\varnothing$, $\eta\circ F$ has compact support by properness.

For $0<d\le1$, define $k_\ell\coloneqq d(1-2^{-\ell-1})$ and
\[
Y_\ell\coloneqq \int_{\Omega_\ell}(|A_M|^\alpha-k_\ell)^p\,d\mu,
\]
where $p\coloneqq 2n/(n-2)$ and
$\Omega_\ell\coloneqq M\cap B_{1/2+2^{-\ell}}\cap\{|A_M|^\alpha>k_\ell\}$.
Choose $\eta_\ell\in C_c^1(B_{1/2+2^{-\ell}})$ with
$\eta_\ell|_{B_{1/2+2^{-\ell-1}}}=1$, $0\le\eta_\ell\le1$, and
$|\nabla\eta_\ell|\le C2^\ell$.
Then
\begin{equation}\label{eq:DG-mid}
\begin{aligned}
Y_{\ell+1}^{2/p}
\le&\Bigl(\int_M
 \bigl[(|A_M|^\alpha-k_{\ell+1})^+(\eta_\ell\circ F)\bigr]^p
 \,d\mu\Bigr)^{2/p}\\
\le&C\int_M\bigl|\nabla\bigl[
 (|A_M|^\alpha-k_{\ell+1})^+(\eta_\ell\circ F)\bigr]\bigr|^2\,d\mu\\
\le&C\int_{\Omega_\ell}2^\ell
 \Bigl(1-\frac{k_\ell}{|A_M|^\alpha}\Bigr)
 |\nabla(|A_M|^\alpha)|^2(\eta_\ell\circ F)^2  + 4^\ell (|A_M|^\alpha-k_\ell)^2\,d\mu\\
\le&C8^\ell\int_{\Omega_\ell}(|A_M|^\alpha-k_\ell)^2\,d\mu
 +C2^\ell k_\ell^{2+2/\alpha}|\Omega_\ell|\\
\le&Cb^\ell\Bigl[\int_{\Omega_\ell}(|A_M|^\alpha-k_\ell)^2\,d\mu
 +d^{2+2/\alpha}|\Omega_\ell|\Bigr],
\end{aligned}
\end{equation}
where we used \eqref{eq:MS} in the second inequality,
$1-k_\ell/|A_M|^\alpha\ge2^{-\ell-2}$ on
$\{|A_M|^\alpha>k_{\ell+1}\}$ in the third, and
\eqref{eq:level-cacc} with $k=k_\ell$ in the fourth;
$C=C(n,\delta,\alpha)\ge1$ and $b=b(n,\delta,\alpha)\ge8$.
For $\ell\ge1$, on $\Omega_\ell$,
\begin{equation}\label{eq:DG-level-gap}
|A_M|^\alpha-k_{\ell-1}\ge d2^{-\ell-1},
\quad 0\le |A_M|^\alpha-k_\ell\le |A_M|^\alpha-k_{\ell-1}.
\end{equation}
By \eqref{eq:DG-mid},
\[
Y_{\ell+1}^{2/p}
\le Cb^\ell\Bigl[\Bigl(\frac{2^{\ell+1}}d\Bigr)^{p-2}
 +d^{2+2/\alpha}\Bigl(\frac{2^{\ell+1}}d\Bigr)^p\Bigr]Y_{\ell-1} 
\le Cb^\ell d^{2-p}Y_{\ell-1}.
\]
Here we used \eqref{eq:DG-level-gap} in the first inequality and
$d\le1$ in the second, increasing $C,b$ if necessary.
Consequently,
\begin{equation}\label{eq:DG-rec}
d^{-p}Y_{\ell+1}\le Cb^\ell(d^{-p}Y_{\ell-1})^{p/2}
\qquad\forall\,\ell\ge1.
\end{equation}
Choose $\eta\in C_c^1(B_2)$ with $\eta|_{B_{3/2}}=1$ and
$|\nabla\eta|\le C$.
Using \eqref{eq:MS} and \eqref{eq:level-cacc} with $k=d/4$
as in \eqref{eq:DG-mid},
\[
\begin{aligned}
Y_0^{2/p}
\le&\Bigl(\int_M\bigl[(|A_M|^\alpha-d/2)^+(\eta\circ F)\bigr]^p\,d\mu\Bigr)^{2/p}\\
\le&C\int_M\bigl|\nabla\bigl[(|A_M|^\alpha-d/2)^+(\eta\circ F)\bigr]\bigr|^2\,d\mu\\
\le&C\|A_M\|_{L^{2\alpha}(M\cap B_2)}^{2\alpha}
 +Cd^{2+2/\alpha}|\{|A_M|^\alpha>d/4\}\cap B_2|\\
\le&C(1+d^{2/\alpha})\|A_M\|_{L^{2\alpha}(M\cap B_2)}^{2\alpha}
\le C_0\|A_M\|_{L^{2\alpha}(M\cap B_2)}^{2\alpha},
\end{aligned}
\]
where we used
\[
\bigl|\{|A_M|^\alpha>d/4\}\cap B_2\bigr|
\le\frac{16}{d^2}\int_{M\cap B_2}|A_M|^{2\alpha}\,d\mu
\]
in the fourth inequality and $d\le1$ in the last inequality.
Choose $\eps_1,\eps_0>0$ so that
$C\eps_1^{p/2-1}b^{n-1}\le1$ and
$(C_0\eps_0)^{p/2}\le\eps_1 b^{-(n-2)/2}$,
where $C,b$ are as in \eqref{eq:DG-rec}.
If $\|A_M\|_{L^{2\alpha}(M\cap B_2)}^{2\alpha}\le\eps_0d^2$, then
$Y_1\le Y_0\le(C_0\eps_0)^{p/2}d^p
\le\eps_1 d^p b^{-(n-2)/2}$.
Iterating \eqref{eq:DG-rec}
over the even and odd indices, we obtain
\[
Y_j\le\eps_1 d^p b^{-(n-2)j/2}\longrightarrow0\quad\text{as }j\to\infty.
\]
Thus $\int_{M\cap B_{1/2}}(|A_M|^\alpha-d)_+^p\,d\mu=0$ and
$\sup_{M\cap B_{1/2}}|A_M|\le d^{1/\alpha}$.
Since $d\in(0,1]$ is arbitrary with
$\|A_M\|_{L^{2\alpha}(M\cap B_2)}^{2\alpha}\le\eps_0d^2$,
\[
\sup_{M\cap B_{1/2}}|A_M|
\le\eps_0^{-1/(2\alpha)}\|A_M\|_{L^{2\alpha}(M\cap B_2)}.
\]
By the change of variables $z\mapsto(z-x)/r$, we obtain
\eqref{eq:eps-sup}.
\end{proof}

\section{Stable cones}\label{sec:cones}

We first prove the weighted Hardy inequality \eqref{eq:cone-hardy-main}.

\begin{proposition}\label{prop:cone-hardy}
Let $C^m\subset\C^m$, $m\ge3$, be a special Lagrangian
cone such that $\Sigma\coloneqq C\cap\Sph^{2m-1}$ is smooth and compact.
If $C$ is $\delta$-Simons stable and $\delta\ge(m-3)/m$, then
\begin{equation}\label{eq:cone-hardy-main}
\delta(\delta+2)\int_C r^{-2}|A_C|^{2\delta}\phi^2\,d\mu_C
\le\int_C |A_C|^{2\delta}|\nabla\phi|^2\,d\mu_C
\end{equation}
for all $\phi\in\operatorname{Lip}_c(C\setminus\{0\})$, where $r$ is
the distance from the vertex.
\end{proposition}

\begin{proof}
If $\delta=0$ or $A_C\equiv0$, \eqref{eq:cone-hardy-main} is immediate.
Assume $\delta>0$. With $g_C=dr^2+r^2g_\Sigma$,
$h_\Sigma$ as in \eqref{eq:cone-tensor-definition}, and
$|A_C|(r,p)=r^{-1}|h_\Sigma|(p)$, we have
\begin{align}
|\nabla A_C|^2=&r^{-4}\bigl(|\nabla^\Sigma h_\Sigma|^2+4|h_\Sigma|^2\bigr),
\label{eq:radial-A}\\
\bigl|\nabla|A_C|\bigr|^2=&r^{-4}\bigl(|h_\Sigma|^2+
\bigl|\nabla^\Sigma|h_\Sigma|\bigr|^2\bigr).
\label{eq:radial-f}
\end{align}
Applying \cite[Lemma~5.17, (5.39)]{Fox2021} to the trace free Codazzi
tensor $h_\Sigma$ on $\Sigma$, we obtain
\begin{equation}\label{eq:cone-link-kato}
|\nabla^\Sigma h_\Sigma|^2\ge\frac{m+3}{m}
\bigl|\nabla^\Sigma|h_\Sigma|\bigr|^2.
\end{equation}
Since $2-\delta\le(m+3)/m$,
\begin{equation}\label{eq:cone-kato}
\begin{aligned}
|\nabla A_C|^2
= r^{-4}\bigl(|\nabla^\Sigma h_\Sigma|^2+4|h_\Sigma|^2\bigr) 
\ge&r^{-4}\Bigl((2-\delta)\bigl|\nabla^\Sigma|h_\Sigma|\bigr|^2
 +4|h_\Sigma|^2\Bigr)\\
=&(2-\delta)\bigl|\nabla|A_C|\bigr|^2
 +(2+\delta)\frac{|A_C|^2}{r^2},
\end{aligned}
\end{equation}
where we used \eqref{eq:radial-A}--\eqref{eq:radial-f} in the equalities
and \eqref{eq:cone-link-kato} in the inequality.

Combining \eqref{eq:exact-simons} with \eqref{eq:cone-kato},
\begin{equation}\label{eq:cone-simons-delta}
\frac12\Delta |A_C|^2+\Pcal |A_C|^2
= |\nabla A_C|^2 
\ge (2-\delta)\bigl|\nabla|A_C|\bigr|^2+(2+\delta)r^{-2}|A_C|^2. 
\end{equation}
For $\phi\in C_c^1(C\cap\{|A_C|>0\})$, multiplying
\eqref{eq:cone-simons-delta} by $|A_C|^{2\delta-2}\phi^2$ and integrating
by parts over $C$, we derive
\[
\begin{aligned}
\int_C\Pcal |A_C|^{2\delta}\phi^2
\ge&\delta\int_C |A_C|^{2\delta-2}\bigl|\nabla|A_C|\bigr|^2\phi^2 +2\int_C |A_C|^{2\delta-1}\phi\langle\nabla|A_C|,\nabla\phi\rangle +(2+\delta)\int_Cr^{-2}|A_C|^{2\delta}\phi^2.
\end{aligned}
\]
With $\psi=|A_C|^\delta\phi$ in \eqref{eq:delta-stability},
\begin{equation}\label{eq:cone-stability-test}
\begin{aligned}
&\delta(\delta+2)\int_Cr^{-2}|A_C|^{2\delta}\phi^2\,d\mu_C\\
\le&\delta\int_C\Pcal|A_C|^{2\delta}\phi^2\,d\mu_C
 -\delta^2\int_C|A_C|^{2\delta-2}\bigl|\nabla|A_C|\bigr|^2\phi^2\,d\mu_C -2\delta\int_C|A_C|^{2\delta-1}\phi
 \langle\nabla|A_C|,\nabla\phi\rangle\,d\mu_C\\
\le&\int_C|A_C|^{2\delta}|\nabla\phi|^2\,d\mu_C.
\end{aligned}
\end{equation}
Thus \eqref{eq:cone-hardy-main} holds on $\{|A_C|>0\}$.

For $\phi\in\operatorname{Lip}_c(C\setminus\{0\})$,
\eqref{eq:cone-hardy-main} follows by repeating
\eqref{eq:cone-stability-test} with the regularization in
\eqref{eq:regularized-power} and \eqref{eq:cone-kato}, letting
$\eps\downarrow0$, and approximating in $W^{1,2}$.
\end{proof}

We can now prove Theorem~\ref{thm:intro-cone}.

\begin{proof}[Proof of Theorem~\ref{thm:intro-cone}]
Suppose $A_C\not\equiv0$. In polar coordinates, for $\eta\in C_c^1(0,\infty)$,
\begin{equation}\label{eq:cone-radial-integrals}
\begin{aligned}
\int_Cr^{-2}|A_C|^{2\delta}\eta(r)^2
&=\Bigl(\int_\Sigma|h_\Sigma|^{2\delta}\,d\mu_\Sigma\Bigr)
\int_0^\infty r^{m-3-2\delta}\eta^2\,dr,\\
\int_C|A_C|^{2\delta}|\nabla(\eta\circ r)|^2
&=\Bigl(\int_\Sigma|h_\Sigma|^{2\delta}\,d\mu_\Sigma\Bigr)
\int_0^\infty r^{m-1-2\delta}(\eta')^2\,dr.
\end{aligned}
\end{equation}
Since $\delta>\delta_{\rm cone}(m)\ge(m-3)/m$, substituting
$\eta\circ r$ into \eqref{eq:cone-hardy-main} and using
\eqref{eq:cone-radial-integrals}, we obtain
\begin{equation}\label{eq:radial-hardy-comp}
\delta(\delta+2)\int_0^\infty r^{m-3-2\delta}\eta^2\,dr
\le\int_0^\infty r^{m-1-2\delta}(\eta')^2\,dr.
\end{equation}
For $\psi\in C_c^1(\R)$, take
\[
\begin{aligned}
\eta(r)\coloneqq&r^{-(m-2-2\delta)/2}\psi(\log r).
\end{aligned}
\]
Then
\begin{equation}\label{eq:radial-energy-identity}
\begin{aligned}
\int_0^\infty r^{m-3-2\delta}\eta^2\,dr
=&\int_0^\infty\psi(\log r)^2\,\frac{dr}{r}
=\int_\R\psi^2\,dt,\\
\int_0^\infty r^{m-1-2\delta}(\eta')^2\,dr
=&\int_0^\infty
 \Bigl(\psi'(\log r)-\frac{m-2-2\delta}{2}\psi(\log r)\Bigr)^2
 \,\frac{dr}{r}\\
=&\int_\R\Bigl(\psi'-\frac{m-2-2\delta}{2}\psi\Bigr)^2\,dt\\
=&\int_\R(\psi')^2\,dt
 +\frac{(m-2-2\delta)^2}{4}\int_\R\psi^2\,dt.
\end{aligned}
\end{equation}
For fixed $0\ne\zeta\in C_c^1(\R)$ and $R>0$, take
$\psi(t)\coloneqq\zeta(t/R)$, so $\psi'(t)=R^{-1}\zeta'(t/R)$.
By \eqref{eq:radial-energy-identity} and \eqref{eq:radial-hardy-comp},
with $t=Rs$,
\[
\begin{aligned}
0\le&\frac{\int_\R(\psi')^2\,dt}{\int_\R\psi^2\,dt}
+\frac{(m-2-2\delta)^2}{4}-\delta(\delta+2)\\
=&\frac{R^{-1}\int_\R(\zeta'(s))^2\,ds}
 {R\int_\R\zeta(s)^2\,ds}
+\frac{(m-2)^2}{4}-m\delta\\
=&R^{-2}\frac{\int_\R(\zeta')^2\,ds}{\int_\R\zeta^2\,ds}
+\frac{(m-2)^2}{4}-m\delta
\longrightarrow\frac{(m-2)^2}{4}-m\delta<0
\quad\text{as }R\to\infty,
\end{aligned}
\]
a contradiction since $\delta>\delta_{\rm cone}(m)$.

Since $A_C=0$, for each connected component $\Sigma_j$ of $\Sigma$
we may choose a special Lagrangian $m$-plane $P_j$ of the prescribed
phase such that
\[
\Sigma_j\longrightarrow P_j\cap\Sph^{2m-1}
\]
is a finite covering of degree $d_j\in\mathbb N$.
Only finitely many $\Sigma_j$ occur because $\Sigma$ is compact.
Thus $C=\sum_jd_j\llbracket P_j\rrbracket$ is a finite union of special
Lagrangian $m$-planes with the same phase and positive integer multiplicities.
\end{proof}

In the following example, we prove that the nonflat Clifford torus
cone is $\delta_{\rm cone}(m)$-Simons stable, so the coefficient in
Theorem~\ref{thm:intro-cone} is sharp.

\begin{example}\label{ex:clifford-cone}
Let $m\ge3$ and let
\[
\mathrm{Cliff}^{m}\coloneqq 
\Bigl\{
\frac r{\sqrt m}(e^{i\theta_1},\dots,e^{i\theta_m}):
r\ge0,\ \theta_j\in\R,\ \sum_{j=1}^m\theta_j=0
\Bigr\}
\]
be the Clifford torus cone. Then
\begin{equation}\label{eq:HL-best-delta}
|A_{\mathrm{Cliff}^{m}}|^2=\frac{(m-1)(m-2)}{r^2},\quad
\Pcal=\frac m{r^2},\quad
\delta_{\rm cone}(m)=\frac{(m-2)^2}{4m}.
\end{equation}
Consequently, for $\delta>0$, the nonflat cone $\mathrm{Cliff}^{m}$ is
$\delta$-Simons stable if and only if $\delta\le\delta_{\rm cone}(m)$.
\end{example}

\begin{proof}
The first two identities in \eqref{eq:HL-best-delta} follow from
$\Rm_\Sigma=0$ for $\Sigma\coloneqq \mathrm{Cliff}^{m}\cap\Sph^{2m-1}$
and $g_{\mathrm{Cliff}^{m}}=dr^2+r^2g_\Sigma$.
By \eqref{eq:radial-energy-identity} with $\delta=0$,
\[
\begin{aligned}
\int_{\mathrm{Cliff}^{m}}|\nabla\phi|^2
&=\int_\Sigma\int_0^\infty
\bigl(r^{m-1}|\partial_r\phi|^2+r^{m-3}|\nabla^\Sigma\phi|^2\bigr)
\,dr\,d\mu_\Sigma \ge\frac{(m-2)^2}{4m}\int_{\mathrm{Cliff}^{m}}\Pcal\phi^2.
\end{aligned}
\]
For sharpness, substituting
$\phi_R(r)\coloneqq r^{-(m-2)/2}\chi((\log r)/R)$ with fixed
$0\ne\chi\in C_c^1(\R)$ into \eqref{eq:radial-energy-identity}, we derive
\[
\begin{aligned}
\frac{\int_{\mathrm{Cliff}^{m}}|\nabla\phi_R|^2\,d\mu_{\mathrm{Cliff}^{m}}}
 {\int_{\mathrm{Cliff}^{m}}\Pcal\phi_R^2\,d\mu_{\mathrm{Cliff}^{m}}}
=&\frac{(m-2)^2}{4m}
+\frac{1}{mR^2}\frac{\int_\R(\chi')^2\,dt}{\int_\R\chi^2\,dt} 
\longrightarrow\frac{(m-2)^2}{4m}\quad\text{as }R\to\infty.
\end{aligned}
\]
\end{proof}

We now characterize stability at $\delta=\delta_{\rm cone}(n)$.

\begin{theorem}\label{thm:critical-cone-non4}
Let $n\ge3$ and let $C^n\subset\C^n$ be a special Lagrangian cone
such that $\Sigma\coloneqq C\cap\Sph^{2n-1}$ is smooth, compact, and
connected. If $n\ne4$, then $C$ is $\delta_{\rm cone}(n)$-Simons stable
if and only if $\nabla^\Sigma h_\Sigma=0$.
For $n=4$, $C$ is $1/4$-Simons stable if and only if either
$h_\Sigma\equiv0$ or
\begin{equation}\label{eq:critical-four-characterization}
|h_\Sigma|>0,\qquad
|\nabla^\Sigma h_\Sigma|^2
=\frac74\bigl|\nabla^\Sigma|h_\Sigma|\bigr|^2
\quad\text{on }\Sigma.
\end{equation}
\end{theorem}

\begin{proof}
If $A_C\equiv0$, then $C$ is $\delta_{\rm cone}(n)$-Simons stable and
$h_\Sigma=0$. We may therefore assume $A_C\not\equiv0$.
Assume that $C$ is $\delta_{\rm cone}(n)$-Simons stable.
We first prove,  for $n\ge3$,
\begin{equation}\label{eq:critical-cone-equality}
|\nabla^\Sigma h_\Sigma|^2
=\bigl(2-\delta_{\rm cone}(n)\bigr)
 \bigl|\nabla^\Sigma|h_\Sigma|\bigr|^2
\quad\text{on }\Sigma\cap\{|h_\Sigma|>0\}.
\end{equation}

In fact, with $\delta\coloneqq \delta_{\rm cone}(n)$, choose
$0\ne\chi\in C_c^1(\R)$ and for $R>1$, take
\[
\phi_R\coloneqq r^{-(n-2-2\delta)/2}\chi((\log r)/R).
\]
Using \eqref{eq:cone-kato}, applying the argument of
Proposition~\ref{prop:cone-hardy} to \eqref{eq:exact-simons} with
$\phi=\phi_R$, and letting $\eps\downarrow0$ by Fatou's lemma, we obtain
\[
\begin{aligned}
&\delta(\delta+2)\int_C r^{-2}|A_C|^{2\delta}\phi_R^2\,d\mu_C\\
&+\delta\int_{C\cap\{|A_C|>0\}}|A_C|^{2\delta-2}\phi_R^2 \cdot\Bigl(|\nabla A_C|^2-(2-\delta)\bigl|\nabla|A_C|\bigr|^2
-(2+\delta)r^{-2}|A_C|^2\Bigr)\,d\mu_C\\
\le&\int_C |A_C|^{2\delta}|\nabla\phi_R|^2\,d\mu_C.
\end{aligned}
\]
By \eqref{eq:radial-A}--\eqref{eq:radial-f},
\eqref{eq:radial-energy-identity} with $m=n$, and $(n-2)^2/4=n\delta$,
after division by $R\int_\R\chi^2\,ds$,
\[
\begin{aligned}
0&\le \delta\int_{\Sigma\cap\{|h_\Sigma|>0\}}
 |h_\Sigma|^{2\delta-2}
 \Bigl(|\nabla^\Sigma h_\Sigma|^2-(2-\delta)
       \bigl|\nabla^\Sigma|h_\Sigma|\bigr|^2\Bigr)\,d\mu_\Sigma\\
&\le \frac1{R^2}
 \frac{\int_\R(\chi')^2\,ds}{\int_\R\chi^2\,ds}
 \int_\Sigma|h_\Sigma|^{2\delta}\,d\mu_\Sigma.
\end{aligned}
\]
Thus \eqref{eq:critical-cone-equality} holds as $R\to\infty$.

For $n\ne4$, by \eqref{eq:critical-cone-equality} and
\eqref{eq:cone-link-kato} with $m=n$,
\[
0\ge \frac{(n-4)^2}{4n}\bigl|\nabla^\Sigma|h_\Sigma|\bigr|^2.
\]
Thus $\nabla^\Sigma|h_\Sigma|=\nabla^\Sigma h_\Sigma=0$ wherever $|h_\Sigma|>0$,
so $|h_\Sigma|$ is a positive constant on each component of $\{|h_\Sigma|>0\}$.
Each component is therefore closed in $\Sigma$. Since $\Sigma$ is connected
and $h_\Sigma\not\equiv0$, we have $|h_\Sigma|>0$ and
$\nabla^\Sigma h_\Sigma=0$ throughout $\Sigma$.

Conversely, if $\nabla^\Sigma h_\Sigma=0$,  $\Pcal=nr^{-2}$ by
\eqref{eq:radial-A}--\eqref{eq:radial-f} and \eqref{eq:exact-simons}.
For $\phi\in C_c^1(C\setminus\{0\})$, applying
\eqref{eq:radial-energy-identity} with $m=n$ and $\delta=0$ to
$r\mapsto\phi(r,p)$ on $(0,\infty)$ and integrating over
$p\in\Sigma$, we obtain
\[
\int_C|\nabla\phi|^2\,d\mu_C
\ge\frac{(n-2)^2}{4}\int_Cr^{-2}\phi^2\,d\mu_C
=\delta_{\rm cone}(n)\int_C\Pcal\phi^2\,d\mu_C.
\]

Let $n=4$ and assume that $C$ is $1/4$-Simons stable.
Substituting
$r^{-1}\chi((\log r)/R)\zeta(p)$, with
$0\ne\chi\in C_c^1(\R)$ and $\zeta\in C^1(\Sigma)$, into
\eqref{eq:delta-stability} and letting $R\to\infty$, we obtain, by
\eqref{eq:radial-energy-identity},
\begin{equation}\label{eq:critical-four-link-form}
\int_\Sigma\Bigl(|\nabla^\Sigma\zeta|^2
 +(1-\tfrac14\Pcal)\zeta^2\Bigr)\,d\mu_\Sigma\ge0
\qquad\forall\,\zeta\in C^1(\Sigma),
\end{equation}
where $\Pcal$ is restricted to $r=1$. By \eqref{eq:exact-simons} and
\eqref{eq:radial-A},
\begin{equation}\label{eq:critical-four-link-simons}
\frac12\Delta_\Sigma|h_\Sigma|^2
=|\nabla^\Sigma h_\Sigma|^2+(4-\Pcal)|h_\Sigma|^2.
\end{equation}
Applying \eqref{eq:critical-four-link-form} to
$(|h_\Sigma|^2+\eps^2)^{1/8}$ and integrating
\eqref{eq:critical-four-link-simons} by parts, we have
\[
\begin{aligned}
0\le&\int_\Sigma\Bigl(
 |\nabla^\Sigma(|h_\Sigma|^2+\eps^2)^{1/8}|^2
 +(1-\tfrac14\Pcal)(|h_\Sigma|^2+\eps^2)^{1/4}
 \Bigr)\,d\mu_\Sigma\\
=&\frac14\int_\Sigma(|h_\Sigma|^2+\eps^2)^{-3/4}
 \Bigl((4-\Pcal)\eps^2-|\nabla^\Sigma h_\Sigma|^2 +\frac{7|\nabla^\Sigma(|h_\Sigma|^2)|^2}
 {16(|h_\Sigma|^2+\eps^2)}\Bigr)\,d\mu_\Sigma\\
\le&\eps^2\int_\Sigma(|h_\Sigma|^2+\eps^2)^{-3/4}\,d\mu_\Sigma
\le\eps^{1/2}\operatorname{Vol}(\Sigma),
\end{aligned}
\]
where we used \eqref{eq:cone-link-kato} with $m=4$ in the last line.
By \eqref{eq:P-estimates}, $\Pcal$ is bounded on $\Sigma$. Hence the functions
$(|h_\Sigma|^2+\eps^2)^{1/8}$ are bounded in $W^{1,2}(\Sigma)$ and
converge uniformly to $|h_\Sigma|^{1/4}$ as $\eps\downarrow0$.
By weak lower semicontinuity, equality holds in
\eqref{eq:critical-four-link-form} for $\zeta=|h_\Sigma|^{1/4}$.
Substituting $|h_\Sigma|^{1/4}+t\eta$, $\eta\in C^\infty(\Sigma)$,
into \eqref{eq:critical-four-link-form} by density and differentiating
at $t=0$, we obtain
\[
(-\Delta_\Sigma+1-\tfrac14\Pcal)|h_\Sigma|^{1/4}=0
\quad\text{weakly on }\Sigma.
\]
By Harnack's inequality, $|h_\Sigma|>0$ on $\Sigma$, since
$h_\Sigma\not\equiv0$; hence \eqref{eq:critical-four-characterization}
follows from \eqref{eq:critical-cone-equality}.

Conversely, suppose $n=4$ and \eqref{eq:critical-four-characterization}
holds. By \eqref{eq:critical-four-link-simons},
\[
\begin{aligned}
\bigl(\Delta_C+\tfrac14\Pcal\bigr)
 \bigl(r^{-1}|h_\Sigma|^{1/4}\bigr)
=&r^{-3}\bigl(\Delta_\Sigma-1+\tfrac14\Pcal|_\Sigma\bigr)
 |h_\Sigma|^{1/4}\\
=&0.
\end{aligned}
\]
Integrating by parts with the positive function
$r^{-1}|h_\Sigma|^{1/4}$, we obtain
\[
\begin{aligned}
Q_{1/4,C}(\phi)
=&\int_C r^{-2}|h_\Sigma|^{1/2}
 \Bigl|\nabla\bigl(r|h_\Sigma|^{-1/4}\phi\bigr)\Bigr|^2\,d\mu_C\\
\ge&0\qquad\forall\,\phi\in C_c^1(C\setminus\{0\}).
\end{aligned}
\]
\end{proof}

For the four dimensional case of Corollary~\ref{cor:critical-graphical-cone},
we use the following contraction estimate \cite[Remark~5.19]{Fox2021}.

\begin{lemma}\label{lem:kato-equality-algebra}
Let $\widehat h$ be a trace free symmetric $(0,3)$-tensor on $\R^n$, $n\ge3$.
Then
\[
|\widehat h(X,\cdot,\cdot)|^2\le\frac n{n+2}|\widehat h|^2
\qquad\forall\,X\in\Sph^{n-1}.
\]
\end{lemma}

\begin{proof}
Fix $X\in\Sph^{n-1}$ and complete $X$ to an orthonormal basis
$e_1=X,e_2,\ldots,e_n$ of $\R^n$. Let
$\widehat h_{ijk}\coloneqq \widehat h(e_i,e_j,e_k)$, with Greek
indices ranging from $2$ to $n$. Expanding the squares and using
$\sum_i\widehat h_{iik}=0$, we obtain
\begin{equation}\label{eq:contraction-decomposition}
\begin{aligned}
&\frac n{n+2}|\widehat h|^2-|\widehat h(e_1,\cdot,\cdot)|^2\\
=&\frac{n-2}{n+1}\sum_\alpha\widehat h_{11\alpha}^2
 +\frac{2(n-1)}{n+2}\sum_{\alpha,\beta}
 \Bigl(\widehat h_{1\alpha\beta}
 +\frac{\widehat h_{111}}{n-1}\delta_{\alpha\beta}\Bigr)^2\\
&+\frac n{n+2}\sum_{\alpha,\beta,\gamma}
 \Bigl(\widehat h_{\alpha\beta\gamma}
 +\frac{\widehat h_{11\alpha}\delta_{\beta\gamma}
 +\widehat h_{11\beta}\delta_{\alpha\gamma}
 +\widehat h_{11\gamma}\delta_{\alpha\beta}}{n+1}\Bigr)^2
\ge0.
\end{aligned}
\end{equation}
\end{proof}

We now prove that a cone of finite $\delta_{\rm cone}(n)$-index is flat
when $\Sigma$ is diffeomorphic to $\Sph^{n-1}$.

\begin{corollary}\label{cor:critical-graphical-cone}
Let $C^n\subset\C^n$, $n\ge3$, be a special Lagrangian cone such that
$C\cap\Sph^{2n-1}$ is parametrized by a smooth immersion of a manifold
$\Sigma$ diffeomorphic to $\Sph^{n-1}$. If
\[
\operatorname{ind}_{\delta_{\rm cone}(n)}(C\setminus\{0\})<\infty,
\]
then $C$ is a single special Lagrangian plane.

In particular, let $n\ge3$, and let
$u\in C^\infty(\R^n\setminus\{0\})$ be homogeneous of degree two. If
$\{(x,Du(x)):x\ne0\}\cup\{0\}$ is a special Lagrangian cone of finite
$\delta_{\rm cone}(n)$-index, then $u$ is a homogeneous quadratic polynomial.
\end{corollary}

\begin{proof}
Suppose $Q_{\delta_{\rm cone}(n),C}(\phi)<0$ for some
$\phi\in C_c^1(C\setminus\{0\})$. Choose $0<a<b$ with
$\operatorname{spt}\phi\subset\{a<r<b\}$ and $\lambda>b/a$.
The functions
$\phi_j(r,p)\coloneqq\phi(\lambda^{-j}r,p)$, $j\ge1$, have disjoint
supports, so for $N\ge1$ and $0\ne(c_1,\ldots,c_N)\in\R^N$,
\[
Q_{\delta_{\rm cone}(n),C}\Bigl(\sum_{j=1}^Nc_j\phi_j\Bigr)
=Q_{\delta_{\rm cone}(n),C}(\phi)
 \sum_{j=1}^N\lambda^{j(n-2)}c_j^2<0.
\]
Thus
\[
\operatorname{ind}_{\delta_{\rm cone}(n)}(C\setminus\{0\})
\ge\dim\operatorname{span}\{\phi_1,\ldots,\phi_N\}=N
\qquad\forall\,N\ge1,
\]
contradicting the hypothesis of finite index. Hence $C$ is
$\delta_{\rm cone}(n)$-Simons stable.

Assume $A_C\not\equiv0$. Suppose first that
$\nabla^\Sigma h_\Sigma=0$. Differentiating the Gauss equation, we obtain
$\nabla^\Sigma\Rm_\Sigma=0$.
Since $\Sigma$ is compact and simply connected, $\Sigma$ is a symmetric
space by \cite[Corollary~B.6]{Ballmann2006}. Applying
\cite[Corollary, p.~60]{Wolf1969} with
$H^*(\Sigma;\mathbb Z)\cong H^*(\Sph^{n-1};\mathbb Z)$, we conclude
that $(\Sigma,g_\Sigma)$ is a round sphere. Fix $p\in\Sigma$. Since
$\nabla^\Sigma h_\Sigma=0$, parallel transport around any loop based at $p$
preserves $(h_\Sigma)_p$. The holonomy group of the round sphere is
$SO(T_p\Sigma)$ \cite[Remark~B.36]{Ballmann2006}.
Because $\dim\Sigma=n-1\ge2$, for each nonzero $v\in T_p\Sigma$ this
group contains the rotation by $\pi$ in a two dimensional subspace
containing $v$, fixing the orthogonal complement. This rotation sends
$v$ to $-v$, so
\[
(h_\Sigma)_p(v,v,v)=(h_\Sigma)_p(-v,-v,-v)=-(h_\Sigma)_p(v,v,v)=0.
\]
This equality holds for every $v\in T_p\Sigma$. Choose an orthonormal
basis $e_1,\ldots,e_{n-1}$ of $T_p\Sigma$ and let $v=\sum_i x_i e_i$.
By symmetry,
\[
\begin{aligned}
0=(h_\Sigma)_p(v,v,v)
=&\sum_i(h_\Sigma)_{iii}(p)x_i^3
 +3\sum_{i\ne j}(h_\Sigma)_{iij}(p)x_i^2x_j +6\sum_{i<j<k}(h_\Sigma)_{ijk}(p)x_ix_jx_k,
\end{aligned}
\]
where $(h_\Sigma)_{ijk}(p)=(h_\Sigma)_p(e_i,e_j,e_k)$.
Since $(x_1,\ldots,x_{n-1})\in\R^{n-1}$ is arbitrary,
$(h_\Sigma)_p=0$ at each $p\in\Sigma$, and therefore
\[
|A_C|^2(r,p)=r^{-2}|(h_\Sigma)_p|^2=0\qquad\forall\,r>0,\ \forall\,p\in\Sigma,
\]
contradicting $A_C\not\equiv0$.

By Theorem~\ref{thm:critical-cone-non4}, it remains to consider $n=4$
and $\nabla^\Sigma h_\Sigma\not\equiv0$, with
\eqref{eq:critical-four-characterization} on $\Sigma$.
We claim that there exists $\varphi\in C^\infty(\Sigma)$ satisfying
$|\nabla^\Sigma\varphi|=1$ on $\Sigma$.
Indeed, at a point where $\nabla^\Sigma|h_\Sigma|\ne0$, take
$e_1=\nabla^\Sigma|h_\Sigma|/|\nabla^\Sigma|h_\Sigma||$.
Since $\nabla^\Sigma h_\Sigma$ is fully symmetric and trace free,
\cite[(5.23) and~(5.44)]{Fox2021} with $n=k=3$ and
Lemma~\ref{lem:kato-equality-algebra} imply
\begin{equation}\label{eq:kato-equality-chain}
\begin{aligned}
|h_\Sigma|^2\bigl|\nabla^\Sigma|h_\Sigma|\bigr|^2
\le&\Bigl(\frac14|h_\Sigma|^2
 +\frac{15}{28}|h_\Sigma(e_1,\cdot,\cdot)|^2\Bigr)
 |\nabla^\Sigma h_\Sigma|^2\\
\le&\frac47|h_\Sigma|^2|\nabla^\Sigma h_\Sigma|^2
=|h_\Sigma|^2\bigl|\nabla^\Sigma|h_\Sigma|\bigr|^2,
\end{aligned}
\end{equation}
where the last equality follows from \eqref{eq:critical-four-characterization}.
Thus $|h_\Sigma(e_1,\cdot,\cdot)|^2=3|h_\Sigma|^2/5$.
In fact,
\[
\max_{|X|=1}|h_\Sigma(X,\cdot,\cdot)|^2=\frac35|h_\Sigma|^2
\quad\text{on }\Sigma.
\]
Indeed, on a connected component $U$ of the set where the inequality
is strict, $\nabla^\Sigma|h_\Sigma|=0$ by
\eqref{eq:kato-equality-chain}, and then $\nabla^\Sigma h_\Sigma=0$
by \eqref{eq:critical-four-characterization}.
Both $|h_\Sigma|$ and
$\max_{|X|=1}|h_\Sigma(X,\cdot,\cdot)|^2$ are constant on $U$,
so the strict inequality holds on $\overline U$. Hence $U$ is also
closed in $\Sigma$, contradicting connectedness and
$\nabla^\Sigma h_\Sigma\not\equiv0$.
By equality in \eqref{eq:contraction-decomposition} with $n=3$, at
each point an orthonormal frame $e_1,e_2,e_3$ can be chosen such that
\begin{equation}\label{eq:critical-four-component-form}
(h_\Sigma)_{111}=-\frac2{\sqrt{10}}|h_\Sigma|,
\qquad
(h_\Sigma)_{122}=(h_\Sigma)_{133}=\frac1{\sqrt{10}}|h_\Sigma|,
\end{equation}
with all other components determined by symmetry or equal to zero.
Here $(h_\Sigma)_{ijk}=h_\Sigma(e_i,e_j,e_k)$. In these frames,
\[
\Bigl(\sum_{a,b}(h_\Sigma)_{iab}(h_\Sigma)_{jab}\Bigr)_{i,j=1}^3
=\frac{|h_\Sigma|^2}{5}\operatorname{diag}(3,1,1).
\]
The eigenspace for $3|h_\Sigma|^2/5$ is one dimensional; choosing
$e_1$ with $h_\Sigma(e_1,e_1,e_1)<0$ defines a smooth unit vector
field on $\Sigma$.

To prove that $e_1$ is a gradient, differentiate
\eqref{eq:critical-four-component-form} in a local orthonormal frame
$e_1,e_2,e_3$. By the symmetry of $\nabla^\Sigma h_\Sigma$, for distinct
$a,b\in\{2,3\}$,
\[
e_a|h_\Sigma|
=\sqrt{10}(\nabla^\Sigma_{e_a}h_\Sigma)_{1bb}
=\sqrt{10}(\nabla^\Sigma_{e_b}h_\Sigma)_{1ab}=0.
\]
Similarly, for $a,b\in\{2,3\}$,
\[
\begin{aligned}
-4|h_\Sigma|\langle\nabla^\Sigma_{e_1}e_1,e_a\rangle
=&\sqrt{10}(\nabla^\Sigma_{e_1}h_\Sigma)_{11a} 
= \sqrt{10}(\nabla^\Sigma_{e_a}h_\Sigma)_{111}
=-2e_a|h_\Sigma|=0,\\
-4|h_\Sigma|\langle\nabla^\Sigma_{e_a}e_1,e_b\rangle
=& \sqrt{10}(\nabla^\Sigma_{e_a}h_\Sigma)_{11b} 
=\sqrt{10}(\nabla^\Sigma_{e_1}h_\Sigma)_{1ab}
=(e_1|h_\Sigma|)\delta_{ab}.
\end{aligned}
\]
It follows that
\[
\langle\nabla^\Sigma_Xe_1,Y\rangle
=\langle\nabla^\Sigma_Ye_1,X\rangle
\qquad\forall\,X,Y\in T\Sigma.
\]
Hence $e_1$ is locally the gradient of a smooth function.
Since $\Sigma\simeq\Sph^3$ is simply connected, there exists a global
$\varphi\in C^\infty(\Sigma)$ with $\nabla^\Sigma\varphi=e_1$,
proving the claim. On the compact manifold $\Sigma$, $\varphi$ attains
a maximum at some $p\in\Sigma$, where
$1=|e_1(p)|^2=|\nabla^\Sigma\varphi(p)|^2=0$, a contradiction.

The assumption $A_C\not\equiv0$ is therefore impossible both when
$\nabla^\Sigma h_\Sigma=0$ and when
$\nabla^\Sigma h_\Sigma\not\equiv0$. Hence $A_C\equiv0$.
The immersion of $\Sigma$ is then a local isometry onto
$P\cap\Sph^{2n-1}$ for a special Lagrangian $n$-plane $P$.
The local isometry $\Sigma\to P\cap\Sph^{2n-1}$ is proper because
$\Sigma$ is compact, and hence is a finite covering. Since
$P\cap\Sph^{2n-1}\simeq\Sph^{n-1}$ is simply connected for $n\ge3$,
$\Sigma$ maps diffeomorphically onto $P\cap\Sph^{2n-1}$, and $C$ is
the single special Lagrangian plane $P$.

For $C=\{(x,Du(x)):x\in\R^n\setminus\{0\}\}\cup\{0\}$,
$Du(tx)=tDu(x)$ for $t>0$. The projection
$(x,Du(x))\mapsto x/|x|$ from $\Sigma$ to $\Sph^{n-1}$ has smooth inverse
\[
\omega\longmapsto
\frac{(\omega,Du(\omega))}{\sqrt{1+|Du(\omega)|^2}},
\qquad \omega\in\Sph^{n-1}.
\]
Since $\Sigma$ is diffeomorphic to $\Sph^{n-1}$, we may apply the first
part of the corollary to $C$, obtaining $A_C=0$. Thus $Du$ is linear
and $u(x)=\frac12\langle x,Du(x)\rangle$ is a homogeneous quadratic polynomial.
\end{proof}

\section{Immersed compactness}\label{sec:immersed-compactness}

\begin{definition}\label{def:imm-reg}
Let $T$ be an integral special Lagrangian current and $x\in\spt T$.
We say that $x\in\ImmReg T$ if there exists $r>0$ such that
\begin{equation}\label{eq:immersed-sheet-representation}
T\llcorner B_r(x)=\sum_{\alpha=1}^Nq_\alpha\llbracket L_\alpha\rrbracket,
\end{equation}
where $q_\alpha\in\mathbb N$ and each $L_\alpha$ is a smooth proper
special Lagrangian immersion without boundary over $B_r(x)$, calibrated
by the same form as $T$.
Define
\[
r_{\rm reg}^T(x)\coloneqq \sup\Bigl\{r\in(0,1]:
r\max_{1\le\alpha\le N}\|A_{L_\alpha}\|_{L^\infty(L_\alpha)}\le1\Bigr\},
\]
where the supremum is over decompositions
\eqref{eq:immersed-sheet-representation}, with $\sup\varnothing=0$.
\end{definition}

We first extend immersions with locally bounded $|A_M|$ across $S$
when $\dim_{\mathcal H}S<n-2$.

\begin{lemma}\label{lem:bounded-curvature-removal}
Let $U\subset\C^n$, $n\ge3$, be open, let $S\subset U$ be relatively
closed with $\dim_{\mathcal H}S<n-2$, and let $F:M^n\to U\setminus S$
be a smooth proper special Lagrangian immersion without boundary,
with one phase on all components. Assume that the current of $F$ extends to
an integral current $T$ with $\partial T=0$ in $U$ and that
\begin{equation}\label{eq:removal-curvature-assumption}
\sup_{F^{-1}(K\setminus S)}|A_M|<\infty\qquad\forall\,K\Subset U.
\end{equation}
Then $F$ extends across $S$ to a smooth proper special Lagrangian
immersion representing $T$. Moreover, if $F$ is $\delta$-Simons stable, the
extension of $F$ across $S$ is also $\delta$-Simons stable.
\end{lemma}

\begin{proof}
Since $\dim_{\mathcal H}S<n-2$, $\|T\|(S)=0$; hence $T$ is calibrated
on $U$ and, because $\partial T=0$, $|T|$ is stationary in $U$.
Fix $x_0\in S\cap\spt T$ and
$\overline B_{4\rho}(x_0)\Subset U$. By
\eqref{eq:stationary-monotonicity},
\begin{equation}\label{eq:removal-local-mass}
\mu(F^{-1}(B_t(z)))\le Ct^n
\quad\forall\,z\in\overline B_{2\rho}(x_0),\ \forall\,t\in(0,\rho).
\end{equation}
For each $\ell\ge1$, cover $S\cap\overline B_{2\rho}(x_0)$ by finitely
many balls $B_{r_i}(z_i)$, with $z_i\in S\cap\overline B_{2\rho}(x_0)$,
$\max_i r_i<\rho/(10\ell)$ and
$\sum_i r_i^{n-2}<\ell^{-1}$. Using the construction preceding
\eqref{eq:singular-cutoff-gradient}, choose $0\le\chi_\ell\le1$ on
$U$, vanishing near $S\cap\overline B_{2\rho}(x_0)$ and converging
to one on $U\setminus S$ as $\ell\to\infty$.
By \eqref{eq:removal-local-mass},
\begin{equation}\label{eq:removal-cutoff-estimates}
\int_M|\nabla(\chi_\ell\circ F)|^q\,d\mu
\le C\sum_i r_i^{n-q}\longrightarrow0
\quad\text{as }\ell\to\infty,\qquad q=1,2.
\end{equation}

In Euclidean coordinates, $F=(F_1,\ldots,F_{2n})$. By
\eqref{eq:removal-curvature-assumption} with
$K=\overline B_{4\rho}(x_0)$,
\begin{equation}\label{eq:removal-projection-derivative}
|\nabla\pi_{dF(TM)}|^2=2|A_M|^2\le C
\quad\text{on }F^{-1}(B_{4\rho}(x_0)),
\end{equation}
where $\pi_{dF(TM)}$ is orthogonal projection onto $dF(TM)$.
For $\psi\in C_c^1(B_\rho(x_0)\times\mathbb R^{(2n)^2})$ and
$i=1,\ldots,2n$, the vector field
$(\chi_\ell\circ F)\psi(F,\pi_{dF(TM)})\nabla F_i$ has compact support
because $F$ is proper and $\chi_\ell=0$ near
$S\cap\overline B_\rho(x_0)$. Integrating by parts with
$(\chi_\ell\circ F)\psi(F,\pi_{dF(TM)})\nabla F_i$ over $M$ and using
$\Delta F_i=0$, we obtain
\begin{equation}\label{eq:removal-projection-identity}
\begin{aligned}
&\int_M(\chi_\ell\circ F)
 \Bigl(\sum_{j=1}^{2n}\langle\nabla F_i,\nabla F_j\rangle
 \partial_{x_j}\psi(F,\pi_{dF(TM)}) 
 +D_2\psi(F,\pi_{dF(TM)})
 [\nabla_{\nabla F_i}\pi_{dF(TM)}]\Bigr)\,d\mu\\
=&-\int_M\psi(F,\pi_{dF(TM)})
 \langle\nabla(\chi_\ell\circ F),\nabla F_i\rangle\,d\mu,
\end{aligned}
\end{equation}
where $D_2$ denotes differentiation in the second variable.
By $|\nabla F_i|\le1$ and \eqref{eq:removal-cutoff-estimates},
\[
\begin{aligned}
&\Bigl|\int_M\psi(F,\pi_{dF(TM)})
 \langle\nabla(\chi_\ell\circ F),\nabla F_i\rangle\,d\mu\Bigr| 
\le \|\psi\|_\infty\int_M|\nabla(\chi_\ell\circ F)|\,d\mu
\longrightarrow0\quad\text{as }\ell\to\infty.
\end{aligned}
\]
Letting $\ell\to\infty$ in \eqref{eq:removal-projection-identity}, we
verify that $|T|\llcorner B_\rho(x_0)$ is a curvature varifold with
$|\nabla\pi_{dF(TM)}|\le C$ on $B_\rho(x_0)$, in the sense of
\cite[Definition~4.1 and Remark~4.2]{Aiex2024}.

Under $z\mapsto(z-x_0)/r$, for $R>0$,
$\psi\in C_c^1(B_R(0)\times\R^{(2n)^2})$, $i=1,\ldots,2n$, and
$0<r<\rho/R$, by $|\nabla F_i|\le1$,
\eqref{eq:removal-projection-derivative}, and \eqref{eq:removal-local-mass},
\begin{equation}\label{eq:removal-scaled-curvature-term}
\begin{aligned}
&\Bigl|r^{1-n}\int_M
D_2\psi\Bigl(\frac{F-x_0}{r},\pi_{dF(TM)}\Bigr)
 [\nabla_{\nabla F_i}\pi_{dF(TM)}]\,d\mu\Bigr|\\
\le&C r^{1-n}\|D_2\psi\|_\infty
 \mu\bigl(F^{-1}(B_{Rr}(x_0))\bigr)
\le C R^n r\|D_2\psi\|_\infty\longrightarrow0
 \quad\text{as }r\downarrow0.
\end{aligned}
\end{equation}
Rescaling \eqref{eq:removal-projection-identity} after $\ell\to\infty$
and then letting $r\downarrow0$ along a sequence with tangent cone $W$,
we obtain, by \eqref{eq:removal-scaled-curvature-term},
\[
\int_{B_R(0)}\sum_{j=1}^{2n}
 (\pi_{T_zW})_{ij}\,\partial_{z_j}\psi(z,\pi_{T_zW})\,d\|W\|(z)=0.
\]
Hence $W$ is a curvature varifold with zero generalized curvature
in the sense of \cite[Definition~4.1]{Aiex2024}. Applying
\cite[Theorem~4.3]{Aiex2024}, we conclude that $W$ is locally a finite
sum of affine planes with positive integer multiplicities. Since $W$ is invariant
under positive dilations about zero,
$W=\sum_{\alpha=1}^Nq_\alpha|P_\alpha|$, where $P_\alpha$ are
$n$-planes through zero and $q_\alpha\in\mathbb N$.
Using $\mathcal H^n(P_\alpha\cap B_1)=\omega_n$, we obtain
\begin{equation}\label{eq:removal-integer-density}
\Theta^n(\|T\|,x_0)
=\frac{\|W\|(B_1)}{\omega_n}
=\sum_{\alpha=1}^Nq_\alpha\in\mathbb N.
\end{equation}
By \eqref{eq:removal-curvature-assumption} and
\eqref{eq:removal-local-mass}, for any fixed $p>n$,
\begin{equation}\label{eq:removal-density-curvature}
r^{1-n/p}\Bigl(\int_{F^{-1}(B_r(x_0))}|A_M|^p\,d\mu\Bigr)^{1/p}
\le Cr\longrightarrow0\quad\text{as }r\downarrow0.
\end{equation}
Applying \cite[Theorem~4.11 and Remark~4.12]{Aiex2024} with
\eqref{eq:removal-projection-derivative}, \eqref{eq:removal-integer-density},
and \eqref{eq:removal-density-curvature}, we conclude that $|T|$ is a finite
Lipschitz multiple graph over one $n$-plane $P_0$ in a neighborhood of $x_0$.

Choose $r_\perp>0$ and a ball $D\subset x_0+P_0$ so that
$\mathcal Z=D+B_{r_\perp}^{P_0^\perp}(0)$ has closure contained in the neighborhood where
$|T|$ is a multiple graph over $P_0$.
Shrink $D$ so that every point of $\spt T\cap\mathcal Z$ has normal
distance less than $r_\perp/2$ from $x_0+P_0$. Let $\pi$ be orthogonal
projection onto $x_0+P_0$. The set $D\cap\pi(S\cap\overline{\mathcal Z})$
is relatively closed in $D$ and has Hausdorff dimension less than $n-2$.
Since $|T|$ is a Lipschitz multiple graph over $P_0$ and 
$dF$ is continuous, 
\[
\bigl\|\bigl(\pi_{P_0}|_{dF(T_xM)}\bigr)^{-1}\bigr\|\le C
\qquad\forall\,x\in F^{-1}(\mathcal Z).
\]
Hence $\pi\circ F$, restricted to
$F^{-1}(\mathcal Z\cap\pi^{-1}(D\setminus\pi(S\cap\overline{\mathcal Z})))$, is a local
diffeomorphism to $D\setminus\pi(S\cap\overline{\mathcal Z})$.
For $K\Subset D\setminus\pi(S\cap\overline{\mathcal Z})$, $F(M)\cap\mathcal Z\cap\pi^{-1}(K)$
lies in the compact set
$K+\overline B_{r_\perp/2}^{P_0^\perp}\Subset U\setminus S$.
Hence $\pi\circ F$ is proper by properness of $F$.

For a smooth loop $\gamma:\Sph^1\to D\setminus\pi(S\cap\overline{\mathcal Z})$, choose a smooth
$\Gamma:\overline{\mathbb D}\to D$ with $\Gamma|_{\Sph^1}=\gamma$, where
$\mathbb D\subset\R^2$ is the unit disk. Since
\[
\begin{aligned}
&\dim_{\mathcal H}\{e-\Gamma(z):e\in D\cap\pi(S\cap\overline{\mathcal Z}),\ z\in\overline{\mathbb D}\} 
\le \dim_{\mathcal H}\bigl(D\cap\pi(S\cap\overline{\mathcal Z})\bigr)+2<n,
\end{aligned}
\]
there exist arbitrarily small $v\in P_0$ such that
$(\Gamma(\overline{\mathbb D})+v)\cap D\cap\pi(S\cap\overline{\mathcal Z})=\varnothing$.
For sufficiently small $v$, $\Gamma(z)+tv\in D\setminus\pi(S\cap\overline{\mathcal Z})$ for
$z\in\partial\mathbb D$ and $0\le t\le1$.
For $x,y\in D\setminus\pi(S\cap\overline{\mathcal Z})$, translating
the line segment from $x$ to $y$, we obtain paths in
$D\setminus\pi(S\cap\overline{\mathcal Z})$ from $x$ to $y$ with lengths approaching $|x-y|$.
Thus $D\setminus\pi(S\cap\overline{\mathcal Z})$ is connected and simply connected.
Since $\pi\circ F$ is a proper local diffeomorphism, each component of
$F^{-1}(\mathcal Z\cap\pi^{-1}(D\setminus\pi(S\cap\overline{\mathcal Z})))$ maps diffeomorphically
onto $D\setminus\pi(S\cap\overline{\mathcal Z})$, and there are finitely many components.
Hence $F$ over $\mathcal Z\cap\pi^{-1}(D\setminus\pi(S\cap\overline{\mathcal Z}))$ consists of smooth graphs
$u_a:D\setminus\pi(S\cap\overline{\mathcal Z})\to P_0^\perp$.
By $|Du_a|\le C$ and \eqref{eq:removal-curvature-assumption},
$|D^2u_a|\le C$.
For $x,y\in D\setminus\pi(S\cap\overline{\mathcal Z})$, integrating
$Du_a$ and $D^2u_a$ along paths from $x$ to $y$ with lengths
approaching $|x-y|$, we obtain
\[
|u_a(x)-u_a(y)|\le C|x-y|,\qquad
|Du_a(x)-Du_a(y)|\le C|x-y|.
\]
Each $u_a$ extends to $C^{1,1}(D)$. Since
$\dim_{\mathcal H}\bigl(D\cap\pi(S\cap\overline{\mathcal Z})\bigr)<n-2$,
$\mathcal H^n\bigl(D\cap\pi(S\cap\overline{\mathcal Z})\bigr)=0$; hence the minimal graph
equations on $D\setminus\pi(S\cap\overline{\mathcal Z})$ hold almost everywhere on $D$:
\begin{equation}\label{eq:minimal-graph-equations}
\sum_{i,j=1}^n\bigl(I+(Du_a)^{\mathsf T}Du_a\bigr)^{-1}_{ij}
\partial_{ij}u_a=0\qquad\text{almost everywhere on }D.
\end{equation}
Because $u_a\in C^{1,1}(D)$ and $|Du_a|\le C$, the matrix
$(I+(Du_a)^{\mathsf T}Du_a)^{-1}$ is uniformly elliptic and H\"older
continuous. Schauder estimates for \eqref{eq:minimal-graph-equations}
and bootstrapping therefore imply $u_a\in C^\infty(D)$.
By continuity, $\operatorname{graph}u_a$ over $D$ is special Lagrangian
with the phase of $F$.
Over $D\cap\pi(S\cap\overline{\mathcal Z})$, $\operatorname{graph}u_a$ agrees with the corresponding sheet
of $F$ outside $S$, since $\pi\circ F$ is a local diffeomorphism on
$F^{-1}(\mathcal Z)$ and $D\setminus\pi(S\cap\overline{\mathcal Z})$ is dense in $D$.

To construct a smooth special Lagrangian extension of $F$, we first
extend $F$ continuously to the metric completion of each connected
component of $M$, using $|F(p)-F(q)|\le d_g(p,q)$.
Denote by $\widehat F:\widehat M\to U$ the disjoint union of the
continuous extensions restricted to the inverse image of $U$;
thus $\widehat M\supset M$ and $\widehat F|_M=F$.
By properness of $F:M\to U\setminus S$, a Cauchy sequence $(p_j)$
with $F(p_j)\to y\in U\setminus S$ as $j\to\infty$ converges in $M$.
Thus $\widehat F^{-1}(U\setminus S)=M$.

Fix $p\in\widehat M\setminus M$ and choose $p_j\in M$ converging to
$p$. Apply the local graph construction at $\widehat F(p)\in S$,
and choose a cylinder $\mathcal Z'\Subset\mathcal Z$ containing
$\widehat F(p)$. For all large $j$, $F(p_j)\in\mathcal Z'$.
If $p_j,p_k$ lie in different components of $F^{-1}(\mathcal Z)$,
every curve in $M$ from $p_j$ to $p_k$ leaves $F^{-1}(\mathcal Z)$,
so
\[
d_g(p_j,p_k)\ge\dist(\overline{\mathcal Z'},\partial\mathcal Z)>0.
\]
Since $(p_j)$ is Cauchy, all $p_j$ with $j\ge j_0$ lie in one
component of $F^{-1}(\mathcal Z)$, where $F$ is represented by
$x\mapsto x+u_a(x)$. By $|Du_a|\le C$ and the paths in
$D\setminus\pi(S\cap\overline{\mathcal Z})$ with lengths approaching
the Euclidean distance,
\[
|\pi F(p_j)-\pi F(p_k)|\le d_g(p_j,p_k)
\le C|\pi F(p_j)-\pi F(p_k)|\qquad\forall\,j,k\ge j_0.
\]
Hence $\pi\circ F$ extends to a coordinate map near $p$, and in
the projection coordinates $\widehat F(x)=x+u_a(x)$.
On overlaps, the changes of projection coordinates are compositions
of the smooth graph parametrizations with orthogonal projections.
Therefore $\widehat M$ is smooth and $\widehat F$ is a smooth special
Lagrangian immersion extending $F$.

To prove that $\widehat F$ is proper, fix $K\Subset U$.
Cover $K\cap\spt T$ by finitely many smaller cylinders, each with
finitely many extended graphs. Then $\widehat F^{-1}(K)$ is closed
in a finite union of compact subsets of the extended graph domains,
and hence is compact.

We also verify that $\widehat F$ represents $T$.
Since $\widehat F|_M=F$ and $\widehat F^{-1}(U\setminus S)=M$,
\[
\bigl(\widehat F_\#\llbracket\widehat M\rrbracket\bigr)
\llcorner(U\setminus S)=T\llcorner(U\setminus S).
\]
Moreover, $\int_{\widehat F^{-1}(S)}d\mu=0$ and $\|T\|(S)=0$, so $\widehat F_\#\llbracket\widehat M\rrbracket=T\quad\text{in }U$. 

Finally, suppose that $F$ is $\delta$-Simons stable. Fix
$\phi\in C_c^1(\widehat M)$ and a compact $K\Subset U$ containing
$\widehat F(\spt\phi)$ in its interior. Construct $\chi_\ell$ from a
finite covering of $S\cap K$ as above, taking $\chi_\ell=1$ if
$S\cap K=\varnothing$. Then
$(\chi_\ell\circ F)\phi\in\operatorname{Lip}_c(M)$.
By \eqref{eq:lipschitz-stability} and
\eqref{eq:removal-cutoff-estimates} with $q=2$,
\[
\begin{aligned}
\delta\int_{\widehat M}\Pcal\phi^2\,d\mu
=&\lim_{\ell\to\infty}
 \delta\int_M\Pcal(\chi_\ell\circ F)^2\phi^2\,d\mu\\
\le&\lim_{\ell\to\infty}
 \int_M\bigl|\nabla\bigl((\chi_\ell\circ F)\phi\bigr)\bigr|^2\,d\mu
=\int_{\widehat M}|\nabla\phi|^2\,d\mu.
\end{aligned}
\]
\end{proof}

\begin{proposition}\label{thm:flat-union}
Let $n\ge3$, $\delta>\delta_{\rm K}(n)$, and
$\sigma<n-4+6/(n+1)$. Suppose $S_j\subset B_4$ are relatively closed with
$\dim_{\mathcal H}S_j\le\sigma$, and let
$F_j:M_j^n\to B_4\setminus S_j$ be smooth proper
$\delta$-Simons stable special Lagrangian immersions without boundary,
with phase $\Theta_j$ on all components. Assume that the induced currents
extend to integral currents $T_j$ with $\partial T_j=0$ in $B_4$ and
$\|T_j\|(B_4)\le\Lambda<\infty$. Then the following hold: 
\begin{enumerate}[label=\textup{(\roman*)}]
\item If $|T_j|\to\sum_{\alpha=1}^Nq_\alpha|P_\alpha|$ locally as
varifolds in $B_4$ as $j\to\infty$, where $P_1,\ldots,P_N$ are distinct affine special
Lagrangian planes with the same phase, then, on each compact subset of
$B_4$ and for all large $j$, $F_j$ extends across $S_j$ and converges
smoothly as immersions.
\item There exists a subsequence such that $\Theta_j\to\Theta$ as
$j\to\infty$ and
\[
T_j\rightharpoonup T\ \text{locally as currents},\quad
|T_j|\to|T|\ \text{locally as varifolds as }j\to\infty.
\]
Here $T$ is a special Lagrangian integral current of phase $\Theta$.
Moreover, along the same subsequence, the extensions of $F_j$ across
the removable singularities in $S_j$ converge smoothly as immersions
on $\ImmReg T$ to a proper $\delta$-Simons stable special Lagrangian
immersion $G:N\to B_4^{\mathrm{reg}}$ without boundary in
$B_4^{\mathrm{reg}}$, where $B_4^{\mathrm{reg}}\coloneqq B_4\setminus\ImmSing T$
and
\[
G_\#\llbracket N\rrbracket=T\llcorner B_4^{\mathrm{reg}}.
\]
\end{enumerate}
\end{proposition}

\begin{proof}
Since $\dim_{\mathcal H}S_j<n$, $\|T_j\|(S_j)=0$ and
$\Re(e^{-i\Theta_j}\Omega)$ calibrates $T_j$ on all of $B_4$.
Because $\partial T_j=0$, $|T_j|$ is stationary in $B_4$.
We first prove \textup{(i)}.
Since $\sigma<n-4+6/(n+1)=n-2-2\delta_{\rm K}(n)$, choose
$\delta_{\rm K}(n)<\alpha<\min\{\delta,1,(n-2-\sigma)/2\}$.
Then $\dim_{\mathcal H}S_j\le\sigma<n-2-2\alpha$.
First suppose $|T_j|\to q|P|$ as $j\to\infty$ for a special Lagrangian plane $P$ and
$q\in\mathbb N$. We claim that
\begin{equation}\label{eq:plane-curv-zero}
\sup_{F_j^{-1}(B_1)}|A_{M_j}|\longrightarrow0\quad\text{as }j\to\infty.
\end{equation}

Indeed, choose $\eta\in C_c^1(B_4)$ with $\eta|_{B_3}=1$ and
$|\nabla\eta|\le C$. By \eqref{eq:P-estimates} and the cutoff argument
for \eqref{eq:singular-curvature-L2}, with $F=F_j$ and $S=S_j$,
\begin{equation}\label{eq:curvature-L2-bound}
\begin{aligned}
\frac\delta n\int_{F_j^{-1}(B_3)}|A_{M_j}|^2\,d\mu_j
\le&\delta\int_{M_j}\Pcal_j(\eta\circ F_j)^2\,d\mu_j 
\le \int_{M_j}|\nabla(\eta\circ F_j)|^2\,d\mu_j
\le C^2\Lambda.
\end{aligned}
\end{equation}
Passing to a subsequence, we may assume that
$(F_j)_\#(|A_{M_j}|^2d\mu_j)\rightharpoonup\nu$ in $B_3$ as $j\to\infty$. For $\zeta\in C_c^1(B_3)$, by \eqref{eq:delta-stability},
\eqref{eq:P-estimates}, and $|T_j|\to q|P|$ as $j\to\infty$,
\begin{equation}\label{eq:plane-curvature-measure}
\begin{aligned}
\frac{\delta}{n}\int_{B_3}\zeta^2\,d\nu
=&\lim_{j\to\infty}\frac{\delta}{n}
 \int_{M_j}|A_{M_j}|^2(\zeta\circ F_j)^2\,d\mu_j\\
\le&\lim_{j\to\infty}\int_{M_j}
 |\nabla(\zeta\circ F_j)|^2\,d\mu_j
=q\int_P|\nabla_P\zeta|^2\,d\mathcal H^n.
\end{aligned}
\end{equation}
If $F_j(p_j)\to y\in B_3$ as $j\to\infty$, then
\eqref{eq:stationary-monotonicity} and the integer multiplicity imply
\begin{equation}\label{eq:stationary-support-mass}
\omega_n\rho^n
\le\|T_j\|(B_\rho(F_j(p_j)))
\le\|T_j\|(B_{2\rho}(y))
\qquad\forall\,j\ge j_0(\rho),
\end{equation}
for $\overline B_{2\rho}(y)\Subset B_3$.
If $y\notin P$, choose $2\rho<\dist(y,P)$. By
\eqref{eq:stationary-support-mass} and $|T_j|\to q|P|$,
\[
0<\omega_n\rho^n\le\|T_j\|(B_{2\rho}(y))
\longrightarrow q\mathcal H^n(P\cap B_{2\rho}(y))=0
\quad\text{as }j\to\infty,
\]
a contradiction. Thus compact subsets of
$B_3\setminus P$ are disjoint from $F_j(M_j)$ for all large $j$, and
$\operatorname{spt}\nu\subset P$.
Moreover, \eqref{eq:stationary-monotonicity} and
$\mu_j(F_j^{-1}(B_4))\le\Lambda$ imply
$\mu_j(F_j^{-1}(B_r(x)))\le C\Lambda r^n$ for $x\in B_2$ and $0<r<1$.
By H\"older's inequality,
\begin{equation}\label{eq:L2-to-alpha}
r^{2\alpha-n}\int_{F_j^{-1}(B_r(x))}|A_{M_j}|^{2\alpha}
\le C\Bigl(r^{2-n}\int_{F_j^{-1}(B_r(x))}|A_{M_j}|^2\Bigr)^\alpha.
\end{equation}
For $\eps_0$ in Proposition~\ref{thm:eps-reg}, choose $\eta_*>0$ with
$C\eta_*^\alpha<\eps_0$.
Set
\[
E_\nu\coloneqq \Big\{x\in P\cap B_2:
\limsup_{r\downarrow0}r^{2-n}\nu(\overline B_r(x))
\ge\frac{\eta_*}{2^n}\Big\}.
\]
For $x\in E_\nu$, there exist arbitrarily small radii $\rho$ with
$\nu(\partial B_\rho(x))=0$ and
$\nu(B_\rho(x))\ge2^{-n-1}\eta_*\rho^{n-2}$.
By Vitali's $5r$-covering lemma, for $r_0>0$ we choose mutually
disjoint balls $B_{\rho_i}(x_i)$ with
$E_\nu\subset\bigcup_iB_{5\rho_i}(x_i)$ and $\rho_i<r_0$; using
$\{B_{5\rho_i}(x_i)\}_i$ in the definition of Hausdorff measure and
letting $r_0\downarrow0$, we obtain
\begin{equation}\label{eq:Hn2S}
\mathcal H^{n-2}(E_\nu)
\le C\sum_i\rho_i^{n-2}
\le C\eta_*^{-1}\nu(B_{5/2})<\infty,
\end{equation}
where $C=C(n)>0$.

Now fix $x\in(P\cap B_{3/2})\setminus E_\nu$ and choose a small ball $B_{2r}(x)$
over which the immersions $F_j$ have no boundary and
$r^{2-n}\nu(\overline B_{2r}(x))<\eta_*/2$.
For large $j$,
$\int_{F_j^{-1}(B_{2r}(x))}|A_{M_j}|^2<\eta_*r^{n-2}$, so we apply
\eqref{eq:L2-to-alpha} and Proposition~\ref{thm:eps-reg} to obtain
$\sup_{F_j^{-1}(B_{r/2}(x))}|A_{M_j}|\le C/r$.
By Lemma~\ref{lem:bounded-curvature-removal}, $F_j$ extends across
$S_j$ as a smooth proper special Lagrangian immersion over $B_{r/2}(x)$.
Using $\mu_j(F_j^{-1}(B_4))\le\Lambda$, we apply
\cite[Corollary~7.13]{Breuning2015} to the proper extensions of $F_j$
and then Schauder estimates to the graph functions satisfying
\eqref{eq:minimal-graph-equations}. After passing to a subsequence,
the local sheets converge in $C^\infty$ on $B_{r/4}(x)$. Every local smooth limit of $F_j$
on $B_{r/4}(x)$ has image in $P$, so
\[
\sup_{F_j^{-1}(B_{r/4}(x))}|A_{M_j}|\longrightarrow0
\quad\text{as }j\to\infty.
\]
Since $\mu_j(F_j^{-1}(B_4))\le\Lambda$, we have
$\nu(B_{r/4}(x))=0$. Therefore,
\[
\operatorname{spt}\nu\cap B_{3/2}\subset E_\nu.
\]
By \eqref{eq:Hn2S} and \cite[Theorem~2.27]{HKM2006}, choose
functions $\eta_\ell\in C_c^1(P\cap B_{3/2})$ with $\eta_\ell=1$ near
$\operatorname{spt}\nu\cap\overline B_{4/3}$ and
$\int_P|\nabla\eta_\ell|^2\to0$ as $\ell\to\infty$. Extending each $\eta_\ell$ to $C_c^1(B_3)$, constant in normal
directions near $P$, and substituting $\zeta=\eta_\ell$ into
\eqref{eq:plane-curvature-measure}, we derive
\[
\frac{\delta}{n}\nu(B_{4/3})
\le q\int_P|\nabla\eta_\ell|^2\longrightarrow0\quad\text{as }\ell\to\infty.
\]
It follows that $\nu(B_{4/3})=0$. As $j\to\infty$, weak convergence
implies $\int_{F_j^{-1}(B_{5/4})}|A_{M_j}|^2\to0$; hence
$\int_{F_j^{-1}(B_{5/4})}|A_{M_j}|^{2\alpha}\to0$ by H\"older's inequality
and $\mu_j(F_j^{-1}(B_4))\le\Lambda$.
Proposition~\ref{thm:eps-reg}, applied with $r=1/8$ at $x\in B_1$
on $B_{1/4}(x)\subset B_{5/4}$, now implies
\[
\sup_{F_j^{-1}(B_1)}|A_{M_j}|
\le C\Bigl(\int_{F_j^{-1}(B_{5/4})}|A_{M_j}|^{2\alpha}\Bigr)^{1/(2\alpha)}
\longrightarrow0\quad\text{as }j\to\infty.
\]

Suppose now that $|T_j|\to\sum_{\alpha=1}^Nq_\alpha|P_\alpha|$
in the varifold sense as $j\to\infty$. By \eqref{eq:curvature-L2-bound}, we may pass to a subsequence
such that $(F_j)_\#(|A_{M_j}|^2d\mu_j)\rightharpoonup\nu$ as $j\to\infty$. For $\zeta\in C_c^1(B_3)$, as in \eqref{eq:plane-curvature-measure},
\begin{equation}\label{eq:union-curvature-measure}
\frac{\delta}{n}\int_{B_3}\zeta^2\,d\nu
\le\sum_\alpha q_\alpha\int_{P_\alpha}
 |\nabla_{P_\alpha}\zeta|^2\,d\mathcal H^n.
\end{equation}
On $(\bigcup_\alpha P_\alpha)\setminus\bigcup_{\alpha<\beta}(P_\alpha\cap P_\beta)$,
we argue as in the proof of \eqref{eq:plane-curv-zero}. If $F_j(p_j)$ stays in a compact subset of
$B_3\setminus\bigcup_\alpha P_\alpha$, \eqref{eq:stationary-support-mass}
contradicts $|T_j|\to\sum_\alpha q_\alpha|P_\alpha|$ as $j\to\infty$.
Hence $A_{M_j}\to0$ locally on
$B_3\setminus\bigcup_{\alpha<\beta}(P_\alpha\cap P_\beta)$ as $j\to\infty$, and
$\operatorname{spt}\nu\cap B_{3/2}\subset\bigcup_{\alpha<\beta}(P_\alpha\cap P_\beta)$.

Assume for contradiction that two distinct special Lagrangian planes of
the same phase intersect in dimension $n-1$. Translating a point of
$P_1\cap P_2$ to the origin, choose orthonormal bases compatible with
the orientations induced by the common calibration so that
$P_1=\operatorname{span}\{e_1,\ldots,e_n\}$ and
$P_2=\operatorname{span}\{e_1,\ldots,e_{n-1},
\cos\gamma\,e_n+\sin\gamma\,Je_n\}$.
Since $dz_1\wedge\cdots\wedge dz_n$ is complex linear,
\[
\begin{aligned}
&(dz_1\wedge\cdots\wedge dz_n)
(e_1,\ldots,e_{n-1},\cos\gamma\,e_n+\sin\gamma\,Je_n)\\
=&e^{i\gamma}(dz_1\wedge\cdots\wedge dz_n)(e_1,\ldots,e_n).
\end{aligned}
\]
$P_1$ and $P_2$ have the same phase, so $e^{i\gamma}=1$ and hence
$P_1=P_2$, contradicting $P_1\ne P_2$. Thus
$\dim(P_\alpha\cap P_\beta)\le n-2$ for $\alpha\ne\beta$, and
\begin{equation}\label{eq:plane-intersection-neighborhood-volume}
\mathcal H^n\bigl(P_\alpha\cap B_2\cap\{\dist\bigl(\cdot,\bigcup_{\beta<\gamma}(P_\beta\cap P_\gamma)\bigr)<\rho\}\bigr)
\le C_\alpha\rho^2.
\end{equation}
Choose $\chi\in C_c^1(B_{3/2})$ with $\chi|_{B_{4/3}}=1$ and define
$\eta_\eps$ by
\[
\eta_\eps(x)\coloneqq \chi(x)
\begin{cases}
1,&\dist\bigl(x,\bigcup_{\alpha<\beta}(P_\alpha\cap P_\beta)\bigr)\le\eps^2,\\
\dfrac{\log(\eps/\dist\bigl(x,\bigcup_{\alpha<\beta}(P_\alpha\cap P_\beta)\bigr))}{\log(1/\eps)},
 &\eps^2<\dist\bigl(x,\bigcup_{\alpha<\beta}(P_\alpha\cap P_\beta)\bigr)<\eps,\\
0,&\dist\bigl(x,\bigcup_{\alpha<\beta}(P_\alpha\cap P_\beta)\bigr)\ge\eps.
\end{cases}
\]
Since $\spt\nu\cap B_{3/2}\subset
\bigcup_{\alpha<\beta}(P_\alpha\cap P_\beta)$ and $\chi|_{B_{4/3}}=1$,
applying \eqref{eq:union-curvature-measure} to smooth approximations of
$\eta_\eps$, and using \eqref{eq:plane-intersection-neighborhood-volume}
with the coarea formula, we obtain
\[
\begin{aligned}
\frac\delta n\nu(B_{4/3})
\le&\frac\delta n\int_{B_3}\eta_\eps^2\,d\nu
\le\sum_\alpha q_\alpha\int_{P_\alpha}|\nabla_{P_\alpha}\eta_\eps|^2 
\le C\Bigl(\frac1{|\log\eps|}+\eps^2\Bigr)\longrightarrow0
\quad\text{as }\eps\downarrow0,
\end{aligned}
\]
where $C=C(P_1,\ldots,P_N,q_1,\ldots,q_N)$.
Thus $\nu(B_{4/3})=0$, and by weak convergence and H\"older's inequality,
\[
\int_{F_j^{-1}(B_{5/4})}|A_{M_j}|^{2\alpha}
\le\Lambda^{1-\alpha}
\Bigl(\int_{F_j^{-1}(B_{5/4})}|A_{M_j}|^2\Bigr)^\alpha\longrightarrow0
\quad\text{as }j\to\infty.
\]
For $B_{4s}(x)\Subset B_4$, repeat the logarithmic cutoff argument for
$(F_j-x)/s$ in \eqref{eq:union-curvature-measure}. Then, by
Proposition~\ref{thm:eps-reg} with $r=1/8$,
\[
\sup_{F_j^{-1}(B_s(x))}|A_{M_j}|\longrightarrow0
\quad\text{as }j\to\infty.
\]
By Lemma~\ref{lem:bounded-curvature-removal}, $F_j$ extends across
$S_j$ as a smooth proper special Lagrangian immersion over $B_s(x)$.
Applying Schauder estimates for \eqref{eq:minimal-graph-equations}
to the graph functions of $F_j$, we obtain smooth convergence as
immersions on compact
subsets of $B_4$.

For \textup{(ii)}, apply the dilation $z\mapsto z/2$, using
$F_j,S_j,T_j$ for the rescaled immersions, sets and currents on $B_2$.
Applying Federer--Fleming compactness
\cite[Chapter~6, Theorem~3.11]{Simon2014GMT} and integral varifold
compactness \cite[Chapter~8, Theorem~5.8 and Remark~5.9]{Simon2014GMT}
with $\|T_j\|(B_2)\le2^{-n}\Lambda$ and $\partial T_j=0$, we pass to a
subsequence such that $\Theta_j\to\Theta$ and
\[
T_j\rightharpoonup T\quad\text{locally as currents},\quad
|T_j|\to V\quad\text{locally as varifolds as }j\to\infty.
\]
For $0\le\chi\in C_c^\infty(B_2)$,
\begin{equation}\label{eq:calibrated-mass-chain}
\begin{aligned}
\int_{B_2}\chi\,d\|V\|
&=\lim_{j\to\infty}T_j\bigl(\chi\Re(e^{-i\Theta_j}\Omega)\bigr)
=T\bigl(\chi\Re(e^{-i\Theta}\Omega)\bigr)\\
&\le\int_{B_2}\chi\,d\|T\|
\le\liminf_{j\to\infty}\int_{B_2}\chi\,d\|T_j\|
=\int_{B_2}\chi\,d\|V\|.
\end{aligned}
\end{equation}
Thus $\|V\|=\|T\|$ and $T$ is calibrated by
$\Re(e^{-i\Theta}\Omega)$. Since $V$ is integral, $V=|T|$.

Let $B_2^{\mathrm{reg}}\coloneqq B_2\setminus\ImmSing T$.
If $T\llcorner B_2^{\mathrm{reg}}=0$, take $N=\varnothing$.
Otherwise, we first use \textup{(i)} to prove that, for each
$K\Subset B_2^{\mathrm{reg}}$,
\begin{equation}\label{eq:regular-set-local-curvature}
\sup_{F_j^{-1}(K)}|A_{M_j}|\le C(K)\qquad\forall\,j\ge j_0(K).
\end{equation}
Indeed, if not, choose $p_j\in M_j$ with
$F_j(p_j)\in K$, $|A_{M_j}|(p_j)\to\infty$, and
$F_j(p_j)\to x\in K$ as $j\to\infty$, where $x\in\spt T$ by \eqref{eq:stationary-support-mass}.
By Definition~\ref{def:imm-reg}, the tangent
cone of $T$ at $x$ is a finite union of planes $P_1,\ldots,P_N$
with positive integer multiplicities $q_1,\ldots,q_N$.
After passing to a subsequence, choose $r_j\downarrow0$ diagonally so that
\begin{equation}\label{eq:regular-set-rescaling-data}
\begin{gathered}
r_j|A_{M_j}|(p_j)\longrightarrow\infty,\qquad
\frac{|F_j(p_j)-x|}{r_j}\longrightarrow0,\\
|(D_{1/r_j}\circ\tau_{-x})_\#T_j|
\longrightarrow\sum_{\alpha=1}^Nq_\alpha|P_\alpha|
\quad\text{as }j\to\infty.
\end{gathered}
\end{equation}
By \eqref{eq:stationary-monotonicity},
\[
\bigl\|(D_{1/r_j}\circ\tau_{-x})_\#T_j\bigr\|(B_4)
=r_j^{-n}\|T_j\|(B_{4r_j}(x))\le C.
\]
By \eqref{eq:regular-set-rescaling-data}, $F_j(p_j)\in B_{r_j}(x)$
for all large $j$. Applying \textup{(i)} to $(F_j-x)/r_j$, we obtain
\[
r_j|A_{M_j}|(p_j)
\le r_j\sup_{F_j^{-1}(B_{r_j}(x))}|A_{M_j}|
\longrightarrow0\quad\text{as }j\to\infty,
\]
contradicting \eqref{eq:regular-set-rescaling-data}.

For each $K\Subset B_2^{\mathrm{reg}}$ and all large $j$, by
\eqref{eq:regular-set-local-curvature} and
Lemma~\ref{lem:bounded-curvature-removal}, we extend $F_j$ smoothly
across $S_j$ near $K$. The extended immersions are again denoted by $F_j$.
By Schauder estimates for \eqref{eq:minimal-graph-equations}, applied
to the graph functions of $F_j$, for $k\ge0$,
\[
\sup_{F_j^{-1}(K)}|\nabla^kA_{M_j}|\le C(K,k).
\]
Applying \cite[Corollaries~1.4 and~7.13]{Breuning2015} to the proper
extensions $F_j$ on an exhaustion of $B_2^{\mathrm{reg}}$ and choosing
a diagonal subsequence, we obtain a proper special Lagrangian immersion
$G:N\to B_2^{\mathrm{reg}}$ of phase $\Theta$ without boundary,
open sets $U_j\Subset U_{j+1}\nearrow B_2^{\mathrm{reg}}$ and
$N_j\nearrow N$, and diffeomorphisms
\[
\Phi_j:N_j\longrightarrow F_j^{-1}(U_j),\qquad
F_j\circ\Phi_j\longrightarrow G
\quad\text{in }C^\infty_{\rm loc}(N)\text{ as }j\to\infty.
\]
By Lemma~\ref{lem:bounded-curvature-removal}, the extended immersions
$F_j$ represent $T_j$. Combining $F_j\circ\Phi_j\to G$ with
$T_j\rightharpoonup T$, we therefore have
\[
G_\#d\mu_N=\|T\|\llcorner B_2^{\mathrm{reg}},
\quad
G_\#\llbracket N\rrbracket=T\llcorner B_2^{\mathrm{reg}}.
\]

For $\phi\in C_c^1(N)$ and all large $j$, $\spt\phi\Subset N_j$.
Define
$\phi_j\coloneqq \phi\circ\Phi_j^{-1}$ on $F_j^{-1}(U_j)$ and
$\phi_j\coloneqq 0$ on $M_j\setminus F_j^{-1}(U_j)$.
Since $\Phi_j(\spt\phi)\Subset F_j^{-1}(U_j)$, $\phi_j\in C_c^1(M_j)$.
By smooth convergence and \eqref{eq:P-estimates},
$\Pcal_j\circ\Phi_j\to\Pcal_G$ uniformly on $\spt\phi$ as $j\to\infty$, including
points where $A_N=0$. Since the extensions $F_j$ are
$\delta$-Simons stable by Lemma~\ref{lem:bounded-curvature-removal},
using $\phi_j$ in \eqref{eq:lipschitz-stability}, we have
\begin{equation}\label{eq:regular-limit-stability}
\begin{aligned}
\delta\int_N\Pcal_G\phi^2\,d\mu_N
&=\lim_{j\to\infty}\delta\int_{M_j}\Pcal_j\phi_j^2\,d\mu_j 
 \le\lim_{j\to\infty}\int_{M_j}|\nabla\phi_j|^2\,d\mu_j
=\int_N|\nabla\phi|^2\,d\mu_N.
\end{aligned}
\end{equation}
Therefore, after rescaling back to $B_4$, the extensions of $F_j$ converge smoothly
to the proper $\delta$-Simons stable immersion
$G:N\to B_4^{\mathrm{reg}}$ representing
$T\llcorner B_4^{\mathrm{reg}}$.
\end{proof}

We next prove cone rigidity for limits of the immersions in
Proposition~\ref{thm:flat-union}.

For $\sigma<n-4+6/(n+1)$, let $\mathcal C(n,\delta,\sigma)$ consist
of integral currents $T$ for which there exist open sets
$U_j\supset B_j(0)$ and smooth proper
$\delta$-Simons stable special Lagrangian immersions
$F_j:M_j\to U_j\setminus S_j$, without boundary over $U_j\setminus S_j$,
where $S_j$ is relatively closed, $\dim_{\mathcal H}S_j\le\sigma$,
and all components have a common phase $\Theta_j$, converging modulo
$2\pi$ to the phase of $T$.
The currents induced by $F_j$ on $U_j\setminus S_j$ extend to integral
currents $T_j$ with $\partial T_j=0$ in $U_j$, and
\[
T_j\rightharpoonup T\quad\text{locally as currents},\qquad
|T_j|\longrightarrow|T|\quad\text{locally as varifolds as }j\to\infty.
\]
$\mathcal C(n,\delta,\sigma)$ is invariant under translations and positive
dilations and, by diagonal selection, closed under simultaneous local
current and varifold limits, tangent cones and iterated tangent cones.
For local limits as in Proposition~\ref{thm:flat-union}\textup{(ii)}, interior
tangent cones belong to $\mathcal C(n,\delta,\sigma)$: choose tangent radii tending
to zero and, at each radius, one immersion from a sequence with current
limit $T$ and varifold limit $|T|$, rescaled onto expanding balls.

\begin{proposition}\label{prop:closure-cone-reduction}
Fix $n\ge3$, $\delta>\delta_{\rm K}(n)$, and
$\sigma<n-4+6/(n+1)$. Let
$C\in\mathcal C(n,\delta,\sigma)$ be a cone invariant under translations
along an $s$-plane, and set $m\coloneqq n-s$. If $m\le2$, or if $m\ge3$ and
$\delta>\delta_{\rm cone}(m)$, then $C$ is a finite planar union with
a common phase.
\end{proposition}

\begin{proof}
Let $E$ be the $s$-dimensional translation space. Since $T_xC$ is
Lagrangian and contains $E$ and $x$ for $\|C\|$-almost every $x$, we have
\[
\langle x,Jv\rangle=0\quad\forall\,v\in E.
\]
Thus, after a unitary change of coordinates, the slice decomposition is
\begin{equation}\label{eq:closure-cone-product}
C=\llbracket E\rrbracket\times C_0,\quad
C_0\subset(E\oplus JE)^\perp\simeq\C^m.
\end{equation}
Here $C_0$ has the same special Lagrangian phase as $C$.
If $m=0$, $C$ is an integer multiple of $\llbracket E\rrbracket$. If $m=1$, let $P_+$ and $P_-$
be the two opposite rays in the unique special Lagrangian line of the
prescribed phase. Then
\[
C_0=q_+\llbracket P_+\rrbracket+q_-\llbracket P_-\rrbracket,
\quad \partial C_0=0,
\]
so $q_+=q_-$ and $C_0=q_+\llbracket P_+\cup P_-\rrbracket$.

Induct on $m\ge2$. Suppose that $x\in\ImmSing C\setminus E$.
Since $C$ is calibrated and $\partial C=0$, $C$ is locally area minimizing.
By \cite[Chapter~7, Theorem~3.1]{Simon2014GMT}, choose a tangent cone
\[
\begin{gathered}
C_{x,r_j}\coloneqq (D_{1/r_j}\circ\tau_{-x})_\#C
\rightharpoonup C_\infty,\\
|C_{x,r_j}|\longrightarrow|C_\infty|,\qquad
r_j\downarrow0\quad\text{as }j\to\infty.
\end{gathered}
\]
Each $C_{x,r_j}$ is invariant under translations along $E$. For
$v\in E$, $t\in\R$, and $j\ge j_0(t)$, conical invariance gives
$(\tau_{v+tx})_\#C_{x,r_j}=(D_{1-r_jt})_\#C_{x,r_j}$.
Since $1-r_jt\to1$ as $j\to\infty$, local current convergence implies
\[
\begin{aligned}
(\tau_{v+tx})_\#C_\infty
=&\lim_{j\to\infty}(\tau_{v+tx})_\#C_{x,r_j} 
= \lim_{j\to\infty}(D_{1-r_jt})_\#C_{x,r_j}=C_\infty.
\end{aligned}
\]
Since $x\notin E$, the factor transverse to $E\oplus\R x$ has dimension
at most $m-1$. For $3\le k<m$, $\delta_{\rm cone}(k)\le\delta_{\rm cone}(m)<\delta$,
so $C_\infty$ is a finite planar union by induction.

Choose for each $j$ a sequence of immersions whose integral current
extensions converge to $C_{x,r_j}$ and whose varifolds converge to
$|C_{x,r_j}|$, with initial singular sets of dimension at most $\sigma$
and masses on $B_4$ at most
$2\sup_k\|C_{x,r_k}\|(B_5)<\infty$, using local convergence of the
mass measures to $\|C_{x,r_j}\|$ and
\eqref{eq:stationary-monotonicity}.
We now prove $x\in\ImmReg C$. Fix
$\Lambda>2\sup_k\|C_{x,r_k}\|(B_5)$.
On varifolds $V$ in an open set $U$ with $\|V\|(U)\le\Lambda<\infty$,
let $d_U$ denote any metric inducing weak varifold convergence.

We claim that, for the current limit $T$ in
Proposition~\ref{thm:flat-union}\textup{(ii)},
there exist $\eps_{\rm flat}=\eps_{\rm flat}(n,\delta,\Lambda,\sigma)>0$ and
$c_{\rm flat}=c_{\rm flat}(n,\delta,\Lambda,\sigma)>0$ such that, for
$0\in\spt T$ and a finite planar integral cone $D$ with a common phase
and $\|D\|(B_4)\le\Lambda$,
\begin{equation}\label{eq:uniform-flat-neighborhood}
d_{B_2}(|T|,|D|)<\eps_{\rm flat}
\quad\Longrightarrow\quad r_{\rm reg}^T(0)\ge c_{\rm flat}.
\end{equation}

Indeed, if no such $\eps_{\rm flat},c_{\rm flat}>0$ exist, choose examples $T^{(k)}$
and finite planar integral cones $D_k$ with a common phase and
$\|D_k\|(B_4)\le\Lambda$, such that
$0\in\spt T^{(k)}$, $\|T^{(k)}\|(B_4)\le\Lambda$,
$d_{B_2}(|T^{(k)}|,|D_k|)\to0$, and
$r_{\rm reg}^{T^{(k)}}(0)\to0$ as $k\to\infty$.
The number of planes in $D_k$, counted with multiplicity, is
\[
\frac{\|D_k\|(B_4)}{\omega_n4^n}\le\frac{\Lambda}{\omega_n4^n}.
\]
After passing to a subsequence in the
Lagrangian Grassmannian, let
$P_1,\ldots,P_N$ be the distinct limit planes, with a common phase.
Then
\[
D_\infty=\sum_{\alpha=1}^Nq_\alpha\llbracket P_\alpha\rrbracket,
\quad q_\alpha\in\mathbb N.
\]

For each $k$, choose immersions $L_{k,\ell}$ with initial singular sets of dimension at most $\sigma$. The integral current extensions $T_{k,\ell}$ of $L_{k,\ell}$ satisfy
$\partial T_{k,\ell}=0$ in $B_4$ and converge as currents to $T^{(k)}$
and as varifolds to $|T^{(k)}|$. Choose points
$x_{k,\ell}\in\spt T_{k,\ell}$ with $x_{k,\ell}\to0$ as $\ell\to\infty$.
If $r_{\rm reg}^{T_{k,\ell}}(x_{k,\ell})\ge c>0$ uniformly, then
$r_{\rm reg}^{T^{(k)}}(0)\ge c/2$ by local immersion compactness
\cite[Corollary~7.13]{Breuning2015} and Definition~\ref{def:imm-reg}.
Choose $\ell(k)\to\infty$ as $k\to\infty$ such that
\[
\begin{aligned}
|(\tau_{-x_{k,\ell(k)}})_\#T_{k,\ell(k)}|
&\longrightarrow|D_\infty|\quad\text{as varifolds on }B_2,\\
r_{\rm reg}^{T_{k,\ell(k)}}(x_{k,\ell(k)})&\longrightarrow0\quad\text{as }k\to\infty.
\end{aligned}
\]
Applying Proposition~\ref{thm:flat-union}\textup{(i)} to
$2(L_{k,\ell(k)}-x_{k,\ell(k)})$, we obtain
$r_{\rm reg}^{T_{k,\ell(k)}}(x_{k,\ell(k)})\ge c>0$, a contradiction.
Thus \eqref{eq:uniform-flat-neighborhood} holds.

Applying \eqref{eq:uniform-flat-neighborhood} with
$T=C_{x,r_j}$ and $D=C_\infty$, using
$|C_{x,r_j}|\to|C_\infty|$ on $B_2$ as $j\to\infty$, we obtain
$r_{\rm reg}^{C_{x,r_j}}(0)\ge c_{\rm flat}>0$ for $j\ge j_0$ and,
after scaling back, $r_{\rm reg}^{C}(x)\ge r_jc_{\rm flat}>0$,
contradicting $x\in\ImmSing C$. Therefore
$C\setminus E\subset\ImmReg C$.

Using the exhaustion construction in
\cite[Corollary~7.13]{Breuning2015} for a sequence of immersions $F_j$
defining $C$ on compact subsets of $\C^n\setminus E$, with
$\sup_{F_j^{-1}(K)}|\nabla^kA_{M_j}|\le C(K,k)$ for
$K\Subset\C^n\setminus E$, $k\ge0$, and all large $j$,
we obtain a manifold $N$ and a proper immersion
$G:N\to\C^n\setminus E$ satisfying
\[
G_\#\llbracket N\rrbracket=C\llcorner(\C^n\setminus E).
\]
By \eqref{eq:regular-limit-stability}, the immersion $G$ is
$\delta$-Simons stable.
By \eqref{eq:closure-cone-product} and conical invariance,
$E\oplus\R G(p)\subset dG(T_pN)$ for almost every $p\in N$, and
hence for all $p$ on each smooth sheet.

Define $X_v$ for $v\in E$ by $dG(X_v)=v$. Decomposing
$\overline\nabla v=0$ into tangent and normal parts,
\[
\nabla X_v=0,\quad A_N(X_v,\cdot)=0.
\]
Consequently, $[X_v,X_w]=0$ for $v,w\in E$. If $\gamma'=X_v\circ\gamma$, then
\[
G(\gamma(t))=G(\gamma(0))+tv.
\]
For $a<b$, $\{G(\gamma(0))+tv:a\le t\le b\}\Subset\C^n\setminus E$.
Since $G$ is proper, $\gamma$ remains in a compact subset of $N$ on
every bounded time interval, so the flow of $X_v$ extends to all $t\in\R$.
Let $N_0\coloneqq(\operatorname{proj}_E\circ G)^{-1}(0)$.
For each $p\in N$, the flow of
$X_{-\operatorname{proj}_E G(p)}$ for unit time reaches $N_0$;
the displacement along $E$ is unique because
$d(\operatorname{proj}_E\circ G)(X_v)=v$.
Since $\nabla X_v=0$, the commuting flows preserve the metric. Thus
\[
N=\R^s\times N_0,\quad
G(y,p)=y+G_0(p),
\]
where $G_0:N_0\to(E\oplus JE)^\perp\setminus\{0\}$ is proper.

For the radial direction, let $Y$ satisfy $dG_0(Y)=G_0$. Along $\gamma'=Y\circ\gamma$,
$G_0(\gamma(t))=e^tG_0(\gamma(0))$.
As for $X_v$, the flow exists for all $t\in\R$ by properness of $G_0$. With $\Sigma\coloneqq G_0^{-1}(\Sph^{2m-1})$,
\[
N_0=(0,\infty)\times\Sigma,\quad
G_0(r,p)=rG_0(1,p),\quad
g_{N_0}=dr^2+r^2g_\Sigma.
\]
When $s>0$, for $0\ne\chi\in C_c^1(\R^s)$ and $\phi\in C_c^1(N_0)$,
substituting $\chi(y/R)\phi(p)$ into \eqref{eq:regular-limit-stability}
and dividing by
$R^s\int_{\R^s}\chi^2\,dy$, we have, since $\Pcal_G=\Pcal_{G_0}$,
\[
\begin{aligned}
\delta\int_{N_0}\Pcal_{G_0}\phi^2\,d\mu_{N_0}
\le&\int_{N_0}|\nabla\phi|^2\,d\mu_{N_0}
+\frac{\int_{\R^s}|\nabla\chi|^2\,dy}
 {R^2\int_{\R^s}\chi^2\,dy}\int_{N_0}\phi^2\,d\mu_{N_0}\\
\longrightarrow&\int_{N_0}|\nabla\phi|^2\,d\mu_{N_0}
\quad\text{as }R\to\infty.
\end{aligned}
\]
For $s=0$, $G=G_0$, and we use \eqref{eq:regular-limit-stability}
on $N_0$. Thus $G_0$ is $\delta$-Simons stable.
By properness of $G_0$, $\Sigma$ is compact. Since $G_0$ is a special
Lagrangian cone without boundary, $\Sigma$ is minimal Legendrian
and without boundary.

If $m=2$, each component of $\Sigma$ is a closed geodesic and hence $A_{N_0}=0$.
If $m\ge3$, applying Theorem~\ref{thm:intro-cone} to
$G_0:(0,\infty)\times\Sigma_j\to(E\oplus JE)^\perp$ for each component
$\Sigma_j$ of $\Sigma$, we again have $A_{N_0}=0$.
For each connected component $\Sigma_j$, choose a special Lagrangian
$m$-plane $P_j$ with the phase of $C_0$ such that
\[
G_0(1,\cdot)|_{\Sigma_j}:\Sigma_j\longrightarrow
P_j\cap\Sph^{2m-1}
\]
is a finite covering of degree $d_j\in\mathbb N$.
Since $\Sigma$ is compact, only finitely many $\Sigma_j$ occur.
Also, $\dim E=n-m<n$, so $\|C\|(E)=0$.
Therefore $C=\sum_jd_j\llbracket E\oplus P_j\rrbracket$ is a finite planar
union with a common phase.
\end{proof}

We are now ready to prove Theorem~\ref{thm:intro-regularity-compactness}.

\begin{proof}[Proof of Theorem~\ref{thm:intro-regularity-compactness}]
By Proposition~\ref{thm:flat-union}\textup{(ii)}, there exists a subsequence
such that $T_j\rightharpoonup T$ as currents and $|T_j|\to|T|$ as varifolds
locally in $B_2$. Moreover, along the same subsequence, the extensions
of $F_j$ converge smoothly as immersions on $\ImmReg T$ to a
$\delta$-Simons stable immersion. It remains to estimate $\ImmSing T$.

Fix $\sigma$ as in \eqref{eq:initial-singular-dimension}.
By diagonal selection, the collection of local current and varifold
limits of immersions satisfying the hypotheses of
Theorem~\ref{thm:intro-regularity-compactness} with
$\dim_{\mathcal H}S_j\le\sigma$ is closed
under further local limits, translations and positive dilations, and the
tangent cones of the local current limits in
Proposition~\ref{thm:flat-union}\textup{(ii)} belong to $\mathcal C(n,\delta,\sigma)$.
For local current and varifold limits $T^{(\ell)},T$ of immersions
satisfying the hypotheses of Theorem~\ref{thm:intro-regularity-compactness}
with $\dim_{\mathcal H}S_j\le\sigma$, we claim that
\begin{equation}\label{eq:immersed-singular-persistence}
\begin{gathered}
T^{(\ell)}\rightharpoonup T,\quad |T^{(\ell)}|\to|T|
\quad\text{as }\ell\to\infty,\\
x_\ell\in\ImmSing T^{(\ell)},\quad x_\ell\to x\text{ as }\ell\to\infty
\quad\Longrightarrow\quad x\in\ImmSing T.
\end{gathered}
\end{equation}
Indeed, \eqref{eq:stationary-support-mass} implies $x\in\spt T$.
If $x\in\ImmReg T$, the tangent cone of $T$ at $x$ is a finite planar union
$C$ with a common phase. Choose $\Lambda_0>2\|C\|(B_5)$.
Since $r^{-n}\|T\|(B_{5r}(x))\to\|C\|(B_5)$ as $r\downarrow0$,
local convergence of the mass measures implies, for sufficiently
small $r>0$ and all large $\ell$,
\[
r^{-n}\max\bigl\{\|T\|(B_{4r}(x)),
\|T^{(\ell)}\|(B_{4r}(x_\ell))\bigr\}\le\Lambda_0.
\]
Choose $r>0$ so small that
$d_{B_2}(|(D_{1/r}\circ\tau_{-x})_\#T|,|C|)<\eps_{\rm flat}/2$, with
$\eps_{\rm flat}$ in \eqref{eq:uniform-flat-neighborhood} for
$\Lambda=\Lambda_0$.
For large $\ell$,
\[
d_{B_2}\bigl(|(D_{1/r}\circ\tau_{-x_\ell})_\#T^{(\ell)}|,|C|\bigr)
<\eps_{\rm flat}.
\]
Applying
\eqref{eq:uniform-flat-neighborhood}, we have
$r_{\rm reg}^{T^{(\ell)}}(x_\ell)\ge c_{\rm flat}r>0$, contradicting
$x_\ell\in\ImmSing T^{(\ell)}$.

Since the local current limits are calibrated integral currents with
zero boundary and hence locally area minimizing, we obtain conical tangent
limits by \cite[Chapter~7, Theorem~3.1]{Simon2014GMT}.
By \eqref{eq:immersed-singular-persistence}, we may apply the dimension
reduction theorem \cite[Appendix~A, Theorem~A.4]{Simon2014GMT} to the local current and
varifold limits of immersions satisfying the hypotheses of
Theorem~\ref{thm:intro-regularity-compactness} with
$\dim_{\mathcal H}S_j\le\sigma$, using $\ImmSing$ as the singular set.
If any immersed singular point occurs, we obtain an
integer $d\in\{0,\ldots,n-1\}$ such that
\[
\dim_{\mathcal H}(\ImmSing T)\le d
\]
for every local limit, and a cone $C\in\mathcal C(n,\delta,\sigma)$ invariant
under translations along a $d$-plane $E$, with
$\ImmSing C=E$.
If $n-d<n_\delta$, then $n-d\le2$ or
$\delta>\delta_{\rm cone}(n-d)$ by
\eqref{eq:compactness-critical-dimension}. By
Proposition~\ref{prop:closure-cone-reduction}, $C$ is a finite planar
union. By Definition~\ref{def:imm-reg},
\[
0\in E=\ImmSing C=\varnothing,
\]
a contradiction. Hence
\[
0\le d\le n-n_\delta.
\]
In particular, $S=\ImmSing T\cap B_1$ is empty for $n<n_\delta$, and
$\dim_{\mathcal H}S\le n-n_\delta$ for $n>n_\delta$.
For $n=n_\delta$, $S$ is locally finite in $B_1$ by
\cite[Theorem~A.4]{Simon2014GMT}.
\end{proof}

We next prove Corollary~\ref{cor:curvature-and-growth-bernstein}.

\begin{proof}[Proof of Corollary~\ref{cor:curvature-and-growth-bernstein}]
Here $n<n_\delta$. If \eqref{eq:curvature-mass-bound} fails, choose immersions $F_j$
satisfying the hypotheses of the corollary with
\[
\sup_{F_j^{-1}(B_1)}|A_{M_j}|\longrightarrow\infty\quad\text{as }j\to\infty.
\]
Applying Theorem~\ref{thm:intro-regularity-compactness} on interior
balls in $B_2$, we pass to a subsequence for which the currents induced
by $F_j$ converge locally to an integral current $T$ and
$\ImmSing T\cap B_2=\varnothing$.
By Proposition~\ref{thm:flat-union}\textup{(ii)},
\[
\sup_{F_j^{-1}(\overline B_1)}|A_{M_j}|\le C\qquad\forall\,j\ge j_0,
\]
a contradiction.

For the complete immersion, applying the Michael--Simon inequality
and the coarea formula on intrinsic balls, we obtain, for some $c(n)>0$,
\[
\mu(B_s^g(p))\ge c(n)s^n\qquad\forall\,p\in M,\ \forall\,s>0.
\]
If $F$ is not proper, some sequence escaping compact subsets of $M$
has image in a fixed Euclidean ball. By completeness, we can
choose disjoint intrinsic unit balls around a subsequence; the images
of the unit balls lie in one Euclidean ball, contradicting the volume growth
hypothesis. Hence $F$ is proper. Applying
\eqref{eq:curvature-mass-bound} to $F/R$ on
$F^{-1}(B_{2R}(0))$, we have
\[
\sup_{F^{-1}(B_R(0))}|A_M|\le\frac{C(n,\delta,2^n\Lambda)}R
\longrightarrow0\quad\text{as }R\to\infty.
\]
Thus $F$ is flat.
\end{proof}

\section{Regularity and compactness in dimension three}\label{sec:small-singular-sets}

\begin{theorem}\label{thm:three-dimensional-null-set-compactness}
Let $\delta>1/4$ and let $S_j\subset B_2\subset\C^3$ be
relatively closed with $\mathcal H^1(S_j)=0$.
Let $F_j:M_j^3\to B_2\setminus S_j$ be smooth proper
$\delta$-Simons stable special Lagrangian immersions without
boundary, with phase $\Theta_j$ on all components of $M_j$.
Assume that the currents induced by $F_j$ extend to integral currents $T_j$ with $\partial T_j=0$ in $B_2$ and
\[
\sup_j\|T_j\|(K)<\infty\qquad\forall\,K\Subset B_2.
\]
Then there exists a subsequence such that
\[
\Theta_j\longrightarrow\Theta\pmod{2\pi},\qquad
T_j\rightharpoonup T,\qquad |T_j|\longrightarrow|T|
\quad\text{locally in }B_2\text{ as }j\to\infty.
\]
Moreover, along the same subsequence, the extensions of $F_j$ across
the removable singularities in $S_j$ converge smoothly as immersions on
compact subsets of $B_2$ to a proper $\delta$-Simons stable special
Lagrangian immersion of phase $\Theta$ representing $T$.
\end{theorem}

\begin{proof}
Since $\mathcal H^1(S_j)=0$ and $\partial T_j=0$, $T_j$ is calibrated
on $B_2$ and $|T_j|$ is stationary. Using
$\sup_j\|T_j\|(K)<\infty$ for $K\Subset B_2$ and applying Federer--Fleming compactness
\cite[Chapter~6, Theorem~3.11]{Simon2014GMT} and integral varifold
compactness \cite[Chapter~8, Theorem~5.8 and Remark~5.9]{Simon2014GMT},
we pass to a subsequence with the stated phase and current limits.
By \eqref{eq:calibrated-mass-chain}, the varifold limit is $|T|$.
Since $\mathcal H^1(S_j)=0$, the covering balls in the construction of
\eqref{eq:removal-cutoff-estimates} may be chosen with $\sum_i r_i\to0$.
For $n=3$ and $q=2$, \eqref{eq:removal-cutoff-estimates} then follows
from \eqref{eq:stationary-monotonicity}. For $\phi$ compactly supported
on a smooth multiplicity one sheet of $T_j$, substitute
$(\chi_\ell\circ F_j)(\phi\circ F_j)$ into
\eqref{eq:delta-stability}. Letting $\ell\to\infty$ and using
\eqref{eq:removal-cutoff-estimates}, we obtain stability across $S_j$
on that sheet. For local current and varifold limits of $T_j$,
stability on the multiplicity one regular part follows from
\eqref{eq:multiplicity-one-limit-stability} below.
Lemma~\ref{lem:low-density-SLag-regularity} therefore
implies smooth embeddedness at points of density less than two. If $\Theta^3(\|T\|,x)<2$,
by upper semicontinuity choose $r>0$ and $j_0$ with
\[
B_r(x)\cap\{y:\Theta^3(\|T_j\|,y)\ge2\}=\varnothing
\qquad\forall\,j\ge j_0.
\]
By Lemma~\ref{lem:low-density-SLag-regularity},
$T_j\llcorner B_r(x)$ is smoothly embedded with multiplicity one.
Applying Allard's regularity theorem
\cite[Chapter~5, Theorem~5.2]{Simon2014GMT} and then Schauder estimates
to \eqref{eq:minimal-graph-equations}, we extend $F_j$ across
$S_j\cap B_{r/2}(x)$ and obtain smooth convergence on compact subsets
of $B_{r/2}(x)$.

For the induction step, assume compactness below density $Q$, where
$Q\ge2$. If $|T_j|$ converges to a planar cone of total multiplicity at
most $Q$, then $F_j$ extends smoothly
across $S_j\cap B_1$ for all large $j$, and
\begin{equation}\label{eq:null-set-planar-extension-bound}
\sup_{F_j^{-1}(B_1)}|A_{M_j}|\longrightarrow0\qquad\text{as }j\to\infty.
\end{equation}
Indeed, apply Proposition~\ref{prop:three-dimensional-planar-induction-step}
at total multiplicity $Q$ and the induction hypothesis at smaller
multiplicity.
For a fixed planar cone $C$ of total multiplicity at most $Q$ and
$0<\Lambda<\infty$, if $\sup_j\|T_j\|(B_2)\le\Lambda$ and $|T_j|$ lies
in a sufficiently small varifold neighborhood of $|C|$ on $B_2$, then
$F_j$ extends across $S_j\cap B_1$ and
$\sup_{F_j^{-1}(B_1)}|A_{M_j}|\le1$.
Otherwise, choose a sequence converging to $|C|$ for which $F_j$ fails
to extend across $S_j\cap B_1$ with $\sup_{F_j^{-1}(B_1)}|A_{M_j}|\le1$, contradicting
\eqref{eq:null-set-planar-extension-bound}.
By Schauder estimates for \eqref{eq:minimal-graph-equations} and
\cite[Corollaries~1.4 and~7.13]{Breuning2015},
local current and varifold limits of the extensions of $F_j$ are
represented by smooth immersions on $B_1$.

By diagonal selection, the collection of local current and varifold
limits of the three dimensional immersions $F_j:M_j^3\to B_2\setminus S_j$
satisfying the hypotheses of
Theorem~\ref{thm:three-dimensional-null-set-compactness}, with
$\mathcal H^1(S_j)=0$, is closed under further local limits, translations
and positive dilations.

We claim that if $C$ is a cone obtained as a local current and varifold
limit of immersions satisfying the hypotheses of
Theorem~\ref{thm:three-dimensional-null-set-compactness} and
$\Theta^3(\|C\|,0)<Q+1$, then $C$ is a planar union of total multiplicity
at most $Q$.
Indeed, argue by downward induction on the dimension $s$ of a
translation space $E$ of $C$. For $s=3$, $C$ is a plane. For $s=2$,
the factor $C_0$ in \eqref{eq:closure-cone-product} has dimension one
and satisfies $\partial C_0=0$, so $C_0$ consists of two opposite rays
with equal integer multiplicity. Suppose the claim holds for translation spaces of
dimension at least $s+1$. For $x\in\spt C\setminus E$, every tangent
cone $D$ satisfies
\[
(\tau_v)_\#D=D\quad\forall\,v\in E+\mathbb Rx,\qquad
\Theta^3(\|D\|,0)\le\Theta^3(\|C\|,0)<Q+1.
\]
By downward induction, $D$ is a planar union of total multiplicity at most
$Q$. Applying \eqref{eq:null-set-planar-extension-bound} to rescalings
about $x$, we obtain a smooth immersion representing $C$ near $x$.
By an exhaustion and \cite[Corollary~7.13]{Breuning2015}, $C$ is represented
by a proper smooth immersion on $\C^3\setminus E$, stable by
\eqref{eq:regular-limit-stability}.

It follows from \eqref{eq:closure-cone-product} and the argument for
proper smooth immersions in the proof of
Proposition~\ref{prop:closure-cone-reduction} that the factor of $C$ in $(E\oplus JE)^\perp\simeq\C^{3-s}$ is a smooth
$\delta$-Simons stable cone with compact immersion domain over the unit
sphere. For $3-s=3$, apply Theorem~\ref{thm:intro-cone}, since
$\delta>1/4>\delta_{\rm cone}(3)$; for $3-s=2$, each component is a cone
over a closed geodesic.
Thus $C$ is a finite planar union, with integer density at most $Q$.
This proves the claim.

By the claim, at $x\in\spt T$ with $\Theta^3(\|T\|,x)<Q+1$, every tangent cone is a
planar union of total multiplicity at most $Q$.
Rescaling about $x$ and using \eqref{eq:null-set-planar-extension-bound},
we extend $F_j$ smoothly across $S_j\cap B_r(x)$ for some fixed $r>0$
and all large $j$, with
\[
\sup_{F_j^{-1}(B_r(x))}|A_{M_j}|\le r^{-1}.
\]
To conclude smooth immersed compactness below density $Q+1$, cover
each compact subset by finitely many of these balls and take an exhaustion.

For $x\in\spt T$ and $\overline B_r(x)\Subset B_2$,
\eqref{eq:stationary-monotonicity} implies
$\Theta^3(\|T\|,x)\le\|T\|(B_r(x))/(\omega_3r^3)<\infty$.
At each $x\in\spt T$, choose an integer $Q>\Theta^3(\|T\|,x)$ and apply
the estimate $\sup_{F_j^{-1}(B_r(x))}|A_{M_j}|\le r^{-1}$ just proved.
Schauder estimates for \eqref{eq:minimal-graph-equations} control the
higher derivatives of the local graph functions. Applying
\cite[Corollaries~1.4 and~7.13]{Breuning2015} on an exhaustion, after a
further subsequence, the extensions of $F_j$ converge smoothly on compact subsets of $B_2$ to a
proper $\delta$-Simons stable special Lagrangian immersion of phase
$\Theta$ representing $T$, using
\eqref{eq:regular-limit-stability} for stability.
\end{proof}

When $|T_j|$ tends to a plane with integer multiplicity, we first prove
smooth convergence on compact subsets after extending $F_j$ across $S_j$.

\begin{proposition}\label{prop:null-set-plane-regularity}
Let $n\ge3$ and $\delta>\delta_{\rm K}(n)$.
Suppose $S_j\subset B_2$ are relatively closed with
$\mathcal H^{n-2}(S_j)=0$, and let $F_j:M_j^n\to B_2\setminus S_j$
be smooth proper $\delta$-Simons stable special Lagrangian immersions
without boundary, with one phase on all components of each $M_j$.
Assume that the currents induced by $F_j$ extend to integral currents $T_j$ with $\partial T_j=0$ in $B_2$,
$\|T_j\|(B_2)\le\Lambda$, and $|T_j|\to q|P|$ locally as varifolds as $j\to\infty$,
where $P$ is a linear special Lagrangian $n$-plane and $q\in\mathbb N$.
Then, on each compact subset of $B_2$ and for all large $j$, the
immersions extend smoothly across $S_j$ and converge smoothly to $P$
with total multiplicity $q$. In particular,
\[
\sup_{F_j^{-1}(K)}|A_{M_j}|\longrightarrow0
\quad\text{as }j\to\infty\qquad\forall\,K\Subset B_2.
\]
\end{proposition}

\begin{proof}
Fix $j$ and omit the index while deriving \eqref{eq:bounded-tilt-estimate}.
For an orthonormal basis $p_1,\ldots,p_n$ of $P$, let
\[
b\coloneqq\Bigl(\sum_a|\nabla\langle F,Jp_a\rangle|^2\Bigr)^{1/2},
\qquad 0\le b\le\sqrt n.
\]
Choose
\[
\frac{n-2}{n-1}<\alpha<
 \min\Bigl\{1,\frac{n+1}{n-1}\delta\Bigr\}.
\]

We claim that constants $\eps_0,C>0$, depending only on $n,\delta,\alpha$,
exist such that $\int_{F^{-1}(B_1)}b^{2\alpha}\,d\mu\le\eps_0$ implies
\begin{equation}\label{eq:bounded-tilt-estimate}
\sup_{F^{-1}(B_{1/2})}b
\le C\Bigl(\int_{F^{-1}(B_1)}b^{2\alpha}\,d\mu\Bigr)^{1/(2\alpha)}.
\end{equation}

Indeed, at an arbitrary point $x\in M$, choose orthonormal bases
$p_1,\ldots,p_n$ of $P$ and $e_1,\ldots,e_n$ of $T_xM$ such that
$e_i=\cos\theta_i\,p_i+\sin\theta_i\,Jp_i$.
Since $b^2=\sum_i\sin^2\theta_i$ and $\Delta\langle F,Jp_a\rangle=0$,
applying Bochner's formula, we obtain on $\{b>0\}$
\begin{equation}\label{eq:normal-coordinate-bochner}
\begin{aligned}
 b\Delta b+|\nabla b|^2=\frac12\Delta b^2 
=&\sum_a\Bigl(
 |\nabla^2\langle F,Jp_a\rangle|^2
 +\Ric_g\bigl(\nabla\langle F,Jp_a\rangle,
              \nabla\langle F,Jp_a\rangle\bigr)\Bigr)\\
=&\sum_{i,j,a}\cos^2\theta_a\,(h_M)_{ija}^2
 +\sum_a\sin^2\theta_a\,(\Ric_g)_{aa}\\
=&|A_M|^2-\sum_a\sin^2\theta_a\sum_{i,j}(h_M)_{ija}^2
 +\sum_a\sin^2\theta_a\,(\Ric_g)_{aa}\\
=&|A_M|^2+2\sum_i\sin^2\theta_i\,(\Ric_g)_{ii},
\end{aligned}
\end{equation}
where $\nabla\langle F,Jp_a\rangle=\sin\theta_a\,e_a$ and
$\nabla^2\langle F,Jp_a\rangle(e_i,e_j)=\cos\theta_a\,(h_M)_{ija}$.
In the last equality, we used
$(\Ric_g)_{aa}=-\sum_{i,j}(h_M)_{ija}^2$, by
\eqref{eq:ric-matrices} and the symmetry of $h_M$.
For $b>0$, by Cauchy--Schwarz and
Lemma~\ref{lem:kato-equality-algebra},
\begin{equation}\label{eq:normal-coordinate-kato}
\begin{aligned}
|\nabla b|^2
=&b^{-2}\sum_\ell
 \Bigl(\sum_i\sin\theta_i\cos\theta_i\,(h_M)_{\ell ii}\Bigr)^2\\
\le&\frac n{n+2}
 \frac{\sum_i\sin^2\theta_i\cos^2\theta_i}{b^2}|A_M|^2.
\end{aligned}
\end{equation}

To compare the Ricci term in \eqref{eq:normal-coordinate-bochner}
with $\Pcal$, let $\mathring{\Ric}_g\coloneqq \Ric_g+\frac{|A_M|^2}{n}g$ denote the trace free
Ricci tensor, and define the Weyl tensor by
\[
\begin{aligned}
(W_g)_{ijk\ell}\coloneqq &(\Rm_g)_{ijk\ell}
-\frac{1}{n-2}\bigl(
(\Ric_g)_{ik}g_{j\ell}+(\Ric_g)_{j\ell}g_{ik}
-(\Ric_g)_{i\ell}g_{jk}-(\Ric_g)_{jk}g_{i\ell}\bigr)\\
&+\frac{\operatorname{Scal}_g}{(n-1)(n-2)}
\bigl(g_{ik}g_{j\ell}-g_{i\ell}g_{jk}\bigr).
\end{aligned}
\]
Using $\sum_i(\mathring{\Ric}_g)_{ii}=0$ and
$\sum_i(W_g)_{iji\ell}=0$, we obtain
\begin{equation}\label{eq:curvature-norm-decompositions}
\begin{aligned}
|\Ric_g|^2&=|\mathring{\Ric}_g|^2+\frac1n|A_M|^4,\\
|\Rm_g|^2&=|W_g|^2+\frac4{n-2}|\mathring{\Ric}_g|^2
 +\frac2{n(n-1)}|A_M|^4,
\end{aligned}
\end{equation}
By \eqref{eq:curvature-norm-decompositions} and
\eqref{eq:simons-potential},
\begin{equation}\label{eq:sharp-curvature-decomposition}
\Pcal|A_M|^2=\frac{n+1}{n(n-1)}|A_M|^4
 +\frac{n+2}{n-2}|\mathring{\Ric}_g|^2+|W_g|^2.
\end{equation}

At a point with $A_M\ne0$, normalize $|A_M|^2=1$. For $0<b^2\le n$, we have
\[
\begin{aligned}
&1+2\sum_i\sin^2\theta_i\,(\Ric_g)_{ii}
 -\frac n{n-1}|\nabla b|^2
 +\frac{n-1}{n+1}\Pcal b^2\\
\ge&1-\frac{2b^2}{n}
 -2|\mathring{\Ric}_g|
 \Bigl(\sum_i\Bigl(\sin^2\theta_i-\frac{b^2}{n}\Bigr)^2\Bigr)^{1/2}\\
&-\frac n{n-1}\Bigl[
 \frac n{n+2}\Bigl(1-\frac{b^2}{n}\Bigr)
 -\frac{n}{(n+2)b^2}
 \sum_i\Bigl(\sin^2\theta_i-\frac{b^2}{n}\Bigr)^2\Bigr]\\
&+\frac{(n-1)b^2}{n+1}
 \Bigl(\frac{n+1}{n(n-1)}
 +\frac{n+2}{n-2}|\mathring{\Ric}_g|^2\Bigr)\\
=&\frac{n-2}{(n-1)(n+2)}\Bigl(1-\frac{b^2}{n}\Bigr)
 +\frac{n^2}{(n-1)(n+2)b^2}
 \sum_i\Bigl(\sin^2\theta_i-\frac{b^2}{n}\Bigr)^2\\
&-2|\mathring{\Ric}_g|
 \Bigl(\sum_i\Bigl(\sin^2\theta_i-\frac{b^2}{n}\Bigr)^2\Bigr)^{1/2}
 +\frac{(n-1)(n+2)b^2}{(n+1)(n-2)}|\mathring{\Ric}_g|^2\\
\ge&\frac{n-2}{(n-1)(n+2)}\Bigl(1-\frac{b^2}{n}\Bigr)
 +\frac{1}{(n-1)b^2}
 \sum_i\Bigl(\sin^2\theta_i-\frac{b^2}{n}\Bigr)^2,
\end{aligned}
\]
where we used $\tr_g\Ric_g=-1$, Cauchy--Schwarz,
\eqref{eq:normal-coordinate-kato}, and
\eqref{eq:sharp-curvature-decomposition} in the first inequality,
and Young's inequality in the last.
Thus, without the normalization,
\begin{equation}\label{eq:normal-coordinate-joint-defect}
\begin{aligned}
&|A_M|^2+2\sum_i\sin^2\theta_i\,(\Ric_g)_{ii}
 -\frac n{n-1}|\nabla b|^2
 +\frac{n-1}{n+1}\Pcal b^2\\
\ge&\frac{n-2}{(n-1)(n+2)}
 \Bigl(1-\frac{b^2}{n}\Bigr)|A_M|^2.
\end{aligned}
\end{equation}
At $A_M=0$, both sides of \eqref{eq:normal-coordinate-joint-defect} vanish.

By \eqref{eq:normal-coordinate-bochner},
\eqref{eq:normal-coordinate-joint-defect}, and \eqref{eq:P-estimates},
\begin{equation}\label{eq:normal-coordinate-power-coercivity}
\begin{aligned}
&b^{2-2\alpha}(b^\alpha\Delta(b^\alpha)+\delta\Pcal b^{2\alpha})\\
\ge&\alpha\Bigl(\alpha-\frac{n-2}{n-1}\Bigr)|\nabla b|^2
 +\frac{\alpha(n-2)}{(n-1)(n+2)}
 \Bigl(1-\frac{b^2}{n}\Bigr)|A_M|^2\\
&+\Bigl(\delta-\frac{n-1}{n+1}\alpha\Bigr)\Pcal b^2
\ge c|A_M|^2
\end{aligned}
\end{equation}
on $\{b>0\}$, for some $c=c(n,\delta,\alpha)>0$.
The last inequality in \eqref{eq:normal-coordinate-power-coercivity}
uses $\alpha>(n-2)/(n-1)$, $0<b^2\le n$,
$\delta-\frac{n-1}{n+1}\alpha>0$, and \eqref{eq:P-estimates}.
For $0<k\le k_*$, where
$0<k_*\le\min\{1/8,c/(3\delta n^{1-\alpha/2})\}$,
\[
\begin{aligned}
b^{2-2\alpha}\bigl(b^\alpha\Delta(b^\alpha)+\delta\Pcal b^\alpha(b^\alpha-k)\bigr)
\ge&c|A_M|^2-\delta k\Pcal b^{2-\alpha}
\ge\frac c2|A_M|^2.
\end{aligned}
\]
For $\eta\in\operatorname{Lip}_c(M)$, multiplying by
$b^{\alpha-2}(b^\alpha-k)^+\eta^2$ and integrating over $M$, we obtain
\begin{equation}\label{eq:bounded-tilt-caccioppoli}
\begin{aligned}
&\frac c2\int_{\{b^\alpha>k\}}
 \Bigl(1-\frac{k}{b^\alpha}\Bigr)b^{2\alpha-2}|A_M|^2\eta^2\,d\mu\\
\le&\delta\int_M\Pcal\bigl[(b^\alpha-k)^+\eta\bigr]^2\,d\mu
 +\int_{\{b^\alpha>k\}}(b^\alpha-k)\eta^2\Delta(b^\alpha)\,d\mu\\
\le&\int_M\bigl|\nabla\bigl[(b^\alpha-k)^+\eta\bigr]\bigr|^2\,d\mu
 +\int_{\{b^\alpha>k\}}(b^\alpha-k)\eta^2\Delta(b^\alpha)\,d\mu\\
=&\int_M(b^\alpha-k)_+^2|\nabla\eta|^2\,d\mu,
\end{aligned}
\end{equation}
where we used \eqref{eq:delta-stability} in the second inequality
and integration by parts in the equality.
For $\eta\in C_c^1(B_2)$ and each $\ell$, cover $S\cap\spt\eta$
by balls with $\sum_i r_i^{n-2}\to0$ as $\ell\to\infty$ and take
$\chi_\ell$ from the covering construction in the proof of
Proposition~\ref{thm:eps-reg}.
Since $T$ is calibrated and $\partial T=0$, $|T|$ is stationary.
The cutoff argument for \eqref{eq:singular-curvature-L2} gives
$\int_{F^{-1}(K)}|A_M|^2\,d\mu<\infty$ for $K\Subset B_2$.
By \eqref{eq:normal-coordinate-kato}, the gradient of
$(b^\alpha-k)^+(\eta\circ F)$ is square integrable for $k>0$.
Since $b\le\sqrt n$,
\[
\int_M b^{2\alpha}|\nabla\chi_\ell|^2\,d\mu
\le n^\alpha\int_M|\nabla\chi_\ell|^2\,d\mu
\le C\sum_i r_i^{n-2}\longrightarrow0\quad\text{as }\ell\to\infty.
\]
Thus \eqref{eq:bounded-tilt-caccioppoli} holds with $\eta\circ F$, and
\eqref{eq:MS} holds for $(b^\alpha-k)^+(\eta\circ F)$.
By \eqref{eq:normal-coordinate-kato},
$|\nabla(b^\alpha)|^2\le\alpha^2b^{2\alpha-2}|A_M|^2$ on $\{b>0\}$.
Take $k_\ell=d(1-2^{-\ell-1})$ and $\rho_\ell=1/2+2^{-\ell-1}$,
where $0<d\le k_*$, and let
\[
Y_\ell\coloneqq\int_{F^{-1}(B_{\rho_\ell})}(b^\alpha-k_\ell)_+^2\,d\mu.
\]
Choose $\eta_\ell\in C_c^1(B_{\rho_\ell})$ with
$0\le\eta_\ell\le1$, $\eta_\ell|_{B_{\rho_{\ell+1}}}=1$, and
$|\nabla\eta_\ell|\le C2^\ell$.
By \eqref{eq:bounded-tilt-caccioppoli},
\[
\begin{aligned}
&\int_M\bigl|\nabla\bigl[(b^\alpha-k_{\ell+1})^+(\eta_\ell\circ F)\bigr]\bigr|^2\,d\mu\\
\le&C2^\ell\int_{\{b^\alpha>k_\ell\}}
 \Bigl(1-\frac{k_\ell}{b^\alpha}\Bigr)b^{2\alpha-2}|A_M|^2
 (\eta_\ell\circ F)^2\,d\mu+C4^\ell Y_\ell\\
\le&C8^\ell Y_\ell,
\end{aligned}
\]
where $1-k_\ell/b^\alpha\ge2^{-\ell-2}$ on
$\{b^\alpha>k_{\ell+1}\}$ in the first inequality.
By H\"older's inequality and \eqref{eq:MS},
\[
\begin{aligned}
Y_{\ell+1}
\le&\Bigl(\int_M\bigl[(b^\alpha-k_{\ell+1})^+(\eta_\ell\circ F)\bigr]^{\frac{2n}{n-2}}\,d\mu\Bigr)^{\frac{n-2}{n}}\\
&\qquad\times\bigl|\{b^\alpha>k_{\ell+1}\}\cap F^{-1}(B_{\rho_\ell})\bigr|^{2/n}\\
\le&C\int_M\bigl|\nabla\bigl[(b^\alpha-k_{\ell+1})^+(\eta_\ell\circ F)\bigr]\bigr|^2\,d\mu
 \Bigl(\frac{2^{\ell+2}}d\Bigr)^{4/n}Y_\ell^{2/n}\\
\le&C8^\ell\Bigl(\frac{2^{\ell+2}}d\Bigr)^{4/n}Y_\ell^{1+2/n}
\le CB^\ell d^{-4/n}Y_\ell^{1+2/n},
\end{aligned}
\]
where we used $b^\alpha-k_\ell\ge d2^{-\ell-2}$ on
$\{b^\alpha>k_{\ell+1}\}$ in the second inequality;
$C=C(n,\delta,\alpha)$ and $B=B(n,\delta,\alpha)>1$.
Therefore \eqref{eq:bounded-tilt-estimate} follows by the iteration in
the proof of Proposition~\ref{thm:eps-reg}, or
\cite[Section~4]{Bellettini2025}.

Returning to $F_j$ and using varifold convergence, H\"older's inequality
and $\|T_j\|(B_2)\le\Lambda$, we obtain, for $K\Subset B_2$,
\[
\begin{aligned}
\int_{F_j^{-1}(K)}b_j^{2\alpha}\,d\mu_j
\le&\mu_j(F_j^{-1}(K))^{1-\alpha}
 \Bigl(\int_{F_j^{-1}(K)}b_j^2\,d\mu_j\Bigr)^\alpha
\longrightarrow0\qquad\text{as }j\to\infty.
\end{aligned}
\]
Consequently, by \eqref{eq:bounded-tilt-estimate} on interior balls and
\eqref{eq:stationary-support-mass},
\[
\sup_{F_j^{-1}(K)}b_j
+\sup_{F_j^{-1}(K)}\dist(F_j,P)\longrightarrow0
\quad\text{as }j\to\infty.
\]

To extend across $S_j$, take $r_\perp>0$ and a ball $D\subset P$ such that
$\mathcal Z=D+B_{r_\perp}^{P^\perp}(0)\Subset B_2$. After decreasing $D$,
\[
F_j(M_j)\cap\bigl(\overline D+\partial B_{r_\perp}^{P^\perp}(0)\bigr)
=\varnothing\qquad\forall\,j\ge j_0.
\]
Let $S_{j,P}=D\cap\pi_P(S_j\cap\overline{\mathcal Z})$.
Then $S_{j,P}$ is relatively closed and $\mathcal H^{n-2}(S_{j,P})=0$.
The map $\pi_P\circ F_j$ on
$F_j^{-1}(\mathcal Z\cap\pi_P^{-1}(D\setminus S_{j,P}))$ is a proper finite
covering of $D\setminus S_{j,P}$.
For a smooth loop $\gamma:\Sph^1\to D\setminus S_{j,P}$,
choose a smooth $\Gamma:\overline{\mathbb D}\to D$ with
$\Gamma|_{\Sph^1}=\gamma$. Cover $S_{j,P}$ by balls $B_{r_i}$ and
$\Gamma(\overline{\mathbb D})$ by at most $Cr_i^{-2}$ balls of radius $r_i$.
Then
\[
\mathcal H^n\bigl(S_{j,P}-\Gamma(\overline{\mathbb D})\bigr)
\le C\sum_i r_i^{-2}r_i^n
=C\sum_i r_i^{n-2}\longrightarrow0\quad\text{as }\max_i r_i\downarrow0.
\]
Choose an arbitrarily small $v\in P$ such that
$\bigl(\Gamma(\overline{\mathbb D})+v\bigr)\cap S_{j,P}=\varnothing$.
For sufficiently small $v$, $\gamma(z)+tv\in D\setminus S_{j,P}$ for
$z\in\Sph^1$ and $0\le t\le1$.
For $x,y\in D\setminus S_{j,P}$, translating the line segment from
$x$ to $y$, we also obtain paths
in $D\setminus S_{j,P}$ from $x$ to $y$ with lengths approaching $|x-y|$.
Thus $D\setminus S_{j,P}$ is simply connected.

For the lifted graph functions $u_{j,l}:D\setminus S_{j,P}\to P^\perp$,
$\max_l\|Du_{j,l}\|_\infty\to0$ as $j\to\infty$. Integrating
$Du_{j,l}$ along paths arbitrarily close to straight segments, we have
\[
|u_{j,l}(x)-u_{j,l}(y)|
\le\|Du_{j,l}\|_\infty|x-y|\qquad\forall\,x,y\in D\setminus S_{j,P}.
\]
Consequently, each $u_{j,l}$ extends to a Lipschitz map on $D$, with
$\partial\llbracket\operatorname{graph}u_{j,l}\rrbracket=0$ in $D+P^\perp$.
Since $\operatorname{graph}u_{j,l}$ is calibrated almost everywhere,
$|\llbracket\operatorname{graph}u_{j,l}\rrbracket|$ is stationary in $D+P^\perp$.
For interior balls centered on $\operatorname{graph}u_{j,l}$,
\[
\frac{\mathcal H^n(\operatorname{graph}u_{j,l}\cap B_\rho(z))}
 {\omega_n\rho^n}
\le(1+\|Du_{j,l}\|_\infty^2)^{n/2}=1+o(1).
\]
Applying Allard's regularity theorem
\cite[Chapter~5, Theorem~5.2]{Simon2014GMT} and then Schauder estimates
to \eqref{eq:minimal-graph-equations}, we obtain smooth convergence of
the extended graphs to $P$, with total multiplicity
$q$ by varifold convergence. Capacity cutoffs pass
\eqref{eq:delta-stability} to the domains obtained by extending $F_j$
across $S_j$.
\end{proof}

To start the density induction in the proof of
Theorem~\ref{thm:three-dimensional-null-set-compactness}, we prove
smooth embeddedness at points of density less than two.

\begin{lemma}\label{lem:low-density-SLag-regularity}
Let $T$ be a special Lagrangian integral current with $\partial T=0$
in an open subset of $\C^n$, $n\ge3$. Suppose that \eqref{eq:delta-stability} holds on the
multiplicity one embedded regular part of $T$ and
$\delta>\delta_{\rm cone}(n)$. Then $T$ is smoothly embedded with
multiplicity one near each $x\in\spt T$ satisfying
\[
\lim_{r\downarrow0}\frac{\|T\|(B_r(x))}{\omega_n r^n}<2.
\]
\end{lemma}

\begin{proof}
Suppose calibrated integral currents $T_j$ satisfy $\partial T_j=0$,
$T_j\rightharpoonup T$ and $|T_j|\to|T|$ as $j\to\infty$, and \eqref{eq:delta-stability}
holds on the multiplicity one regular part of each $T_j$.
Applying Allard's regularity theorem
\cite[Chapter~5, Theorem~5.2]{Simon2014GMT} and Schauder estimates to
\eqref{eq:minimal-graph-equations}, we obtain smooth convergence of
the local graph functions on compact subsets of the multiplicity one
regular part of $T$. Let $\Omega$ be a relatively compact domain
in this regular part and $\Omega_j$ the corresponding graph domains
of $T_j$. For $\phi\in C_c^1(\Omega)$, let $\phi_j$ be the pullback
of $\phi$ under normal projection $\Omega_j\to\Omega$. By smooth
convergence and \eqref{eq:P-estimates},
\begin{equation}\label{eq:multiplicity-one-limit-stability}
\begin{aligned}
\int_\Omega\bigl(|\nabla\phi|^2-\delta\Pcal\phi^2\bigr)\,d\mu
=&\lim_{j\to\infty}\int_{\Omega_j}
 \bigl(|\nabla\phi_j|^2-\delta\Pcal_j\phi_j^2\bigr)\,d\mu_j
\ge0.
\end{aligned}
\end{equation}

We claim that, for $1\le m\le n$, a special Lagrangian integral cone
$C$ in $\C^m$ with $\partial C=0$, density less than two at the vertex,
and stable multiplicity one regular part is a multiplicity one plane.
For $m=1$, the calibrated rays lie in one line with equal multiplicities
by $\partial C=0$, and both multiplicities are one since $\Theta^1(\|C\|,0)<2$.

Now let $m\ge2$. At a nonzero point $z$ of an
$m$-dimensional cone $C$,
a tangent cone $D$ is invariant along $\R z$. By the Lagrangian
condition and slicing,
\[
D=\llbracket\R z\rrbracket\times D',\qquad
D'\subset\bigl(\R z\oplus J\R z\bigr)^\perp\simeq\C^{m-1}.
\]
Passing $\partial C=0$ to $D$ and slicing
$D=\llbracket\R z\rrbracket\times D'$,
we have $\partial D=\partial D'=0$.
Since $D$ is a tangent cone of $C$ at $z$ and
$D=\llbracket\R z\rrbracket\times D'$, by
\eqref{eq:stationary-monotonicity} and the conical invariance of $C$,
\[
\Theta^{m-1}(\|D'\|,0)=\Theta^m(\|D\|,0)
=\Theta^m(\|C\|,z)\le\Theta^m(\|C\|,0)<2.
\]
By \eqref{eq:multiplicity-one-limit-stability}, $D$ is stable on
the multiplicity one regular part. Since $\Pcal_D=\Pcal_{D'}$,
substituting $\chi(t/R)\phi$ into stability on $D$ and dividing by
$R\int_\R\chi^2\,dt$, for $0\ne\chi\in C_c^1(\R)$ and $\phi$
compactly supported in the multiplicity one regular part of $D'$, we have
\[
\begin{aligned}
\delta\int_{D'}\Pcal_{D'}\phi^2\,d\mu_{D'}
\le&\int_{D'}|\nabla\phi|^2\,d\mu_{D'}
+\frac{\int_\R(\chi')^2\,dt}{R^2\int_\R\chi^2\,dt}
 \int_{D'}\phi^2\,d\mu_{D'}\\
\longrightarrow&\int_{D'}|\nabla\phi|^2\,d\mu_{D'}
\quad\text{as }R\to\infty.
\end{aligned}
\]
Thus $D'$ is $\delta$-Simons stable. By induction, $D'$ and hence $D$ are
multiplicity one planes, so $C$ is smooth on $C\setminus\{0\}$ by
\cite[Chapter~5, Theorem~5.2]{Simon2014GMT}.

For $m=2$, $C\cap\Sph^3$ is a finite union of closed geodesics,
so $C$ is a finite planar union. If $m\ge3$, then
$\delta>\delta_{\rm cone}(n)\ge\delta_{\rm cone}(m)$, and $C$ is again a
finite planar union by Theorem~\ref{thm:intro-cone}.
Since $\Theta^m(\|C\|,0)\in\mathbb N\cap[1,2)=\{1\}$, the planar
union consists of one plane with multiplicity one.

At a point $x$ in the statement, every tangent cone $C$ of $T$ has
density less than two and satisfies \eqref{eq:delta-stability} on the
multiplicity one regular part of $C$, by the graph convergence proved
at the start of this proof.
By the claim, every tangent cone of $T$ at $x$ is a multiplicity one
plane. Applying \cite[Chapter~5, Theorem~5.2]{Simon2014GMT}, we conclude
that $T$ is smoothly embedded with
multiplicity one near $x$.
\end{proof}

For the proof of Proposition~\ref{prop:three-dimensional-planar-induction-step},
we next estimate the $C^2$ norm of each graph over one of finitely many planes.

\begin{lemma}\label{lem:finite-annular-graph-matching}
Let $Q\ge1$ and $A_0\Subset A_1\Subset A_2\subset\mathbb R^6$.
Suppose $\Sigma_j$ is the sum of $Q$ smooth minimal graphs defined across
$\overline A_2$, converging smoothly to linear planes of dimension three.
Assume that the $A_i$ are spherical annuli centered at zero, or
concentric balls whose center lies in every limiting plane of $\Sigma_j$.
Let $H_{j,1},\ldots,H_{j,p}$, $p\le Q$, be distinct linear
planes of dimension three converging to planes in the limit of
$\Sigma_j$, and assume
\[
E_j^2\coloneqq\int_{\Sigma_j\cap A_2}
 \dist\Bigl(x,\bigcup_{i=1}^pH_{j,i}\Bigr)^2\,d\mu_j(x)
\longrightarrow0\quad\text{as }j\to\infty.
\]
Then the $l$th graph in $\Sigma_j$ has a normal graph representation over some
$H_{j,i_j(l)}$ on $A_0$, with
\begin{equation}\label{eq:annular-finite-own-height}
\sum_{l=1}^Q\|u_{j,l}\|_{C^2(H_{j,i_j(l)}\cap A_0)}^2
\le CE_j^2.
\end{equation}
For
\begin{equation}\label{eq:annular-finite-bilateral-integrals}
B_j^2\coloneqq E_j^2+\sum_{i=1}^p\int_{H_{j,i}\cap A_2}
 \dist(y,\spt\Sigma_j)^2\,d\mathcal H^3(y),
\end{equation}
we also have
$\max_i\min_l\|\pi_{H_{j,i}}-\pi_{H_{j,i_j(l)}}\|\le CB_j$.
Here $C$ depends only on $Q,A_0,A_1,A_2$.
\end{lemma}

\begin{proof}
We induct on $p$, using intermediate regions between $A_0$ and
$A_2$ for Schauder estimates in \eqref{eq:minimal-graph-equations} and
changes of graph projection.
For $p=1$, the area formula and the $L^2$ to $C^2$ interior estimate for
\eqref{eq:minimal-graph-equations} imply
\[
\sum_l\|u_{j,l}\|_{C^2(H_{j,1}\cap A_0)}^2
\le C\sum_l\int_{H_{j,1}\cap A_1}|u_{j,l}|^2\,d\mathcal H^3
\le CE_j^2.
\]

Suppose \eqref{eq:annular-finite-own-height} holds for $p-1$
planes, and consider
\[
a_j\coloneqq\min_{1\le i<k\le p}\|\pi_{H_{j,i}}-\pi_{H_{j,k}}\|>0.
\]
If $a_j\to0$ as $j\to\infty$, relabel so that
$\|\pi_{H_{j,p-1}}-\pi_{H_{j,p}}\|=a_j$ and delete $H_{j,p}$.
By the triangle inequality,
\[
\begin{aligned}
&\Bigl(\int_{\Sigma_j\cap A_2}
 \dist\Bigl(x,\bigcup_{i=1}^{p-1}H_{j,i}\Bigr)^2\,d\mu_j(x)\Bigr)^{1/2}\\
\le&E_j+a_j\Bigl(\int_{\Sigma_j\cap A_2}|x|^2\,d\mu_j(x)\Bigr)^{1/2}
\le C(E_j+a_j).
\end{aligned}
\]
If $E_j/a_j\to0$ as $j\to\infty$, the graph functions $u_{j,l}$ over
$H_{j,i_j(l)}$ are $O(a_j)$ in $C^2$ on intermediate regions by induction.
After orthogonal identification of $H_{j,i_j(l)}$ with
$\lim_jH_{j,i_j(l)}$ and passage to a subsequence,
\[
a_j^{-1}u_{j,l}\longrightarrow v_l\quad\text{in }C^1_{\rm loc}\text{ as }j\to\infty,
\qquad \Delta v_l=0.
\]

For fixed $l$, express each $H_{j,k}$ having the same limit as
$H_{j,i_j(l)}$ as a linear graph over $H_{j,i_j(l)}$.
If the linear graph map divided by $a_j$ stays bounded, pass to a linear
limit. Otherwise, normalize to norm one before passing to a limit;
outside the kernel of the normalized linear limit for $H_{j,k}$, the distance to
$H_{j,k}$ divided by $a_j$ tends to infinity.
Planes $H_{j,k}$ whose limits differ from $\lim_jH_{j,i_j(l)}$
are excluded away from the intersections of
$\lim_jH_{j,k}$ and $\lim_jH_{j,i_j(l)}$.
The kernels and plane intersections have zero three dimensional measure.
Let $b_k$ denote the finite linear limits. By $E_j^2/a_j^2\to0$ as
$j\to\infty$ and unique continuation on the connected ball or annulus,
\[
\min_k|v_l-b_k|=0\quad\text{a.e.}
\quad\Longrightarrow\quad
v_l=b_k\quad\text{for one fixed }k.
\]
Choose $i$ for which the linear graph function of $H_{j,i}$,
divided by $a_j$, converges to $b_k=v_l$. Over $H_{j,i}$,
$u_{j,l}=o(a_j)$ in $C^1$.
If $\inf_j a_j>0$, smooth convergence also implies
$\|u_{j,l}\|_{C^1}/a_j\to0$ as $j\to\infty$.
In both cases, we claim
\begin{equation}\label{eq:annular-height-comparison}
\sum_l\int_{H_{j,i_j(l)}\cap A_1}|u_{j,l}|^2\,d\mathcal H^3
\le C\int_{\Sigma_j\cap A_2}
 \dist\Bigl(x,\bigcup_iH_{j,i}\Bigr)^2\,d\mu_j(x)
\le CE_j^2.
\end{equation}
Indeed, for each graph $u_{j,l}$ over $H_{j,i}$ choose a unit
$e_{j,i}\in H_{j,i}$ satisfying
\[
|\pi_{H_{j,k}^\perp}e_{j,i}|
\ge c_p\|\pi_{H_{j,i}}-\pi_{H_{j,k}}\|
\qquad(k\ne i).
\]
Each excluded spherical band has measure at most $Cc_p$; choose
$c_p>0$ with $C(p-1)c_p<\mathcal H^2(\Sph^2)$. On an enlarged convex
graph disk consider
\[
\Phi_m(x)\coloneqq x+\frac{Km}{a_j}|u_{j,l}(x)|e_{j,i},
\qquad m=0,\ldots,p-1.
\]
For $K>10/c_p$, consecutive projected points are separated by at least
$Kc_p|u_{j,l}(x)|$, so
\[
\#\bigl\{m\in\{0,\ldots,p-1\}:
 |\pi_{H_{j,k}^\perp}\Phi_m(x)|<4|u_{j,l}(x)|\bigr\}\le1
\qquad(k\ne i).
\]
By the mean value theorem and $u_{j,l}=o(a_j)$ in $C^1$,
for all large $j$,
\[
\frac34|u_{j,l}(x)|\le |u_{j,l}(\Phi_m(x))|
\le\frac54|u_{j,l}(x)|.
\]
Differentiating $\Phi_m$ and using $\|Du_{j,l}\|_\infty/a_j\to0$,
we also have $3/4\le\operatorname{Jac}\Phi_m\le5/4$.
The maps $\Phi_m$ are injective and stay in the enlarged disk. For each $x$,
\[
\exists m\in\{0,\ldots,p-1\}:\quad
\dist\Bigl(\Phi_m(x)+u_{j,l}(\Phi_m(x)),\bigcup_kH_{j,k}\Bigr)
=|u_{j,l}(\Phi_m(x))|.
\]
Partitioning by an admissible $m$, applying the area formula and summing
over at most $C(Q)$ disks covering $H_{j,i_j(l)}\cap A_1$, we obtain
\eqref{eq:annular-height-comparison}, and hence
\eqref{eq:annular-finite-own-height} on $A_0$ by the $L^2$ to $C^2$
interior estimate for \eqref{eq:minimal-graph-equations}.

If no constants $\beta,\varepsilon>0$, depending only on
$p,Q,A_0,A_1,A_2$, make \eqref{eq:annular-finite-own-height} valid for
$E_j\le\min\{\beta a_j,\varepsilon\}$, choose a sequence for which
$E_j\to0$ and $E_j/a_j\to0$ as $j\to\infty$ and
\eqref{eq:annular-finite-own-height} fails.
There exists a subsequence along which $H_{j,1},\ldots,H_{j,p}$ converge, and
either $a_j\to0$ or $\inf_j a_j>0$; both cases were proved using
\eqref{eq:annular-height-comparison} and the graph equation.
Thus such $\beta,\varepsilon$ exist. In the remaining case $a_j\le E_j/\beta$, apply the induction hypothesis
after deleting one closest plane:
\[
\sum_l\|u_{j,l}\|_{C^2(H_{j,i_j(l)}\cap A_0)}^2
\le C(E_j+a_j)^2\le C_\beta E_j^2.
\]
Every plane after deletion belongs to $\{H_{j,1},\ldots,H_{j,p}\}$,
so the induction on $p$ is complete. For varying limiting planes, pass
to a subsequence, combine equal limits, and represent the $Q$ graphs
of $\Sigma_j$ over $\lim_jH_{j,i}$ on smaller fixed regions.
Thus $C$ depends only on $Q,A_0,A_1,A_2$.

If $i\notin i_j(\{1,\ldots,Q\})$, take a fixed disk $U$ in the
limit of $H_{j,i}$, away from all other limit planes.
The radius of $U$ depends only on $Q,A_0,A_1,A_2$, by the exclusion of
the $p-1$ spherical bands used to choose $e_{j,i}$.
The first index $j$ for which the graphs of $\Sigma_j$ over
$U$ are separated may depend on the angles between the limits of
$H_{j,1},\ldots,H_{j,p}$. By \eqref{eq:annular-finite-own-height} and
the integrals over $H_{j,i}\cap A_2$ in
\eqref{eq:annular-finite-bilateral-integrals},
\[
\int_U\min_l|(A_{j,i}-A_{j,i_j(l)})x|^2\,dx\le CB_j^2,
\]
where $H_{j,i}=\operatorname{graph}A_{j,i}$ and
$H_{j,i_j(l)}=\operatorname{graph}A_{j,i_j(l)}$ over $\lim_jH_{j,i}$.
On a fixed open
disk, for any finite family of linear maps $M_k$,
\begin{equation}\label{eq:finite-linear-map-integral}
\int_U\min_k|M_kp|^2\,dp\ge c\min_k\|M_k\|^2.
\end{equation}
Indeed, choose $c_1>0$, depending only on $U$ and the number of maps,
so that
\[
\sum_k\mathcal H^3\{p\in U:|M_kp|<c_1\|M_k\|\}
\le\tfrac12\mathcal H^3(U).
\]
Integrating $\min_k|M_kp|^2\ge c_1^2\min_k\|M_k\|^2$ over
$U\setminus\bigcup_k\{p:|M_kp|<c_1\|M_k\|\}$ proves
\eqref{eq:finite-linear-map-integral}. Applying
\eqref{eq:finite-linear-map-integral} with
$M_k=A_{j,i}-A_{j,i_j(k)}$, we obtain
$\min_l\|\pi_{H_{j,i}}-\pi_{H_{j,i_j(l)}}\|\le CB_j$.
\end{proof}

We now prove the induction step used in
Theorem~\ref{thm:three-dimensional-null-set-compactness} when the
varifold limit is a finite planar union.

\begin{proposition}\label{prop:three-dimensional-planar-induction-step}
Fix an integer $Q\ge2$ and $\delta>1/4$.
Let $F_j,S_j,T_j$ satisfy the hypotheses of
Proposition~\ref{prop:null-set-plane-regularity} with $n=3$, without
assuming $|T_j|\to q|P|$ as $j\to\infty$. Assume removable extension and smooth convergence below
limiting density $Q$ for every sequence with these hypotheses. Suppose
\[
|T_j|\longrightarrow\sum_{\alpha=1}^Nq_\alpha|P_\alpha|
\quad\text{as }j\to\infty,\qquad q_\alpha\in\mathbb N,
\qquad\sum_\alpha q_\alpha=Q,
\]
where the $P_\alpha$ are distinct linear special Lagrangian
planes of dimension three with a common phase. Then, after extension across $S_j$,
the immersions $F_j$ converge smoothly on compact subsets of $B_2$.
\end{proposition}

\begin{proof}
For $N=1$, apply Proposition~\ref{prop:null-set-plane-regularity}.
If $N\ge2$, the common intersection is either $\{0\}$ or a line,
since distinct planes of dimension three with the same phase intersect in dimension
at most one.
Choose disjoint neighborhoods $\mathcal U_\alpha$ of $P_\alpha$ in the
metric $\|\pi_P-\pi_{P'}\|$. It suffices to find $\rho>0$ and $j_0$ such that
\begin{equation}\label{eq:finite-Q-tangent-neighborhood}
dF_j(T_xM_j)\in\bigcup_\alpha\mathcal U_\alpha
\qquad\forall\,j\ge j_0,\ \forall\,x\in F_j^{-1}(B_\rho).
\end{equation}
Indeed, under \eqref{eq:finite-Q-tangent-neighborhood}, the subsets with
tangent planes in $\mathcal U_\alpha$ are open and closed in the original
immersion domain over $B_\rho\setminus S_j$. Restriction preserves
properness, absence of boundary and stability. By $\mathcal H^2(S_j)=0$
and \eqref{eq:stationary-monotonicity}, the currents of $F_j$ restricted to
$dF_j(T_xM_j)\in\mathcal U_\alpha$ extend across $S_j$ to integral
currents $T_{j,\alpha}$ on $B_\rho$ satisfying
\[
\partial T_{j,\alpha}=0\quad\text{in }B_\rho,\qquad
T_{j,\alpha}\rightharpoonup q_\alpha\llbracket P_\alpha\rrbracket
\quad\text{in }B_\rho\text{ as }j\to\infty.
\]
The restricted immersions therefore extend across $S_j$ and converge
smoothly to $P_\alpha$ by Proposition~\ref{prop:null-set-plane-regularity}.
For a varifold $V$ of dimension three, a finite family of linear planes of dimension three
$H=(H_1,\ldots,H_N)$, and $s>0$, define
\begin{equation}\label{eq:finite-plane-integrals}
\begin{aligned}
\mathcal E(V,H;s)^2\coloneqq s^{-5}\Bigl(&
\int_{B_s}\dist\Bigl(x,\bigcup_l H_l\Bigr)^2\,d\|V\|(x)\\
&+\sum_{l=1}^N\int_{H_l\cap B_s}
 \dist(p,\spt V)^2\,d\mathcal H^3(p)\Bigr).
\end{aligned}
\end{equation}

When a line $L$ is contained in every $H_l$, also define, for $\rho>0$,
\begin{equation}\label{eq:truncated-plane-integrals}
\begin{aligned}
\mathcal Q_L(T,H;\rho)^2\coloneqq\rho^{-5}\Bigl[&
\int_{B_\rho}\dist\Bigl(x,\bigcup_{l=1}^NH_l\Bigr)^2\,d\|T\|(x)\\
&+\sum_{l=1}^N\int_{H_l\cap B_{\rho/2}\cap\{\dist(p,L)>\rho/16\}}
 \dist(p,\spt T)^2\,d\mathcal H^3(p)\Bigr].
\end{aligned}
\end{equation}

Consider points $z_j$ of density at least $Q$ in a smaller closed ball. Every limit $z$ of these points lies in
$\bigcap_\alpha P_\alpha$: otherwise $\sum_\alpha q_\alpha|P_\alpha|$
has density at most $Q-1$ at $z$, contradicting upper semicontinuity.
For a fixed small $r_0>0$,
\[
\bigl(x\mapsto(x-z_j)/r_0\bigr)_\#|T_j|
\longrightarrow\sum_\alpha q_\alpha|P_\alpha|
\quad\text{as }j\to\infty
\]
on $B_2$. By \eqref{eq:stationary-support-mass},
\begin{equation}\label{eq:finite-Q-initial-comparison}
\mathcal E\bigl((\tau_{-z_j})_\#T_j,(P_1,\ldots,P_N);2r_0\bigr)
\longrightarrow0\qquad\text{as }j\to\infty.
\end{equation}
Convergence in \eqref{eq:finite-Q-initial-comparison} is uniform in
$z_j$, since a contrary sequence has a convergent subsequence of
centers and still converges to $\sum_\alpha q_\alpha|P_\alpha|$.

Translate by $-z_j$ and dilate by $r_0^{-1}$, keeping the notation $T_j$. Fix a varifold neighborhood on $B_2$ of
$\sum_\alpha q_\alpha|P_\alpha|$. Shrink $\mathcal U_\alpha$ so that every
planar cone with $q_\alpha$ planes in each $\overline{\mathcal U_\alpha}$,
counted with multiplicity, belongs to this neighborhood.
For all large $j$, $|T_j|$ belongs to the neighborhood and
\[
\Theta^3(\|T_j\|,0)\ge Q,
\qquad \|T_j\|(B_1)<(Q+1/2)\omega_3.
\]
Every family of planes used below has total multiplicity $Q$, and
all planes belong to $\bigcup_\alpha\mathcal U_\alpha$.
Let $B_j=\mathcal E(T_j,H_j;1)$, or
$B_j=\mathcal Q_{L_j}(T_j,H_j;1)$ when the distinct planes meet exactly
in $L_j$. At $B_j=0$, restrict to the constant tangent plane components
and apply Proposition~\ref{prop:null-set-plane-regularity}.

For $B_j>0$ with $B_j\to0$ as $j\to\infty$, we claim that, after a subsequence, there
exist families $H'_j$ of total multiplicity $Q$ and a scale $\theta$
from a fixed finite collection in $(0,1/4)$ such that
\begin{equation}\label{eq:finite-Q-one-scale}
\mathcal E(T_j,H'_j;\theta)\le\tfrac12 B_j
\quad\text{or}\quad
\mathcal Q_{L'_j}(T_j,H'_j;\theta)\le\tfrac12 B_j.
\end{equation}
We choose $H'_j$ so that the projections also satisfy
\begin{equation}\label{eq:finite-Q-one-scale-projections}
\operatorname{dist}_{\mathcal H}
 \bigl(\{\pi_{H'_{j,l}}\}_l,\{\pi_{H_{j,l}}\}_l\bigr)\le CB_j.
\end{equation}
The multiplicities of the planes in $H'_j$ are the degrees of the
normal projections. In the $\mathcal Q_{L'_j}$ alternative of
\eqref{eq:finite-Q-one-scale},
each plane $H'_{j,l}$ contains $L'_j$.

Suppose first that $\bigcap_\alpha P_\alpha=\{0\}$.
Choose $\mathcal U_\alpha$ small enough that
$\bigcap_{\alpha=1}^N H_\alpha=\{0\}$ whenever
$H_\alpha\in\mathcal U_\alpha$ for $\alpha=1,\ldots,N$.
Repeat each $H_{j,i}$ according to the prescribed multiplicity, so that
both graph degrees and reverse integrals
total $Q$.

When $B_j=\mathcal E(T_j,H_j;1)\to0$ as $j\to\infty$, the support of
any current limit of $T_j$ in $B_1$ is
$B_1\cap\bigcup_l\lim_jH_{j,l}$, by \eqref{eq:finite-plane-integrals}.
By constancy on $\lim_jH_{j,l}$, $\Theta^3(\|T_j\|,0)\ge Q$ and
$\|T_j\|(B_1)<(Q+1/2)\omega_3$, the current limit has total multiplicity
$Q$ and density at most $Q-1$ away from zero.
Hence, by upper semicontinuity, for $0<s<1$,
\[
\sup\bigl\{|x|:x\in\overline B_s,\quad
 \Theta^3(\|T_j\|,x)\ge Q\bigr\}\longrightarrow0
\qquad\text{as }j\to\infty.
\]

On a fixed annulus excluding zero, apply
the assumed compactness below density $Q$ to the constant sequences
$T_j,T_j,\ldots$ and extend $F_j$ across $S_j$ on the annulus.
By Proposition~\ref{thm:flat-union}\textup{(i)}, the normal projections
of $F_j$ onto $P_\alpha$ form covers of total degree $Q$. Since a spherical annulus in a
plane of dimension three is simply connected, we obtain $Q$ individual graphs.
Label the $Q$ graphs of $F_j$ consistently on a nested exhaustion
of $B_1\setminus\{0\}$.

On a fixed annulus, choose $i_j(l)$ as in
Lemma~\ref{lem:finite-annular-graph-matching} so that the $l$th graph
has the form $x\mapsto x+u_{j,l}(x)$ over $H_{j,i_j(l)}$, with
$\|u_{j,l}\|_{C^2}\le CB_j$.
For a fixed regular value $b$ in this annulus, integrating
$\operatorname{div}_{M_j}F_j^\top=3$ over $\{|F_j|<b\}$ and using
the area formula for the graphs $u_{j,l}$, we obtain
\[
\begin{aligned}
\bigl|3b^{-3}\|T_j\|(B_b)-3Q\omega_3\bigr|
=&\Bigl|b^{-2}\int_{\{|F_j|=b\}}|\nabla|F_j||\,d\sigma_j
 -3Q\omega_3\Bigr|\\
\le&C\sum_{l=1}^Q\|u_{j,l}\|_{C^1}^2\le CB_j^2,
\end{aligned}
\]
where $d\sigma_j$ is the induced area measure on level sets of $|F_j|$.
Since $\Theta^3(\|T_j\|,0)\ge Q$, \eqref{eq:stationary-monotonicity} then gives
\begin{equation}\label{eq:four-plane-isolated-radial-action}
\int_{B_b}\frac{|x^{\perp T_xT_j}|^2}{|x|^5}\,d\|T_j\|(x)
\le CB_j^2.
\end{equation}

On the immersion domain of $F_j$, let $r=|F_j|$.
On $\{r>0,\ |\nabla r|>0\}$, define $\omega=F_j/r$ and
\[
X\coloneqq |\nabla r|^{-2}\pi_{dF_j(TM_j)}\omega-\omega,
\qquad d\nu_s\coloneqq e^{-2s}|\nabla r|\,d\sigma_j,
\qquad a(s)\coloneqq\int_{\{r=e^s\}}|X|^2\,d\nu_s,
\]
For $\eta\in C_c^1(( -\infty,\log b))$ and $\Phi\in C^1(\Sph^5)$,
using $r^{-3}\eta(\log r)\Phi(\omega)F_j$ in $\delta|T_j|=0$,
we obtain
\begin{equation}\label{eq:three-plane-angular-balance}
\frac d{ds}\int_{\{r=e^s\}}\Phi(\omega)\,d\nu_s
=\int_{\{r=e^s\}}
 \bigl(D\Phi(\omega)X+3\Phi(\omega)|X|^2\bigr)\,d\nu_s.
\end{equation}
The identity holds distributionally in $s$.
By local analyticity and $\Delta r^2=6$, the set $\{|\nabla r|=0\}$
has zero three dimensional measure. The coarea formula and
\eqref{eq:four-plane-isolated-radial-action} give
\begin{equation}\label{eq:angular-radial-integral}
\int_{-\infty}^{\log b}a(s)\,ds
=\int_{B_b}\frac{|x^{\perp T_xT_j}|^2}{|x|^5}\,d\|T_j\|(x)
\le CB_j^2.
\end{equation}
Integrating $\operatorname{div}_{M_j}F_j^\top=3$ over $\{r<e^s\}$ and
using $\Theta^3(\|T_j\|,0)\ge Q$, we also have
\[
\nu_s(\{r=e^s\})=3e^{-3s}\|T_j\|(B_{e^s})\ge3Q\omega_3.
\]

For almost every $0<\varrho<b$, let $K_\varrho$ be the closed angular
support of $\nu_{\log\varrho}$ and
$m_\varrho(s)\coloneqq\int_{\{r=e^s\}}\dist(\omega,K_\varrho)^2\,d\nu_s$.
To substitute $\Phi(\omega)=\dist(\omega,K_\varrho)^2$ into
\eqref{eq:three-plane-angular-balance}, take the smooth approximations
\[
\Phi_\varepsilon(\xi)\coloneqq-\varepsilon\log
\frac{\displaystyle\int_{\{r=\varrho\}}
 e^{-|\xi-\omega|^2/\varepsilon}\,d\nu_{\log\varrho}}
 {\displaystyle\nu_{\log\varrho}(\{r=\varrho\})}.
\]
Then $\Phi_\varepsilon\to\dist(\cdot,K_\varrho)^2$ uniformly as $\varepsilon\downarrow0$,
$0\le\Phi_\varepsilon\le4$ and
$|D\Phi_\varepsilon|^2\le4\Phi_\varepsilon$.
For $m_\varepsilon(s)\coloneqq
\int_{\{r=e^s\}}\Phi_\varepsilon(\omega)\,d\nu_s$, \eqref{eq:three-plane-angular-balance} implies
\[
m_\varepsilon'\le2\sqrt{m_\varepsilon a}+12a
\le\frac{m_\varepsilon}{1+\log(b/\varrho)}
 +(13+\log(b/\varrho))a.
\]
Multiplying the differential inequality for $m_\varepsilon$ by
$e^{-s/(1+\log(b/\varrho))}$, integrating over
$[\log\varrho,\log b]$, and letting $\varepsilon\downarrow0$, we obtain,
by $m_\varrho(\log\varrho)=0$ and \eqref{eq:angular-radial-integral},
\begin{equation}\label{eq:angular-section-distance}
m_\varrho(\log b)
\le C(1+\log(b/\varrho))\int_{\log\varrho}^{\log b}a(s)\,ds
\le C(1+\log(b/\varrho))B_j^2.
\end{equation}

By $\|u_{j,l}\|_{C^2}\le CB_j$ from
Lemma~\ref{lem:finite-annular-graph-matching} and the area formula for
$x\mapsto\pi_{H_{j,i_j(l)}}x/|\pi_{H_{j,i_j(l)}}x|$ on
$\operatorname{graph}u_{j,l}\cap\partial B_b$,
\[
\begin{aligned}
\sum_l\int_{H_{j,i_j(l)}\cap\Sph^5}
 \dist(\omega,K_\varrho)^2\,d\mathcal H^2(\omega)
\le&2b^{-2}\int_{\{r=b\}}
 \dist(F_j/b,K_\varrho)^2|\nabla r|\,d\sigma_j+CB_j^2\\
=&2m_\varrho(\log b)+CB_j^2.
\end{aligned}
\]
Since $\varrho K_\varrho\subset\spt T_j$, integrating in polar
coordinates and using \eqref{eq:angular-section-distance}, we obtain
\begin{equation}\label{eq:isolated-plane-to-current-height}
\sum_l\int_{H_{j,i_j(l)}\cap B_t}\dist(p,\spt T_j)^2\,d\mathcal H^3(p)
\le CB_j^2\int_0^t\varrho^4(1+\log(b/\varrho))\,d\varrho.
\end{equation}
For each $i$, choose $k\in\{i_j(1),\ldots,i_j(Q)\}$ with
$\|\pi_{H_{j,i}}-\pi_{H_{j,k}}\|\le CB_j$, by
Lemma~\ref{lem:finite-annular-graph-matching}.
Rotating $H_{j,k}$ to $H_{j,i}$, we have
\begin{equation}\label{eq:finite-Q-plane-change}
\begin{aligned}
&\int_{H_{j,i}\cap B_t}\dist(p,\spt T_j)^2\,d\mathcal H^3(p) 
\le 2\int_{H_{j,k}\cap B_t}\dist(p,\spt T_j)^2\,d\mathcal H^3(p)
 +CB_j^2t^5.
\end{aligned}
\end{equation}
In the sum over $i$, each integral over $H_{j,k}$ occurs at most $Q$ times.

Substituting
$\Phi(\omega)=\dist(\omega,\bigcup_l(H_{j,l}\cap\Sph^5))^2$ into
\eqref{eq:three-plane-angular-balance}, with
$m(s)\coloneqq\int_{\{r=e^s\}}\Phi(\omega)\,d\nu_s$, we have
$-m'\le2\sqrt{ma}$ and $m(\log b)\le CB_j^2$; hence
\[
m(s)\le CB_j^2(1+\log b-s).
\]
By the coarea formula, $|\nabla r|^{-2}=1+|X|^2$, and
\eqref{eq:angular-radial-integral},
\begin{equation}\label{eq:isolated-current-to-plane-height}
\begin{aligned}
\int_{B_t}\dist\Bigl(x,\bigcup_lH_{j,l}\Bigr)^2\,d\|T_j\|(x)
\le&\int_{-\infty}^{\log t}e^{5s}(m(s)+4a(s))\,ds\\
\le&CB_j^2t^5(1+\log(b/t))\qquad\forall\,t\in(0,b/2).
\end{aligned}
\end{equation}
Combining \eqref{eq:isolated-current-to-plane-height} with
\eqref{eq:isolated-plane-to-current-height} and
\eqref{eq:finite-Q-plane-change}, we obtain
\begin{equation}\label{eq:four-plane-isolated-core-height}
\begin{aligned}
&\int_{B_t}\dist\Bigl(x,\bigcup_lH_{j,l}\Bigr)^2\,d\|T_j\|(x)
 +\sum_{l=1}^Q\int_{H_{j,l}\cap B_t}
 \dist(p,\spt T_j)^2\,d\mathcal H^3(p)\\
\le&CB_j^2t^5(1+\log(b/t))\qquad\forall\,t\in(0,b/2).
\end{aligned}
\end{equation}

For the $l$th graph, keep $H_{j,i_j(l)}$ fixed.
By Lemma~\ref{lem:finite-annular-graph-matching}, the height on each
smaller fixed annulus is $O(B_j)$. If the $l$th graph is represented
over $H_{j,k}$ and $H_{j,k'}$ on overlapping annuli, then
$\|\pi_{H_{j,k}}-\pi_{H_{j,k'}}\|\le CB_j$.
After finitely many changes of projection, the height over
$H_{j,i_j(l)}$ is $O_t(B_j)$ on each fixed annulus.
Dividing by $B_j$ and passing diagonally to a subsequence, we obtain
a harmonic map $v_l$ on $(\Pi_l\cap B_b)\setminus\{0\}$, where
$\Pi_l\coloneqq\lim_{j\to\infty}H_{j,i_j(l)}$.

Express each $H_{j,k}$ converging to $\Pi_l$ as a linear graph over
$H_{j,i_j(l)}$, identifying $H_{j,i_j(l)}$ orthogonally with $\Pi_l$.
After division by $B_j$ and passage to a subsequence, the linear graph
functions either converge to $b_1,\ldots,b_s$, including zero, or have
unbounded operator norm. Normalize the unbounded linear graph functions
by the operator norm, as in the proof of
Lemma~\ref{lem:finite-annular-graph-matching}.
Outside the kernels of the nonzero linear limits and the intersections
of distinct limiting planes, only $b_1,\ldots,b_s$ can minimize the
distance to $v_l(p)$.
The excluded linear subspaces have measure zero.
By Fatou's lemma and \eqref{eq:four-plane-isolated-core-height},
\[
\int_{\Pi_l\cap B_t}\min_k|v_l(p)-b_k(p)|^2\,d\mathcal H^3(p)
\le Ct^5(1+\log(b/t)).
\]
Since each $b_k$ is linear,
\[
\begin{aligned}
\int_{\Pi_l\cap B_t}|v_l|^2\,d\mathcal H^3
\le&2\int_{\Pi_l\cap B_t}\min_k|v_l(p)-b_k(p)|^2\,d\mathcal H^3(p)
 +Ct^5\max_k\|b_k\|^2\\
\le&C't^5(1+|\log t|),
\end{aligned}
\]
where $C'$ may depend on the norms of the $b_k$. By the mean value
estimate,
\[
|v_l(p)|\le C'|p|\sqrt{1+|\log|p||}.
\]
Thus $v_l$ extends harmonically across zero, with $v_l(0)=0$.
The $L^2$ norms of $v_l$ on the fixed annulus where $F_j$ was
represented over $H_{j,i_j(l)}$ are uniformly bounded; hence, by the interior estimates for
the harmonic extensions of $v_l$,
\begin{equation}\label{eq:isolated-harmonic-derivative}
\sum_l|Dv_l(0)|^2\le C,
\end{equation}
with $C$ independent of the norms of the $b_k$.

Let $H'_{j,l}$ be the graph of $B_jDv_l(0)$ over $H_{j,i_j(l)}$,
using the orthogonal identifications with $\Pi_l$.
By \eqref{eq:isolated-harmonic-derivative},
$\|\pi_{H'_{j,l}}-\pi_{H_{j,i_j(l)}}\|\le CB_j$.
Moreover, $\min_l\|\pi_{H_{j,i}}-\pi_{H_{j,i_j(l)}}\|\le CB_j$
for each $i$, by Lemma~\ref{lem:finite-annular-graph-matching}.
Consequently,
\begin{equation}\label{eq:finite-Q-isolated-projections}
\operatorname{dist}_{\mathcal H}
 \bigl(\{\pi_{H'_{j,l}}\}_l,\{\pi_{H_{j,l}}\}_l\bigr)\le CB_j.
\end{equation}
To estimate $\mathcal E(T_j,H'_j;t)$, we apply Taylor's formula and
interior estimates to the harmonic extensions of $v_l$, obtaining
\begin{equation}\label{eq:isolated-harmonic-linear-estimates}
\sum_l\int_{\Pi_l\cap B_t}|v_l-Dv_l(0)p|^2\,d\mathcal H^3(p)
\le Ct^7,
\end{equation}
with $C$ independent of the norms of the $b_k$.
On $B_t\setminus B_\tau$, apply the graph area formula to
$v_l(p)-Dv_l(0)p$. By \eqref{eq:four-plane-isolated-core-height},
\eqref{eq:finite-Q-plane-change}, and \eqref{eq:stationary-monotonicity},
for $0<\tau<t/2$,
\[
\begin{aligned}
\Bigl(\frac\tau t\Bigr)^5
 \frac{\mathcal E(T_j,H'_j;\tau)^2}{B_j^2}
\le&C\Bigl(\frac\tau t\Bigr)^5
 \Bigl(\frac{\mathcal E(T_j,H_j;\tau)^2}{B_j^2}+1\Bigr)\\
\le&C\Bigl(\frac\tau t\Bigr)^5
 \bigl(1+\log(b/\tau)\bigr)\longrightarrow0
 \qquad\text{as }\tau\downarrow0.
\end{aligned}
\]
Using \eqref{eq:isolated-harmonic-linear-estimates} on
$B_t\setminus B_\tau$, and passing first to $j\to\infty$ and then
to $\tau\downarrow0$, we find
\begin{equation}\label{eq:finite-Q-isolated-step}
\limsup_j\mathcal E(T_j,H'_j;t)^2/B_j^2\le Ct^2.
\end{equation}

The assumption $\bigcap_\alpha P_\alpha=\{0\}$ was used only to show
that the points of density at least $Q$ tend to zero. Thus
\eqref{eq:finite-Q-isolated-projections} and
\eqref{eq:finite-Q-isolated-step} also hold whenever those points tend
to zero on a fixed ball, even if the limiting planes share a line.

Now suppose that $\bigcap_\alpha P_\alpha=L$ is a line.
For $p\ge2$, consider distinct $H_{j,1},\ldots,H_{j,p}$ with
\[
H_{j,i}\cap H_{j,k}=L_j\quad(i\ne k),\qquad L_j\to L\quad\text{as }j\to\infty.
\]

In the infima below, $C$ has total multiplicity $Q$, lies in the
prescribed $\mathcal U_\alpha$, and has fewer than $p$ distinct planes
containing $L_j$; two distinct planes in $C$ intersect exactly in $L_j$.
First suppose that $\inf_{\#\{C_l\}<p}\mathcal Q_{L_j}(T_j,C;1)>0$ and
\begin{equation}\label{eq:finite-Q-common-line-fine-gap}
B_j\coloneqq\mathcal Q_{L_j}(T_j,H_j;1)\longrightarrow0,
\qquad
\frac{B_j}{\displaystyle\inf_{\#\{C_l\}<p}\mathcal Q_{L_j}(T_j,C;1)}
\longrightarrow0\quad\text{as }j\to\infty.
\end{equation}
Because at least two of the limiting $P_\alpha$ are distinct,
for all large $j$,
\begin{equation}\label{eq:finite-Q-one-plane-gap}
\inf_{\#\{C_l\}=1}\mathcal Q_{L_j}(T_j,C;1)\ge c>0,
\end{equation}
where $C$ consists of one plane containing $L_j$, counted with
multiplicity $Q$.
Applying \cite[Theorems~4.4 and~4.20]{KrummelWickramasekera2026II}
with \eqref{eq:finite-Q-common-line-fine-gap} and
\eqref{eq:finite-Q-one-plane-gap}, for each fixed $\tau>0$ and all large
$j$ we have, on $\{\dist(x,L_j)>\tau\}$, currents $T_{j,i}$ with
positive projection degrees $q_i$, $\sum_iq_i=Q$.
Count the degrees on fixed outer disks and pass to a subsequence
so that $q_i$ is independent of $j$.
After assigning multiplicity $q_i$ to $H_{j,i}$,
$\mathcal Q_{L_j}(T_j,H_j;1)$ changes by at most a factor $\sqrt Q$.
By \cite[Theorem~4.4(a)]{KrummelWickramasekera2026II}, the smallest
nonzero principal sine between distinct planes $H_{j,i}$ divided by
$B_j$ tends to infinity. Since the graph heights of $T_{j,i}$ over $H_{j,i}$ are at most
$C_\tau B_j$, the supports $\spt T_{j,i}$ are separated on
$\{\dist(x,L_j)>\tau\}$.
Restriction of $F_j$ to each support is therefore open and
closed, proper, boundaryless and stable outside $S_j$.
Applying Proposition~\ref{prop:null-set-plane-regularity} and the minimal
graph estimates on fixed smaller balls, we extend the restrictions of
$F_j$ over $\spt T_{j,i}$ across $S_j$ and obtain normal projection covers
of degrees $q_i$ whose graph heights have $C^k$ norm at most $C_{\tau,k}B_j$.

On $\overline B_{1/2}$, pass to a Hausdorff limit $D$ of
\[
\{y\in\overline B_{1/2}:\Theta^3(\|T_j\|,y)\ge Q\}.
\]
Then $D\subset L$ and $0\in D$. On every fixed ball compactly contained
in $B_{1/2}\setminus D$, $\Theta^3(\|T_j\|,\cdot)<Q$ for all large $j$.
Apply the assumed compactness to each constant sequence
$T_j,T_j,\ldots$ to extend $F_j$ across $S_j$ on that ball.
On balls disjoint from $D$, the extensions of $F_j$ across $S_j$
converge smoothly by Proposition~\ref{thm:flat-union}\textup{(i)}. In particular, on balls about
$L\setminus D$, Lemma~\ref{lem:finite-annular-graph-matching}
represents each completed local graph over one $H_{j,i}$
with $C^2$ height $O(B_j)$.
If one local graph is assigned to both $H_{j,i}$ and $H_{j,k}$,
then, by the $C^1$ estimates over $H_{j,i}$ and $H_{j,k}$
on a smaller fixed disk,
\[
\|\pi_{H_{j,i}}-\pi_{H_{j,k}}\|\le C B_j,
\]
contradicting the separation in \cite[Theorem~4.4(a)]{KrummelWickramasekera2026II}.
Thus, on overlapping subsets of $M_j$, the graph functions of $F_j$
are defined over $H_{j,i}$ with the same $i$.

After passing to a subsequence, let $\Pi_i\coloneqq\lim_jH_{j,i}$.
For $K\Subset\Pi_i\setminus D$, choose
$K\subset\bigcup_{\nu=1}^{N_K}B_{r_\nu}(z_\nu)$ with
$\overline B_{r_\nu}(z_\nu)\cap D=\varnothing$.
The extensions of the local graph functions of $F_j$ agree on the
dense subsets where $F_j$ is defined outside $S_j$, so the identifications on
overlaps are permutations of the $q_i$ lifts of the projection onto
$H_{j,i}$. Choosing a diagonal subsequence on which the permutations
on overlaps are fixed, we obtain locally smooth convergence of the
heights divided by $B_j$ to harmonic maps $v_i$ on degree-$q_i$ covers
$p_i$ of $(\Pi_i\setminus D)\cap B_{1/4}$.

For points $z_j$ of density at least $Q$ with
$z_j\to z\in D\cap B_{1/8}$ as $j\to\infty$, apply
\cite[Theorem~4.20(a),(d) and the proof of Theorem~4.20]{KrummelWickramasekera2026II}
with exponent $1/Q$ and a fixed $0<r_c<1/16$. Then
\[
\int_{B_{r_c}(z_j)}
 \frac{\dist(x-z_j,\bigcup_iH_{j,i})^2}{|x-z_j|^{5-1/Q}}
 \,d\|T_j\|(x)\le CB_j^2.
\]
We claim that, after a further subsequence for this $z$,
\begin{equation}\label{eq:finite-Q-common-center-trace}
\begin{gathered}
|\pi_{L_j^\perp}z_j|\le CB_j,\qquad
B_j^{-1}\pi_{L_j^\perp}z_j\longrightarrow c_z\in L^\perp
\quad\text{as }j\to\infty,\\
\sum_i\int_{p_i^{-1}(\Pi_i\cap B_{r_c}(z))}
 \frac{|v_i(\widetilde x)-\pi_{\Pi_i^\perp}c_z|^2}
 {|p_i(\widetilde x)-z|^{5-1/Q}}\,d\mu_i(\widetilde x)\le C.
\end{gathered}
\end{equation}
Here $d\mu_i$ is the volume measure of the pullback Euclidean metric.
The harmonic maps $v_i$ do not depend on $z$; the further
subsequence of $z_j$ may depend on $z$.
Indeed, to pass from the Lipschitz graph approximations of $T_{j,i}$ to the
smooth projection covers of $F_j$, first restrict to fixed regions away
from $L$. The Lipschitz and smooth graph heights over $H_{j,i}$ are
$O(B_j)$ and agree outside sets of measure $O(B_j^{2+\alpha})$.
After division by $B_j^2$, the squared $L^2$ distance between the
Lipschitz and smooth graph functions is at most
$CB_j^{-2}B_j^2B_j^{2+\alpha}=CB_j^{2+\alpha}\to0$ as $j\to\infty$.
Then exhaust towards $L$, using the independence of $C$ in
\eqref{eq:finite-Q-common-center-trace} from $\tau$.

Cover $D\cap B_{1/8}$ by at most $C\tau^{-1}$ balls of radius $2\tau$
centered at nearby points of density at least $Q$. By
\eqref{eq:finite-Q-common-center-trace},
\begin{equation}\label{eq:finite-Q-density-tube}
\begin{aligned}
&\limsup_j B_j^{-2}
\int_{B_{1/16}\cap\{\dist(x,D)<\tau\}}
 \dist\Bigl(x,\bigcup_iH_{j,i}\Bigr)^2\,d\|T_j\|(x)\\
\le&C\tau^{-1}(\tau^{5-1/Q}+\tau^3)
\le C\tau^2\qquad\forall\,\tau\in(0,1/32).
\end{aligned}
\end{equation}
On lifted balls, applying the mean value estimate with
\eqref{eq:finite-Q-common-center-trace}, we also have
\[
|v_i(\widetilde x)-\pi_{\Pi_i^\perp}c_z|
\le C\dist(x,D)^{1-1/(2Q)},
\qquad x=p_i(\widetilde x),
\]
where $z$ is a nearest point in $D$.
For $z_j=0$, $c_0=0$; by
\eqref{eq:finite-Q-common-center-trace} with $z=0$,
\begin{equation}\label{eq:finite-Q-cover-L2}
\sum_i\int_{p_i^{-1}(\Pi_i\cap B_{r_c/2})}|v_i|^2\,d\mu_i\le C.
\end{equation}

Work on a smaller ball inside $B_{r_c}$.
If $D$ misses a point of $L$ in the larger ball, move each normal
circle along $L$ to a small ball disjoint from $D$ and contract the
circle in that ball. Loops in the complement of a closed subset
of a line are generated by these normal circles, so each loop
contracts when the line has such a gap.
Every fixed homotopy in $B_{r_c}\setminus D$ has positive distance
from $D$ and therefore lifts in the projection covers of $F_j$ for
all large $j$. Thus each cover $p_i$ splits on compact exhaustions.
By \eqref{eq:finite-Q-common-center-trace}, $|v_i|\le C$ near $D$.
Using logarithmic cutoffs, we extend $v_i$ harmonically across $D$,
with $v_i(0)=0$. The constant in \eqref{eq:finite-Q-cover-L2}
is independent of the location and width of the gap, so the interior
Taylor constants are uniform as well.

If $D$ fills the line in the larger ball, each connected cover
is cyclic of degree $d\le Q$. For a cover projection $p$,
average over the $d$ lifts to obtain a single valued harmonic map
$\overline v$, which extends across $L$ with $\overline v(0)=0$.
By \eqref{eq:finite-Q-common-center-trace},
$|v-\overline v\circ p|(\widetilde x)
\le C\dist(p(\widetilde x),L)^{1-1/(2Q)}$.
For cutoff functions $\chi_\varepsilon$ near $L$ on a fixed ball,
\[
\begin{aligned}
&\int_{p^{-1}(B_{r_c}\cap\{\dist(\cdot,L)<2\varepsilon\})}
 |v-\overline v\circ p|^2|\nabla\chi_\varepsilon|^2\,d\mu_i 
\le C\varepsilon^{2-1/Q}\longrightarrow0
\qquad\text{as }\varepsilon\downarrow0.
\end{aligned}
\]
By Caccioppoli's inequality, $v-\overline v\circ p$ has finite angular energy.

The angular frequencies of the degree-$d$ cover $p$ are $k/d$.
For a positive noninteger frequency $\nu$, factor
$(1-s^2)^{\nu/2}$ from the polar eigenfunction, where $s$ is
the axial coordinate on the unit sphere. The polar operator for frequency $\nu$ is
\[
-\frac1{(1-s^2)^\nu}\frac d{ds}
 \Bigl((1-s^2)^{\nu+1}\frac d{ds}\Bigr)+\nu(\nu+1).
\]
For the displayed operator on finite energy functions, an eigenbasis
consists of the orthogonal polynomials for the weight $(1-s^2)^\nu$, with
eigenvalues $(\nu+k)(\nu+k+1)$, $k=0,1,\ldots$. The radial factors
are $R^\ell$ and $R^{-\ell-1}$, where $\ell=\nu+k$.
At zero, \eqref{eq:finite-Q-common-center-trace} excludes
$R^{-\ell-1}$, since
\[
\int_0^{r_c}\frac{(R^{-\ell-1})^2}{R^{5-1/Q}}R^2\,dR
=\int_0^{r_c}R^{-2\ell-5+1/Q}\,dR=\infty.
\]
For a nonzero $R^\ell$ term, \eqref{eq:finite-Q-common-center-trace}
at zero requires
\[
\int_0^{r_c}R^{2\ell-3+1/Q}\,dR<\infty
\quad\Longleftrightarrow\quad
\ell>1-\frac1{2Q}>1-\frac1d.
\]
Since $\ell\in d^{-1}\mathbb N_0$, every nonzero term $R^\ell$ satisfies
$\ell=1$ or $\ell\ge1+1/d\ge1+1/Q$; the degree one part is single valued.

For each connected component $a$, of degree $d_a$ and projection
$p_a$, let $\ell_a$ be the degree one part of $v_a$.
Using
$\int_{p_a^{-1}(B_{r_c/2})}
\langle\ell_a\circ p_a,v_a-\ell_a\circ p_a\rangle\,d\mu_a=0$
and \eqref{eq:finite-Q-cover-L2}, we have
\[
\sum_a d_a\|\ell_a\|^2
\le C\sum_a\int_{p_a^{-1}(B_{r_c/2})}|v_a|^2\,d\mu_a
\le C.
\]
After subtracting $\ell_a\circ p_a$, only terms of degree at least
$1+1/Q$ remain. Since
$(2t/r_c)^{2\ell+3}\le(2t/r_c)^{5+2/Q}$ for every remaining degree
$\ell\ge1+1/Q$ and $0<t<r_c/2$, integrating the angular expansion yields
\begin{equation}\label{eq:finite-Q-cover-linear-remainder}
\begin{aligned}
\sum_a\int_{p_a^{-1}(B_t)}
 |v_a-\ell_a\circ p_a|^2\,d\mu_a
\le&\Bigl(\frac{2t}{r_c}\Bigr)^{5+2/Q}
 \sum_a\int_{p_a^{-1}(B_{r_c/2})}|v_a|^2\,d\mu_a
\le Ct^{5+2/Q}.
\end{aligned}
\end{equation}
If $p_a$ splits, Taylor's formula gives $Ct^7$ in place of
$Ct^{5+2/Q}$ in \eqref{eq:finite-Q-cover-linear-remainder}.
On $L$, $v_a=\overline v_a$; by Taylor's formula,
\[
|v_a(x)-D\overline v_a(0)x|\le C|x|^2
\quad\forall\,x\in L\cap B_{r_c/2}.
\]

To construct $H'_{j,a}$ containing a line $L'_j$, first choose
a small fixed scale $\theta$, then a small $a_0>0$.
If $z=te_L\in D$ for some $0<|t|<a_0$, with $e_L$ a unit
vector along $L$, then
\[
D\overline v_a(0)e_L
=\pi_{\Pi_i^\perp}(c_z/t)+O(a_0)
\]
on every component over $\Pi_i$.
For distinct planes $P_\alpha,P_\beta$, by
\eqref{eq:finite-Q-common-center-trace} and $v_a(z)=O(|t|)$,
\[
\bigl|c_z/t\bigr|
\le C\bigl(|\pi_{P_\alpha^\perp}(c_z/t)|
          +|\pi_{P_\beta^\perp}(c_z/t)|\bigr)\le C,
\]
where we used $P_\alpha\cap P_\beta=L$ and $c_z/t\in L^\perp$.
Change the linear functions $\ell_a$ by $O(a_0)$ and adjust the planes
$H'_{j,a}$ by $o(B_j)$. Assign multiplicity $d_a$ to component $a$,
with $\sum_a d_a=Q$. For $e_{L_j}\in L_j$ and
$w_j\to c_z/t$ as $j\to\infty$, choose
$L'_j\coloneqq\mathbb R(e_{L_j}+B_jw_j)$ and the adjusted planes so that
\[
L'_j\subset H'_{j,a},\qquad
H'_{j,a}\cap H'_{j,b}=L'_j\quad(H'_{j,a}\ne H'_{j,b}).
\]

For component $a$ over $H_{j,i}$, choose
$O_{j,a}:H_{j,i}\to H'_{j,a}$ by rotations through the principal angles.
Using $\sum_a d_a\|\ell_a\|^2\le C$, we have
\[
\|O_{j,a}-I\|\le CB_j,\qquad
\|\pi_{H_{j,i}}O_{j,a}-I\|
=\tfrac12\|O_{j,a}-I\|^2\le CB_j^2.
\]
Since $L'_j\to L$ and $L_j\to L$ as $j\to\infty$, for the fixed $\theta$ and all large $j$,
\begin{equation}\label{eq:total-three-moving-reverse-domain}
O_{j,a}^{-1}\bigl(H'_{j,a}\cap B_{\theta/2}
 \cap\{\dist(p,L'_j)>\theta/16\}\bigr)
\subset H_{j,i}\cap B_{3\theta/4}
 \cap\{\dist(p,L_j)>\theta/32\}.
\end{equation}

Let $u_{j,a,l}$, $l=1,\ldots,d_a$, be the local normal graph values
of component $a$. For $x$ in the set on the left hand side of
\eqref{eq:total-three-moving-reverse-domain},
\[
d_a\,\dist(O_{j,a}x,\spt T_j)^2
\le\sum_{l=1}^{d_a}
 |u_{j,a,l}(x)-(O_{j,a}-I)x|^2.
\]
The inequality for $d_a\,\dist(O_{j,a}x,\spt T_j)^2$ is invariant
under permutation of $u_{j,a,1},\ldots,u_{j,a,d_a}$.
By $\sum_a d_a\|\ell_a\|^2\le C$ and the changes
$O(a_0)$ in $\ell_a$,
\[
\operatorname{dist}_{\mathcal H}
 (\{\pi_{H'_{j,l}}\}_l,\{\pi_{H_{j,l}}\}_l)\le CB_j.
\]
Integrating the inequality for
$d_a\,\dist(O_{j,a}x,\spt T_j)^2$ and using
\eqref{eq:finite-Q-cover-linear-remainder} and
\eqref{eq:finite-Q-density-tube},
\[
\limsup_j\frac{\mathcal Q_{L'_j}(T_j,H'_j;\theta)^2}{B_j^2}
\le C(\theta^{2/Q}+a_0^2).
\]

If $D\cap B_{a_0}=\{0\}$ instead, extend $F_j$ across $S_j$ on a fixed
annulus inside $B_{a_0}$ and apply
Lemma~\ref{lem:finite-annular-graph-matching}.
The graph functions of $F_j$ over $H_{j,l}$ have $C^2$ norm
$O_{a_0}(B_j)$. Repeating
\eqref{eq:four-plane-isolated-radial-action}--\eqref{eq:four-plane-isolated-core-height}
for the graphs of $F_j$ over $H_{j,l}$, we obtain
\eqref{eq:four-plane-isolated-core-height}
on smaller balls.
The Taylor coefficients of $v_a$ are controlled by
\eqref{eq:finite-Q-common-center-trace} with $z=0$, so
\eqref{eq:finite-Q-isolated-step} holds at a second fixed scale chosen
after $a_0$.
For $0<\kappa<1$, choose $\theta,a_0$ with
$C(\theta^{2/Q}+a_0^2)<\kappa^2$ and then choose $t$ sufficiently
small in \eqref{eq:finite-Q-isolated-step}. Under
\eqref{eq:finite-Q-common-line-fine-gap}, for all large $j$,
\[
\mathcal Q_{L'_j}(T_j,H'_j;\theta)<\kappa B_j
\quad\text{or}\quad
\mathcal E(T_j,H'_j;t)<\kappa B_j.
\]
The two scales $\theta,t$ are independent of $j$.

To prove \eqref{eq:finite-Q-one-scale} when
$H_{j,i}\cap H_{j,k}=L_j$ for $i\ne k$, we induct on the number $p$
of distinct $H_{j,i}$, allowing up to $Q$ distinct planes in $H'_j$.
For $p=2$, by \eqref{eq:finite-Q-one-plane-gap},
\[
\frac{B_j}{\inf_{\#\{C_l\}=1}\mathcal Q_{L_j}(T_j,C;1)}
\le\frac{B_j}{c}\longrightarrow0\quad\text{as }j\to\infty.
\]
To pass to general $p$, first fix the scales for
$D\cap B_{a_0}\ne\{0\}$ and $D\cap B_{a_0}=\{0\}$, then choose a
fixed $0<\beta<1$ as follows. If no such $\beta$ exists, choose a
sequence for which
$B_j/\inf_{\#\{C_l\}<p}\mathcal Q_{L_j}(T_j,C;1)\to0$ as $j\to\infty$ and both
alternatives in \eqref{eq:finite-Q-one-scale} fail.
Pass to a Hausdorff limit $D$ of the points of density at least $Q$.
The estimates just proved for $D\cap B_{a_0}\ne\{0\}$ and
$D\cap B_{a_0}=\{0\}$ contradict this failure.
If $B_j>\beta\inf_{\#\{C_l\}<p}\mathcal Q_{L_j}(T_j,C;1)$,
choose $C_j$ with fewer than $p$ distinct planes containing $L_j$ and
\[
\mathcal Q_{L_j}(T_j,C_j;1)\le\frac{2B_j}{\beta}.
\]
On fixed disks away from $L$, apply
Lemma~\ref{lem:finite-annular-graph-matching} to the integrals with
respect to $\|T_j\|$ defining $\mathcal Q_{L_j}(T_j,H_j;1)$ and
$\mathcal Q_{L_j}(T_j,C_j;1)$.
Using the integrals over $H_{j,i}$ and $C_{j,l}$ in
\eqref{eq:truncated-plane-integrals} and the estimate for finite
families of linear maps in
Lemma~\ref{lem:finite-annular-graph-matching}, we obtain
\begin{equation}\label{eq:finite-Q-comparison-match}
\operatorname{dist}_{\mathcal H}
 \bigl(\{\pi_{H_{j,i}}\}_i,\{\pi_{C_{j,l}}\}_l\bigr)
\le C\bigl(B_j+\mathcal Q_{L_j}(T_j,C_j;1)\bigr)\le C_\beta B_j.
\end{equation}
The first inequality in \eqref{eq:finite-Q-comparison-match} also holds
when $B_j=\mathcal E(T_j,H_j;1)$,
since the reverse domains for $\mathcal E(T_j,H_j;1)$ include the
disks used in \eqref{eq:truncated-plane-integrals}.
Applying the induction hypothesis to $C_j$ with $\kappa\beta/4$ and using
\[
\frac{\kappa\beta}{4}\mathcal Q_{L_j}(T_j,C_j;1)
\le\frac{\kappa\beta}{4}\frac{2B_j}{\beta}
=\frac\kappa2B_j<\kappa B_j,
\]
we obtain $H'_j$ at one of the fixed scales $\theta$
such that, for all large $j$,
\[
\mathcal Q_{L'_j}(T_j,H'_j;\theta)\le\frac\kappa2B_j
\quad\text{or}\quad
\mathcal E(T_j,H'_j;\theta)\le\frac\kappa2B_j.
\]
The scales in this application are chosen after $\beta$.

If $\mathcal E(T_j,H_j;1)\to0$ or
$\mathcal Q_{L_j}(T_j,H_j;1)\to0$ as $j\to\infty$, each limit of $H_{j,l}$ belongs to $\{P_1,\ldots,P_N\}$ by the
integrals over $H_{j,l}$ in
\eqref{eq:finite-plane-integrals}--\eqref{eq:truncated-plane-integrals},
and every $P_\alpha$ occurs by the integrals with respect to $\|T_j\|$.
Thus the reduction from $H_j$ to $C_j$ preserves the limiting planes
$P_\alpha$ and total multiplicity $Q$, with multiplicities assigned
by the projection degrees. A single plane is excluded
by \eqref{eq:finite-Q-one-plane-gap}.
We have obtained finitely many possible scales for all $p\le Q$.

It remains to consider $B_j\coloneqq\mathcal E(T_j,H_j;1)$
when no line is contained in every $H_{j,l}$. Define
$\varepsilon_j\coloneqq\inf_C\mathcal E(T_j,C;1)$, where $C$ ranges over
the closure of the families of planes of total multiplicity $Q$ whose
planes contain a common line and lie in the prescribed $\mathcal U_\alpha$.
The closure is taken with respect to convergence of the orthogonal
projection matrices after relabeling, with each plane repeated
according to its multiplicity.
Suppose first that $\varepsilon_j>0$ and $B_j/\varepsilon_j\to0$ as $j\to\infty$.
If $\bigcap_l\lim_jH_{j,l}=\{0\}$, apply
\eqref{eq:finite-Q-isolated-step}. Otherwise $\varepsilon_j\to0$ as $j\to\infty$.
Choose $C_j$ and a line $L_j\subset\bigcap_l C_{j,l}$ such that
\[
\mathcal Q_{L_j}(T_j,C_j;1)
\le\mathcal E(T_j,C_j;1)\le2\varepsilon_j.
\]
In the following infima, $D$ has total multiplicity $Q$, lies in
the prescribed $\mathcal U_\alpha$, and has fewer distinct planes than
$C_j$; the distinct planes contain $L_j$ and intersect pairwise exactly
in $L_j$. If $\mathcal Q_{L_j}(T_j,C_j;1)=0$ for infinitely many $j$,
pass to that subsequence and use the case $B_j=0$ preceding
\eqref{eq:finite-Q-one-scale}. Otherwise assume
$\mathcal Q_{L_j}(T_j,C_j;1)>0$. If, along a subsequence,
\[
\inf_{\#\{D_l\}<\#\{C_{j,l}\}}\mathcal Q_{L_j}(T_j,D;1)
\le C\mathcal Q_{L_j}(T_j,C_j;1),
\]
choose $D_j$ with fewer distinct planes containing $L_j$ such that
\[
\mathcal Q_{L_j}(T_j,D_j;1)
\le C\mathcal Q_{L_j}(T_j,C_j;1)\le C\varepsilon_j.
\]
Replace $C_j$ by $D_j$ and repeat, treating
$\mathcal Q_{L_j}(T_j,C_j;1)=0$ as above. Otherwise, after at most
$Q-2$ replacements, the denominator below is positive and
\[
\begin{gathered}
0<s_j\coloneqq\mathcal Q_{L_j}(T_j,C_j;1)\le C\varepsilon_j,\\
\frac{s_j}
 {\displaystyle\inf_{\#\{D_l\}<\#\{C_{j,l}\}}
   \mathcal Q_{L_j}(T_j,D;1)}\longrightarrow0\quad\text{as }j\to\infty.
\end{gathered}
\]
For two planes, we use \eqref{eq:finite-Q-one-plane-gap}, since
$s_j\le C\varepsilon_j\to0$ as $j\to\infty$. The constant in $s_j\le C\varepsilon_j$ may depend on the subsequence.

Suppose points of density at least $Q$ satisfy $z_j\to z\ne0$ as $j\to\infty$.
If $s_j=o(\varepsilon_j)$, by \eqref{eq:finite-Q-common-center-trace} applied to
$C_j$ and \eqref{eq:finite-Q-comparison-match},
\[
\max_l\dist(z_j,H_{j,l})
\le Cs_j+C|z_j|(B_j+s_j)=o(\varepsilon_j).
\]
If $0<c\le s_j/\varepsilon_j\le C$, the separation between distinct planes in
$C_j$, divided by $\varepsilon_j$, tends to infinity.
Dividing the graph heights of $F_j$ over $C_j$ and the linear graph
maps representing $H_{j,l}$ by $\varepsilon_j$, we obtain a harmonic limit $v$ on each
connected cover and finite linear limits $b_l$ on fixed disks.
For a fixed disk $U$ in a connected cover, Fatou's lemma and
\eqref{eq:finite-plane-integrals} imply
\[
\int_U\min_l|v-b_l|^2\,d\mu
\le C\lim_{j\to\infty}\frac{B_j^2}{\varepsilon_j^2}=0,
\]
where $d\mu$ is the pullback Euclidean volume measure and we used
$B_j=o(\varepsilon_j)$. Hence $v\in\{b_l\}_l$ almost everywhere, and
by unique continuation $v=b_l$ on each connected cover for some $l$.
By the integrals over $H_{j,l}$ in
\eqref{eq:finite-plane-integrals}, the harmonic map $v$ on at least one
connected cover equals $b_l$, for each $l$. Applying
\eqref{eq:finite-Q-common-center-trace} with normalization $\varepsilon_j$, we obtain
\[
\frac{\dist(z_j,H_{j,l})}{\varepsilon_j}\longrightarrow0
\quad\text{as }j\to\infty,
\quad\text{for each }H_{j,l}.
\]

Since $z\ne0$, rotate $H_{j,l}$ to $H'_{j,l}$ so that
\[
z_j\in\bigcap_lH'_{j,l},\qquad
\max_l\|\pi_{H'_{j,l}}-\pi_{H_{j,l}}\|=o(\varepsilon_j).
\]
The rotations from $H_{j,l}$ to $H'_{j,l}$ map $H_{j,l}\cap B_1$ onto
$H'_{j,l}\cap B_1$ in \eqref{eq:finite-plane-integrals}, and
\[
\varepsilon_j\le\mathcal E(T_j,H'_j;1)=o(\varepsilon_j),
\]
contradicting the definition of $\varepsilon_j$. All points of density at least $Q$
therefore tend to zero, so we use \eqref{eq:finite-Q-isolated-step}.

Fix the scale in \eqref{eq:finite-Q-isolated-step} first.
By the $B_j=o(\varepsilon_j)$ case, choose $M>0$ so that
\eqref{eq:finite-Q-isolated-step} holds whenever $\varepsilon_j>MB_j$.
When $0\le\varepsilon_j\le MB_j$, choose $C_j$ and a line
$L_j\subset\bigcap_l C_{j,l}$ such that
\[
\mathcal Q_{L_j}(T_j,C_j;1)
\le\mathcal E(T_j,C_j;1)
\le2\varepsilon_j+B_j\le(2M+1)B_j.
\]
Applying the common-line case of \eqref{eq:finite-Q-one-scale} to
$C_j$, with line $L_j$ and
$\kappa=[4(2M+1)]^{-1}$, and using
\[
\frac{\mathcal Q_{L_j}(T_j,C_j;1)}{4(2M+1)}
\le\frac{\mathcal E(T_j,C_j;1)}{4(2M+1)}
\le\frac{B_j}{4},
\]
we obtain $H'_j$ at one of the fixed scales $\theta$ with
\[
\mathcal Q_{L'_j}(T_j,H'_j;\theta)\le\frac{B_j}{4}
\quad\text{or}\quad
\mathcal E(T_j,H'_j;\theta)\le\frac{B_j}{4}.
\]
By \eqref{eq:finite-Q-comparison-match} and the $L^2$ estimates for
the linear graph maps on fixed disks,
\[
\operatorname{dist}_{\mathcal H}
 \bigl(\{\pi_{H'_{j,l}}\}_l,\{\pi_{H_{j,l}}\}_l\bigr)
\le C\bigl(B_j+\mathcal E(T_j,C_j;1)\bigr)\le C_MB_j.
\]
Before iterating \eqref{eq:finite-Q-one-scale}, we verify that
$|(D_{1/\theta})_\#T_j|$ stays in the fixed varifold neighborhood of
$\sum_\alpha q_\alpha|P_\alpha|$ on $B_2$ and
$\|(D_{1/\theta})_\#T_j\|(B_1)<(Q+1/2)\omega_3$.
By current constancy and the prescribed multiplicities in the
$\mathcal U_\alpha$, for $0<r<1/2$,
\[
\mathcal E(T,H;1)=0\ \text{or}\ \mathcal Q_L(T,H;1)=0
\quad\Longrightarrow\quad
\|T\|(B_r)=Q\omega_3r^3.
\]
For $B_j\to0$ as $j\to\infty$, pass to a subsequence on which
$\theta$ is independent of $j$ and each $H_{j,l}$ converges.
By the compactness argument in Proposition~\ref{thm:flat-union}\textup{(ii)}
and \eqref{eq:calibrated-mass-chain}, $T_j\rightharpoonup T$ and
$|T_j|\to|T|$ on $B_1$ as $j\to\infty$. The current $T$ is a planar
cone of total multiplicity $Q$ by constancy and the prescribed multiplicities.
Since $2\theta<1/2$ and $T$ is conical, on $B_2$ we have
\[
|(D_{1/\theta})_\#T_j|\longrightarrow|(D_{1/\theta})_\#T|=|T|
\quad\text{as }j\to\infty.
\]
The planes comprising $T$ lie in $\overline{\mathcal U_\alpha}$ with
multiplicities $q_\alpha$, so $|T|$ belongs to the fixed varifold
neighborhood of $\sum_\alpha q_\alpha|P_\alpha|$ on $B_2$.
Moreover, $\|T\|(\partial B_\theta)=0$, hence
\[
\|(D_{1/\theta})_\#T_j\|(B_1)
=\theta^{-3}\|T_j\|(B_\theta)
\longrightarrow\theta^{-3}\|T\|(B_\theta)
=Q\omega_3<(Q+1/2)\omega_3\quad\text{as }j\to\infty.
\]
Thus, for all large $j$, $|(D_{1/\theta})_\#T_j|$ belongs to the fixed
varifold neighborhood of $\sum_\alpha q_\alpha|P_\alpha|$ on $B_2$ and
$\|(D_{1/\theta})_\#T_j\|(B_1)<(Q+1/2)\omega_3$. Since the scale list is finite, choose $\varepsilon_*>0$
so that \eqref{eq:finite-Q-one-scale} and both rescaling conditions hold
whenever $\mathcal E(T_j,H_j;1)\le\varepsilon_*$ or
$\mathcal Q_{L_j}(T_j,H_j;1)\le\varepsilon_*$. Decrease
$\varepsilon_*$ further so that, in the operator norm,
\[
2C\varepsilon_*<\min_\alpha\operatorname{dist}
 \bigl(\pi_{P_\alpha},\{\pi_H:H\notin\mathcal U_\alpha\}\bigr).
\]

We now return to the original currents $T_j$ and fix $r_0$ as in
\eqref{eq:finite-Q-initial-comparison}. For each $j$ and a center $z_j$
with $\Theta^3(\|T_j\|,z_j)\ge Q$, translate $z_j$ to zero and start
with $H^{(0)}=(P_\alpha)_\alpha$, counted with multiplicities $q_\alpha$.
By \eqref{eq:finite-Q-initial-comparison},
$\mathcal E(T_j,H^{(0)};r_0)\to0$ uniformly in the center as $j\to\infty$.
Choose $j_0$ such that $\mathcal E(T_j,H^{(0)};r_0)<\varepsilon_*$
for $j\ge j_0$, and fix $j\ge j_0$.
Applying \eqref{eq:finite-Q-one-scale} successively to
$\bigl(x\mapsto x/r_k\bigr)_\#T_j$, starting with $\mathcal E$ at $k=0$,
we construct $r_k,H^{(k)}=(H_l^{(k)})_l$ with
$r_{k+1}=\theta_k r_k$, $0<\theta_{\min}\le\theta_k<1/4$, and
\begin{equation}\label{eq:finite-Q-iterated-comparison}
\begin{aligned}
&\mathcal E(T_j,H^{(k)};r_k)
 \le2^{-k}\mathcal E(T_j,H^{(0)};r_0)\\
&\quad\text{or}\quad
 \mathcal Q_{L^{(k)}}(T_j,H^{(k)};r_k)
 \le2^{-k}\mathcal E(T_j,H^{(0)};r_0).
\end{aligned}
\end{equation}
Here $L^{(k)}\subset\bigcap_lH_l^{(k)}$ is a line in the $\mathcal Q$ case,
and $\theta_{\min}>0$ is the minimum of the finite set of scales in
\eqref{eq:finite-Q-one-scale}. By \eqref{eq:finite-Q-one-scale-projections}
and \eqref{eq:finite-Q-iterated-comparison},
\[
\begin{aligned}
&\sum_{h=k}^{\infty}\operatorname{dist}_{\mathcal H}
 \bigl(\{\pi_{H_l^{(h+1)}}\}_l,\{\pi_{H_l^{(h)}}\}_l\bigr)\\
\le&C\mathcal E(T_j,H^{(0)};r_0)\sum_{h=k}^{\infty}2^{-h}
 =2^{1-k}C\mathcal E(T_j,H^{(0)};r_0).
\end{aligned}
\]
Hence for a family of at most $Q$ planes $H_{\infty,l}$,
\[
\begin{aligned}
&\operatorname{dist}_{\mathcal H}
 \bigl(\{\pi_{H_l^{(k)}}\}_l,\{\pi_{H_{\infty,l}}\}_l\bigr)\\
\le&2^{1-k}C\mathcal E(T_j,H^{(0)};r_0)
 \longrightarrow0\quad\text{as }k\to\infty.
\end{aligned}
\]
For $r_{k+1}<r\le r_k$, by the definitions of $\mathcal E$ and $\mathcal Q$,
\begin{equation}\label{eq:finite-Q-iterated-height}
\begin{aligned}
&r^{-5}\int_{B_r}
 \dist\Bigl(x,\bigcup_lH_{\infty,l}\Bigr)^2\,d\|T_j\|(x)\\
\le&2r^{-5}\int_{B_r}
 \dist\Bigl(x,\bigcup_lH_l^{(k)}\Bigr)^2\,d\|T_j\|(x)\\
&+C4^{-k}\mathcal E(T_j,H^{(0)};r_0)^2
 r^{-5}\int_{B_r}|x|^2\,d\|T_j\|(x)\\
\le&C\Bigl(\frac{r_k}{r}\Bigr)^5
 4^{-k}\mathcal E(T_j,H^{(0)};r_0)^2\\
\le&C\theta_{\min}^{-5}4^{-k}\mathcal E(T_j,H^{(0)};r_0)^2\\
\le&4C\theta_{\min}^{-5}\mathcal E(T_j,H^{(0)};r_0)^2
 \Bigl(\frac r{r_0}\Bigr)^{\log4/\log(1/\theta_{\min})},
\end{aligned}
\end{equation}
where we used \eqref{eq:finite-Q-iterated-comparison} and
$\|T_j\|(B_r)\le Cr^3$ in the second inequality,
$r>\theta_{\min}r_k$ in the third, and
$r>r_{k+1}\ge\theta_{\min}^{k+1}r_0$ in the last.
Thus $\spt\|V\|\subset\bigcup_lH_{\infty,l}$ for every tangent varifold
$V$ of $|T_j|$ at zero.

The choice of $\varepsilon_*$ and
$\mathcal E(T_j,H^{(0)};r_0)<\varepsilon_*$ ensure that every
$H_{\infty,l}$ belongs to $\bigcup_\alpha\mathcal U_\alpha$.
For $x\in M_j$ with $\Theta^3(\|T_j\|,F_j(x))\ge Q$, apply the
iteration centered at $F_j(x)$. The local smooth sheet through $x$
converges under dilation about $F_j(x)$ to $dF_j(T_xM_j)$. Thus every
tangent varifold $V$ of $|T_j|$ at $F_j(x)$ satisfies
\[
dF_j(T_xM_j)\subset\spt V\subset\bigcup_lH_{\infty,l}
\quad\Longrightarrow\quad
dF_j(T_xM_j)\in\bigcup_\alpha\mathcal U_\alpha.
\]
Next take $x_j\in M_j$ with $F_j(x_j)\to0$ as $j\to\infty$ and
$\Theta^3(\|T_j\|,F_j(x_j))<Q$. If
\begin{equation}\label{eq:finite-Q-density-distance}
0<\dist\bigl(F_j(x_j),\{y\in B_2:\Theta^3(\|T_j\|,y)\ge Q\}\bigr)
\longrightarrow0\quad\text{as }j\to\infty,
\end{equation}
dilate about $F_j(x_j)$ by the reciprocal distance.
By \eqref{eq:finite-Q-iterated-height}, the dilations of $T_j$ about
$F_j(x_j)$ have a planar limit. In the rescaled coordinates, the dilations of $T_j$ have density less
than $Q$ on $B_{1/2}$. Applying the assumed removable extension below density $Q$
and then Proposition~\ref{thm:flat-union}\textup{(i)} on this ball, we obtain
$dF_j(T_{x_j}M_j)\in\bigcup_\alpha\mathcal U_\alpha$ for all large $j$.
If the distance in \eqref{eq:finite-Q-density-distance} is
$\ge c>0$ for all large $j$, or
$\{y\in B_2:\Theta^3(\|T_j\|,y)\ge Q\}=\varnothing$, apply compactness below density $Q$ and
Proposition~\ref{thm:flat-union}\textup{(i)} on a fixed ball about
$F_j(x_j)$ without dilation.

By \eqref{eq:finite-Q-tangent-neighborhood} and
Proposition~\ref{prop:null-set-plane-regularity}, the restrictions of
$F_j$ with tangent planes in $\mathcal U_\alpha$ extend across $S_j$
and converge smoothly to $P_\alpha$ near zero. Repeating the argument
after translation and rescaling at points of $\bigcup_\alpha P_\alpha$
and taking a finite cover, we obtain smooth convergence of the extended
immersions on compact subsets of $B_2$.
\end{proof}

\section{Rigidity with finite total curvature}\label{sec:complete}

\begin{theorem}\label{thm:finite-total-rigidity}
Let $F:M^n\to\C^n$, $n\ge3$, be a complete connected proper special
Lagrangian immersion without boundary. If $F$ is $\delta$-Simons stable
for $\delta>\delta_{\rm K}(n)$ and $\int_M|A_M|^n\,d\mu<\infty$, then $F$
is flat.
\end{theorem}

\begin{proof}
Fix $z\in\C^n$, $0<s<t$,
and $\lambda\ge0$, and set $\Omega_\tau\coloneqq F^{-1}(B_\tau(z))$. Assume
$\Delta|A_M|^2\ge-\lambda|A_M|^2$ on $\Omega_t$. For $\beta\ge n/2$ and
$\eta\in\operatorname{Lip}_c(\Omega_t)$, $0\le\eta\le1$, multiplying by
$\eta^2|A_M|^{2\beta-2}$ and integrating by parts over $\Omega_t$, we obtain
\[
\begin{aligned}
\frac{4(\beta-1)}{\beta^2}
\int_{\Omega_t}\eta^2|\nabla(|A_M|^\beta)|^2\,d\mu
\le&\lambda\int_{\Omega_t}\eta^2|A_M|^{2\beta}\,d\mu\\
&-\frac4\beta\int_{\Omega_t}\eta|A_M|^\beta
 \langle\nabla(|A_M|^\beta),\nabla\eta\rangle\,d\mu.
\end{aligned}
\]
By \eqref{eq:MS} and Young's inequality,
\[
\begin{aligned}
 \Bigl(\int_{\Omega_t}(\eta|A_M|^\beta)^{\frac{2n}{n-2}}\,d\mu\Bigr)^{\frac{n-2}{n}} 
\le&C_n\int_{\Omega_t}|\nabla(\eta|A_M|^\beta)|^2\,d\mu\\
\le&2C_n\int_{\Omega_t}\eta^2|\nabla(|A_M|^\beta)|^2
 +|A_M|^{2\beta}|\nabla\eta|^2\,d\mu\\
\le&C_n\beta^2\int_{\Omega_t}|A_M|^{2\beta}
 (|\nabla\eta|^2+\lambda\eta^2)\,d\mu.
\end{aligned}
\]
In particular, with $\eta\in\operatorname{Lip}_c(\Omega_t)$,
$\eta|_{\Omega_s}=1$ and $|\nabla\eta|\le2/(t-s)$,
\begin{equation}\label{eq:finite-total-iteration}
\Bigl(\int_{\Omega_s}|A_M|^{\frac{2n\beta}{n-2}}\,d\mu\Bigr)^{\frac{n-2}{n\beta}}
\le\bigl[C_n\beta^2(\lambda+(t-s)^{-2})\bigr]^{1/\beta}
 \Bigl(\int_{\Omega_t}|A_M|^{2\beta}\,d\mu\Bigr)^{1/\beta}.
\end{equation}
For $0<\rho\le t$, rescaling $F\mapsto(F-z)/\rho$, applying
\eqref{eq:finite-total-iteration} successively with
$\beta=\frac n2(\frac n{n-2})^k$, $t=\frac12+2^{-k-1}$ and
$s=\frac12+2^{-k-2}$ for $k=0,1,\ldots$, and letting $k\to\infty$,
we obtain
\begin{equation}\label{eq:finite-total-linear-mean}
\sup_{\Omega_{\rho/2}}|A_M|^2
\le C\rho^{-2}
 \Bigl(\int_{\Omega_\rho}|A_M|^n\,d\mu\Bigr)^{2/n},
\end{equation}
where $C=C(n,\lambda\rho^2)>0$, since
$\sum_{k\ge0}(k+1)(\frac{n-2}{n})^k<\infty$.

Let
$C_n\ge\max\{C(n,12),2\sqrt{C(n,12)}\}$ and choose
$\varepsilon_n>0$ so that
$C_n\varepsilon_n^{2/n}<1$. Assume that
\[
\int_{F^{-1}(B_R(F(x)))}|A_M|^n\,d\mu<\varepsilon_n.
\]
Since $F$ is proper, choose $y\in F^{-1}(\overline B_R(F(x)))$ with
\[
(R-|F(y)-F(x)|)|A_M|(y)
=\max_{z\in F^{-1}(\overline B_R(F(x)))}
 (R-|F(z)-F(x)|)|A_M|(z).
\]
If $(R-|F(y)-F(x)|)|A_M|(y)>2$, then
$2|A_M|(y)^{-1}<R-|F(y)-F(x)|$, and hence
$B_{|A_M|(y)^{-1}}(F(y))\subset B_R(F(x))$.
By maximality, $|A_M|\le2|A_M|(y)$ on
$F^{-1}(B_{|A_M|(y)^{-1}}(F(y)))$; thus, by
\eqref{eq:exact-simons} and \eqref{eq:P-estimates},
\[
\Delta|A_M|^2\ge-3|A_M|^4\ge-12|A_M|^2(y)|A_M|^2.
\]
Applying \eqref{eq:finite-total-linear-mean} with center $F(y)$,
$\lambda=12|A_M|^2(y)$, and $\rho=|A_M|(y)^{-1}$, we obtain
\[
1\le C_n\Bigl(
\int_{F^{-1}(B_{|A_M|(y)^{-1}}(F(y)))}|A_M|^n\,d\mu
\Bigr)^{2/n}<C_n\varepsilon_n^{2/n}<1,
\]
a contradiction.
Hence
\[
\max_{F^{-1}(\overline B_R(F(x)))}
 (R-|F-F(x)|)|A_M|\le2.
\]
Using $R-|F-F(x)|\ge R/2$ and $|A_M|\le4/R$ on
$F^{-1}(B_{R/2}(F(x)))$, we apply
\eqref{eq:finite-total-linear-mean} with
$\lambda=48/R^2$ and $\rho=R/2$ to find
\[
R|A_M|(x)\le C_n
 \Bigl(\int_{F^{-1}(B_R(F(x)))}|A_M|^n\,d\mu\Bigr)^{1/n}.
\]
For $r(x)\coloneqq |F(x)|$, taking $R=r(x)/2$ and using finite total
curvature and $B_{r(x)/2}(F(x))\subset\{|z|>r(x)/2\}$, we obtain
\begin{equation}\label{eq:finite-total-initial-decay}
r(x)|A_M|(x)\le2C_n
 \Bigl(\int_{\{r>r(x)/2\}}|A_M|^n\,d\mu\Bigr)^{1/n}
\longrightarrow0\quad\text{as }r(x)\to\infty.
\end{equation}
Since $F$ is minimal, $r\Delta r=n-|\nabla r|^2$ and
$|\nabla r|\le1$, so for $0<c<n-2$,
\begin{equation}\label{eq:finite-total-radial-barrier}
\Delta r^{-c}
=cr^{-c-2}\bigl((c+2)|\nabla r|^2-n\bigr)
\le-c(n-2-c)r^{-c-2}\quad(0<c<n-2).
\end{equation}
Since $r^2\Pcal\to0$ as $r\to\infty$ by
\eqref{eq:finite-total-initial-decay} and \eqref{eq:P-estimates},
we may choose $R_0>0$ so that, by
\eqref{eq:finite-total-radial-barrier}, for
$c\in\{3(n-2)/4,\delta_{\rm K}(n)/2\}$ and $r\ge R_0$,
\begin{equation}\label{eq:finite-total-strict-barriers}
(\Delta+\delta_{\rm K}(n)\Pcal)r^{-c}
\le-\tfrac12c(n-2-c)r^{-c-2}<0.
\end{equation}
With $C_0>0$ chosen so that
$|A_M|^{\delta_{\rm K}(n)}\le C_0r^{-3(n-2)/4}$ on $\{r=R_0\}$,
\eqref{eq:finite-total-initial-decay} implies
\[
v\coloneqq \frac{|A_M|^{\delta_{\rm K}(n)}-C_0r^{-3(n-2)/4}}
 {r^{-\delta_{\rm K}(n)/2}}
\longrightarrow0\quad\text{as }r\to\infty.
\]
On $\{|A_M|>0\}$, we also have, by \eqref{eq:scalar-simons},
\begin{equation}\label{eq:critical-ground}
(\Delta+\delta_{\rm K}(n)\Pcal)(|A_M|^{\delta_{\rm K}(n)})
=\delta_{\rm K}(n)|A_M|^{\delta_{\rm K}(n)-2}
 \big(|\nabla A_M|^2-\frac{n+4}{n+1}\bigl|\nabla|A_M|\bigr|^2\big)\ge0.
\end{equation}
If $\sup_{\{r\ge R_0\}}v>0$, then $v\le0$ on $\{r=R_0\}$ and
$v\to0$ at infinity. Properness therefore implies that the positive
maximum is attained at some $x_0\in\{r>R_0\}$, with
\[
v(x_0)=\max_{\{r\ge R_0\}}v>0,\quad
\nabla v(x_0)=0,\quad \Delta v(x_0)\le0,\quad |A_M|(x_0)>0.
\]
At $x_0$, by \eqref{eq:finite-total-strict-barriers} and
\eqref{eq:critical-ground},
\begin{equation}\label{eq:radial-comparison-maximum}
\begin{aligned}
0&<(\Delta+\delta_{\rm K}(n)\Pcal)(|A_M|^{\delta_{\rm K}(n)}-C_0r^{-3(n-2)/4})\\
 &=r^{-\delta_{\rm K}(n)/2}\Delta v
 +v(\Delta+\delta_{\rm K}(n)\Pcal)r^{-\delta_{\rm K}(n)/2}<0.
\end{aligned}
\end{equation}
Therefore,
\begin{equation}\label{eq:finite-total-improved-decay}
|A_M|^{\delta_{\rm K}(n)}\le C_0r^{-3(n-2)/4},\quad
|A_M|\le Cr^{-3(n+1)/4}.
\end{equation}

Fix $p\in M$ and let $\rho\coloneqq \dist_g(p,\cdot)$.
By \eqref{eq:finite-total-improved-decay}, properness, and
$\nabla^2(r^2/2)=g+\langle F,A_M\rangle$, we can choose $D>0$ such that
$\nabla^2(r^2/2)\ge g/2$ on $\{\rho\ge D\}$.
For a minimizing unit speed geodesic $\gamma:[0,\rho(x)]\to M$
from $p$ to $x$ with $\rho(x)>D$, integrating
$(r^2\circ\gamma/2)''\ge1/2$ twice from $D$, with
$|(r^2\circ\gamma)'(D)|\le2(D+|F(p)|)$, we obtain
$r(\gamma(t))\ge t/2-C$ for $D\le t\le\rho(x)$.
With \eqref{eq:finite-total-improved-decay} and
$t/2-C\le r(\gamma(t))\le t+|F(p)|$, since $n\ge3$,
\[
\begin{aligned}
\left|\Bigl(\frac{r^2\circ\gamma}{2}\Bigr)''-1\right|
&=\bigl|\langle F\circ\gamma,A_M(\dot\gamma,\dot\gamma)\rangle\bigr|\\
&\le C\,r(\gamma(t))^{1-3(n+1)/4}
\le Ct^{-2}.
\end{aligned}
\]
Thus $\bigl(r^2\circ\gamma/2\bigr)''=1+O(t^{-2})$.
Integrating $\bigl(r^2\circ\gamma/2\bigr)''$ over $[D,\rho(x)]$
and using $\rho(x)\ge r(x)-|F(p)|$, we obtain
\begin{equation}\label{eq:finite-total-radial-gradient}
\begin{aligned}
|\nabla r|(x)
&\ge r(x)^{-1}\Bigl(\frac{r^2\circ\gamma}{2}\Bigr)'(\rho(x)) \ge\frac{\rho(x)-C}{r(x)}\ge1-\frac C{r(x)}
\end{aligned}
\end{equation}
for all sufficiently large $r(x)$. After increasing $C$ and the radius,
assume $|\nabla r|^2\ge1-C/r>0$. 

Using \eqref{eq:finite-total-improved-decay} and
\eqref{eq:finite-total-radial-gradient}, we compute, for sufficiently large $r$,
\begin{equation}\label{eq:finite-total-sharp-barrier}
\begin{aligned}
&(\Delta+\delta_{\rm K}(n)\Pcal)(r^{2-n}-r^{1-n})\\
=&n(n-2)(|\nabla r|^2-1)r^{-n}-(n-1)\bigl((n+1)|\nabla r|^2-n\bigr)r^{-n-1}
 +\delta_{\rm K}(n)\Pcal(r^{2-n}-r^{1-n})\\
\le&-\tfrac{n-1}{2}r^{-n-1}
 +\delta_{\rm K}(n)\Pcal(r^{2-n}-r^{1-n})\\
\le&-\tfrac{n-1}{2}r^{-n-1}
 +\tfrac32\delta_{\rm K}(n)|A_M|^2r^{2-n}\\
\le&-\tfrac{n-1}{2}r^{-n-1}+Cr^{2-n-3(n+1)/2}<0,
\end{aligned}
\end{equation}
where we used $1-C/r\le|\nabla r|^2\le1$,
\eqref{eq:P-estimates}, and \eqref{eq:finite-total-improved-decay}.
Choose $R_1>\max\{R_0,2\}$ so that
\eqref{eq:finite-total-sharp-barrier} holds for $r\ge R_1$, and
$C_1>0$ with $|A_M|^{\delta_{\rm K}(n)}\le C_1(r^{2-n}-r^{1-n})$ on
$\{r=R_1\}$. By \eqref{eq:finite-total-improved-decay},
\[
\frac{|A_M|^{\delta_{\rm K}(n)}-C_1(r^{2-n}-r^{1-n})}
 {r^{-\delta_{\rm K}(n)/2}}\longrightarrow0\quad\text{as }r\to\infty.
\]
Applying the maximum principle as in \eqref{eq:radial-comparison-maximum},
with \eqref{eq:finite-total-sharp-barrier}, we obtain
\begin{equation}\label{eq:finite-total-sharp-decay}
|A_M|^{\delta_{\rm K}(n)}\le C_1(r^{2-n}-r^{1-n})\le C_1r^{2-n}.
\end{equation}
Since $\delta_{\rm K}(n)=(n-2)/(n+1)$, we conclude that
$|A_M|\le Cr^{-n-1}$.

For almost every sufficiently large
regular value $R$, integrating $\operatorname{div}_M F^T=n$ over
$\{r<R\}$ and using \eqref{eq:finite-total-radial-gradient} with the
coarea formula, we obtain
\[
\begin{aligned}
n\mu\{r<R\}= R\int_{\{r=R\}}|\nabla r|\,d\sigma 
\ge (R-C)\int_{\{r=R\}}\frac{1}{|\nabla r|}\,d\sigma
=(R-C)\frac{d}{dR}\mu\{r<R\}.
\end{aligned}
\]
Consequently, after integrating $\frac{d}{dR}\log\mu\{r<R\}\le n/(R-C)$ in $R$,
\begin{equation}\label{eq:finite-total-derived-volume}
\mu\{r<R\}\le\Lambda R^n\quad\forall\,R\ge1
\end{equation}
for some $\Lambda>0$.

Taking $\alpha=\delta_{\rm K}(n)$ in \eqref{eq:regularized-energy}
and letting $\eps\downarrow0$, we obtain, for
$\phi\in\operatorname{Lip}_c(M)$,
\begin{equation}\label{eq:supercritical-gap}
(\delta-\delta_{\rm K}(n))
\int_M\Pcal|A_M|^{2\delta_{\rm K}(n)}\phi^2\,d\mu
\le\int_M|A_M|^{2\delta_{\rm K}(n)}|\nabla\phi|^2\,d\mu.
\end{equation} 
Choose $\phi_R\in C_c^1(M)$ with $0\le\phi_R\le1$,
$\phi_R|_{\{r\le R\}}=1$, $\phi_R|_{\{r\ge2R\}}=0$, and
$|\nabla\phi_R|\le2/R$.
By \eqref{eq:finite-total-derived-volume} and
\eqref{eq:finite-total-sharp-decay},
\[
\int_M|A_M|^{2\delta_{\rm K}(n)}|\nabla\phi_R|^2\,d\mu
\le\frac C{R^2}R^{4-2n}\mu\{r<2R\}
\le C\Lambda R^{2-n}\longrightarrow0\quad\text{as }R\to\infty.
\]
Taking $\phi=\phi_R$ in \eqref{eq:supercritical-gap}, letting
$R\to\infty$, and using Fatou's lemma and $\Pcal\ge|A_M|^2/n$, we obtain
\[
\int_M|A_M|^{2\delta_{\rm K}(n)+2}\,d\mu
\le n\int_M\Pcal|A_M|^{2\delta_{\rm K}(n)}\,d\mu=0.
\]
Thus $A_M\equiv0$, so $F$ is flat.
\end{proof}

\section{Finite scalar index and global geometry}\label{sec:compactness}

To prove finite total curvature in Theorem~\ref{thm:intro-finite-index},
we first establish \eqref{eq:index-energy-bound}.

\begin{proposition}\label{thm:index-energy}
Let $F:M^n\to B_2\subset\C^n$, $n\ge3$, be a smooth proper special
Lagrangian immersion with one phase on all components and no boundary
over $B_2$. Assume $\mu(M)\le\Lambda<\infty$,
$\operatorname{ind}_{\delta_+}(M)\le I\in\mathbb N\cup\{0\}$, and
\begin{equation}\label{eq:index-energy-parameter}
\delta_+>\max\{\delta_{\rm K}(n),\delta_{\rm cone}(n)\}.
\end{equation}
Then, for $0<\theta<2$,
\begin{equation}\label{eq:index-energy-bound}
\int_{F^{-1}(B_\theta)}|A_M|^n\,d\mu\le C,
\end{equation}
where $C=C(n,\delta_+,\Lambda,I,\theta)>0$.
\end{proposition}

\begin{proof}
If $F$ is $\delta_+$-Simons stable on $F^{-1}(B_{2t}(z))$ and
$\mu(F^{-1}(B_{2t}(z)))\le Lt^n$, applying
Corollary~\ref{cor:curvature-and-growth-bernstein} to $(F-z)/t$, we have
\begin{equation}\label{eq:index-stable-ball}
t\sup_{F^{-1}(B_t(z))}|A_M|\le C,
\quad C=C(n,\delta_+,L).
\end{equation}

We claim that $A_M=0$ for a proper immersion
$F:M^n\to\C^n\setminus\{0\}$ that is $\delta_+$-Simons stable or is
a smooth limit of proper $\delta_+$-Simons stable immersions on
$\{a\le|z|\le b\}$ for all $0<a<b<\infty$. Assume a common phase
on all sheets and
\begin{equation}\label{eq:index-punctured-growth}
\mu\{0<r<t\}\le Lt^n\quad\forall\,t>0,\quad r\coloneqq |F|.
\end{equation}
When $F$ is a smooth limit of
$F_j:M_j\to U_j\subset\C^n\setminus\{0\}$, assume
\[
\{a\le|z|\le b\}\subset U_j,\quad
\mu_j(F_j^{-1}(B_t(0)))\le Lt^n
\]
for fixed $0<a<b<\infty$, $t>0$, and $j\ge j_0(a,b,t)$. Passing to the
limit as in \eqref{eq:regular-limit-stability}, we obtain
\eqref{eq:delta-stability} with $\delta=\delta_+$.
Extend $F_\#\llbracket M\rrbracket$ and the associated varifold $V_F$
by zero across the origin. Choose $0\le\chi_\epsilon\le1$ with
$\chi_\epsilon|_{B_\epsilon}=0$,
$\chi_\epsilon|_{\C^n\setminus B_{2\epsilon}}=1$, and
$|\nabla\chi_\epsilon|\le C/\epsilon$.
For a compactly supported smooth $(n-1)$-form $\omega$ on $\C^n$,
applying Stokes' theorem to $(\chi_\epsilon\circ F)F^*\omega$ and using
\eqref{eq:index-punctured-growth}, we obtain
\[
\begin{aligned}
\Bigl|\int_M(\chi_\epsilon\circ F)F^*(d\omega)\Bigr|
= \Bigl|\int_M d(\chi_\epsilon\circ F)\wedge F^*\omega\Bigr| 
\le&\frac{C\|\omega\|_\infty}{\epsilon}
\mu\{\epsilon<r<2\epsilon\}\\
\le&CL\|\omega\|_\infty\epsilon^{n-1}\longrightarrow0
\quad\text{as }\epsilon\downarrow0.
\end{aligned}
\]
Letting $\epsilon\downarrow0$, we obtain
$\partial(F_\#\llbracket M\rrbracket)=0$ in $\C^n$.
Since all sheets have the same calibration,
$V_F=|F_\#\llbracket M\rrbracket|$ is stationary in $\C^n$.

By \eqref{eq:index-punctured-growth},
$\|(D_{1/R})_\#V_F\|(B_s)\le Ls^n$ for $R,s>0$.
Applying the Federer--Fleming compactness theorem
\cite[Chapter~6, Theorem~3.11]{Simon2014GMT} and Allard's integral
compactness theorem \cite[Chapter~8, Theorem~5.8 and Remark~5.9]{Simon2014GMT},
and using \eqref{eq:calibrated-mass-chain}, we choose $R_j\to\infty$ as
$j\to\infty$ and a calibrated integral current $\mathcal C_\infty$ with
$\partial\mathcal C_\infty=0$ such that
\[
\begin{gathered}
(R_j^{-1}F)_\#\llbracket M\rrbracket
 \rightharpoonup\mathcal C_\infty\quad\text{locally as currents},\\
(D_{1/R_j})_\#V_F
 \longrightarrow|\mathcal C_\infty|\quad\text{locally as varifolds as }j\to\infty.
\end{gathered}
\]
By the stationary monotonicity identity
\cite[Chapter~4, \textup{(3.6)}]{Simon2014GMT},
$R^{-n}\|V_F\|(B_R)$ is nondecreasing and at most $L$.
By the varifold convergence at radii with
$\|\mathcal C_\infty\|(\partial B_s)=0$, followed by approximation
of arbitrary radii,
\[
\frac{\|\mathcal C_\infty\|(B_s)}{s^n}
=\lim_{j\to\infty}\frac{\|V_F\|(B_{R_js})}{(R_js)^n}
=\lim_{R\to\infty}\frac{\|V_F\|(B_R)}{R^n}
\qquad\forall\,s>0.
\]
Consequently, by \eqref{eq:stationary-monotonicity},
\[
\begin{aligned}
\int_{B_t\setminus B_s}
 \frac{|\pi_{(T_z\mathcal C_\infty)^\perp}z|^2}{|z|^{n+2}}
 \,d\|\mathcal C_\infty\|(z)
=&\frac{\|\mathcal C_\infty\|(B_t)}{t^n}
 -\frac{\|\mathcal C_\infty\|(B_s)}{s^n}=0
\quad\forall\,0<s<t.
\end{aligned}
\]
Thus $z\in T_z\mathcal C_\infty$ for
$\|\mathcal C_\infty\|$-almost every $z$, and
$(D_\lambda)_\#\mathcal C_\infty=\mathcal C_\infty$ for $\lambda>0$
by \cite[Chapter~6, Lemma~2.33]{Simon2014GMT}.
For $x\in\spt\mathcal C_\infty\setminus\{0\}$, every tangent cone of
$\mathcal C_\infty$ at $x$ is in $\mathcal C(n,\delta_+,0)$ by a diagonal selection of the
immersions $R_j^{-1}F$ rescaled at $x$, and is invariant under
translations along $\R x$. Applying
Proposition~\ref{prop:closure-cone-reduction} in transverse dimension
at most $n-1$ and then \eqref{eq:uniform-flat-neighborhood}, we
obtain $x\in\ImmReg\mathcal C_\infty$.
By Proposition~\ref{thm:flat-union}\textup{(ii)}, the immersions $R_j^{-1}F$
converge smoothly on $\{a\le|z|\le b\}$, $0<a<b<\infty$, to a proper
immersion representing $\mathcal C_\infty$ on $\C^n\setminus\{0\}$.
By \eqref{eq:regular-limit-stability} for the smoothly converging
immersions $R_j^{-1}F$, the immersion representing $\mathcal C_\infty$
is $\delta_+$-Simons stable. Applying Theorem~\ref{thm:intro-cone} on
each component of the compact immersion domain of
$\mathcal C_\infty\cap\Sph^{2n-1}$, we have $A_{\mathcal C_\infty}=0$.
If $\sup_{\{r\ge R\}}r|A_M|\not\to0$ as $R\to\infty$, choose
$r(x_j)\to\infty$ as $j\to\infty$, with $r(x_j)|A_M|(x_j)\ge\varepsilon>0$.
Applying the same blow-down argument with $R_j=r(x_j)$ and passing
to a subsequence, we obtain smooth convergence on
$\{1/2\le|z|\le2\}$. Since $A_{\mathcal C_\infty}=0$,
\[
|R_j^{-1}F(x_j)|=1,\quad
\varepsilon\le R_j|A_M|(x_j)\longrightarrow0\quad\text{as }j\to\infty.
\]
Applying the maximum principle as in
\eqref{eq:finite-total-strict-barriers}--\eqref{eq:radial-comparison-maximum}
and increasing $C$, we obtain $|A_M|\le Cr^{-3(n+1)/4}$ for $r\ge1$.
Also, $B_{r(x)/2}(F(x))\subset B_{3r(x)/2}(0)\setminus\{0\}$, so
$r|A_M|\le C$ by \eqref{eq:index-stable-ball} with $t=r(x)/4$.

For $\delta_{\rm K}(n)<\zeta<\min\{\delta_+,1,(n-2)/2\}$, taking
$\alpha=\zeta$, $\delta=\delta_+$, and
$\phi\in C_c^1(\{\epsilon<r<2R\})$ in
\eqref{eq:zero-level-energy}, with $\phi=1$ on
$\{2\epsilon\le r\le R\}$,
$|\nabla\phi|\le C/\epsilon$ on $\{\epsilon<r<2\epsilon\}$, and
$|\nabla\phi|\le C/R$ on $\{R<r<2R\}$, we derive
\begin{equation}\label{eq:index-flat-cutoff}
\begin{aligned}
(\delta_+-\zeta)\int_M\Pcal|A_M|^{2\zeta}\phi^2\,d\mu
\le&\int_M|A_M|^{2\zeta}|\nabla\phi|^2\,d\mu\\
\le&\frac C{\epsilon^2}\int_{\{\epsilon<r<2\epsilon\}}r^{-2\zeta}\,d\mu
+\frac C{R^2}\int_{\{R<r<2R\}}r^{-3\zeta(n+1)/2}\,d\mu\\
\le&C\epsilon^{n-2-2\zeta}+CR^{n-2-3\zeta(n+1)/2}.
\end{aligned}
\end{equation}
Since $n-2-2\zeta>0$ and $n-2-3\zeta(n+1)/2<0$, letting
$\epsilon\downarrow0$ and $R\to\infty$, and using
$\Pcal\ge|A_M|^2/n$, we conclude by Fatou's lemma
\[
\int_M\Pcal|A_M|^{2\zeta}\,d\mu=0,\quad A_M\equiv0.
\]

For the annular estimate, fix $z\in\C^n$, and let $r\coloneqq |F-z|$ and
$\mathcal A_{s,t}(z)\coloneqq \{s<|F-z|<t\}$. Suppose
$F|_{\mathcal A_{s,R}(z)}$ is proper, without boundary, and
$\delta_+$-Simons stable, with
$\mu(\mathcal A_{s,t}(z))\le Lt^n$ for $s<t\le R$. We claim that, for
$C=C(n,\delta_+,L)$,
\begin{equation}\label{eq:index-stable-annulus}
\int_{\mathcal A_{2s,R/2}(z)}|A_M|^n\,d\mu\le C\quad(R>4s).
\end{equation}
Indeed, by \eqref{eq:index-stable-ball}, $r|A_M|\le C$ on
$\mathcal A_{2s,R/2}(z)$. We will show that, for each $\varepsilon>0$,
a constant $K=K(n,\delta_+,L,\varepsilon)>2$ exists such that
\[
r|A_M|\le\varepsilon\quad(Ks\le r\le R/K).
\]
Otherwise, choose proper $\delta_+$-Simons stable immersions $F_j$
without boundary over the closed annuli $\{s_j\le|F_j-z_j|\le R_j\}$,
with $\mu(\mathcal A_{s_j,t}(z_j))\le Lt^n$ for $s_j<t\le R_j$,
and points $x_j$ such that, for $r_j\coloneqq |F_j-z_j|$,
\begin{equation}\label{eq:stable-annulus-rescaling-data}
\begin{gathered}
r_j(x_j)/s_j\to\infty,\quad R_j/r_j(x_j)\to\infty
\quad\text{as }j\to\infty,\quad
r_j(x_j)|A_{M_j}|(x_j)>\varepsilon.
\end{gathered}
\end{equation}
Then, for fixed $0<a<b<\infty$,
\[
\frac{s_j}{r_j(x_j)}<a<b<\frac{R_j}{r_j(x_j)}
\qquad\forall\,j\ge j_0(a,b).
\]
Thus $(F_j-z_j)/r_j(x_j)$ has no boundary over $\{a\le|z|\le b\}$.
By \eqref{eq:index-stable-ball} for $(F_j-z_j)/r_j(x_j)$ and
Schauder estimates for the graph functions of $(F_j-z_j)/r_j(x_j)$
satisfying \eqref{eq:minimal-graph-equations},
\[
\sup_{\{a r_j(x_j)\le r_j\le b r_j(x_j)\}}
 r_j(x_j)^{k+1}|\nabla^k A_{M_j}|\le C(a,b,k)
\]
for fixed $0<a<b<\infty$, $k\ge0$, and all large $j$.
Applying the compactness argument of Proposition~\ref{thm:flat-union}\textup{(ii)}
on an exhaustion of $\C^n\setminus\{0\}$, we obtain a proper smooth
limit $F_\infty:M_\infty\to\C^n\setminus\{0\}$ satisfying
\eqref{eq:index-punctured-growth}. The immersion $F_\infty$ is flat by
\eqref{eq:index-flat-cutoff}. Using \eqref{eq:stable-annulus-rescaling-data}
at $x_j$ and the smooth convergence of $(F_j-z_j)/r_j(x_j)$, we obtain
\[
\Bigl|\frac{F_j(x_j)-z_j}{r_j(x_j)}\Bigr|=1,\quad
\varepsilon<r_j(x_j)|A_{M_j}|(x_j)\longrightarrow0\quad\text{as }j\to\infty,
\]
a contradiction.

\emph{Claim.}
Let $0<s<R$ be regular values such that $\{s\le r\le R\}$ is compact
and disjoint from $\partial M$. 
If, for some $\Lambda>0$,
\begin{align}
\sup_{\{s\le r\le R\}}r|A_M|&\le\varepsilon_n,
\label{eq:annular-small-curvature}\\
\mu(\mathcal A_{s,t}(z))&\le\Lambda t^n\quad\forall\,t\in(s,R],
\label{eq:annular-volume-assumption}
\end{align}
then, for $T\ge1$ with $T^2s<R$,
\begin{equation}\label{eq:annular-trimmed-energy}
\int_{\mathcal A_{Ts,R/T}(z)}|A_M|^n\,d\mu
\le C_n\Lambda T^{-n/2}
\Bigl(\|r|A_M|\|_{L^\infty(\{r=s\})}^n
+\|r|A_M|\|_{L^\infty(\{r=R\})}^n\Bigr),
\end{equation}
where $\varepsilon_n,C_n>0$ are constants depending only on $n$.

Indeed, by \eqref{eq:critical-ground},
$(\Delta+\delta_{\rm K}(n)\Pcal)(|A_M|^{\delta_{\rm K}(n)})\ge0$
on $\{|A_M|>0\}$.
By \eqref{eq:annular-small-curvature} and
\eqref{eq:finite-total-radial-barrier}, we can choose $\varepsilon_n>0$ so that
$(\Delta+\delta_{\rm K}(n)\Pcal)r^{-3\delta_{\rm K}(n)/2}<0$ and
$(\Delta+\delta_{\rm K}(n)\Pcal)r^{-\delta_{\rm K}(n)/2}<0$.
For
\[
\begin{aligned}
\psi\coloneqq&\|r|A_M|\|_{L^\infty(\{r=s\})}^{\delta_{\rm K}(n)}
 s^{\delta_{\rm K}(n)/2}r^{-3\delta_{\rm K}(n)/2} +\|r|A_M|\|_{L^\infty(\{r=R\})}^{\delta_{\rm K}(n)}
 R^{-\delta_{\rm K}(n)/2}r^{-\delta_{\rm K}(n)/2},
\end{aligned}
\]
we have $\psi\ge|A_M|^{\delta_{\rm K}(n)}$ on
$\{r=s\}\cup\{r=R\}$. By \eqref{eq:finite-total-radial-barrier} and
\eqref{eq:annular-small-curvature},
\[
(\Delta+\delta_{\rm K}(n)\Pcal)\psi\le0.
\]
Applying the maximum principle on $\mathcal A_{s,R}(z)$ to
$(|A_M|^{\delta_{\rm K}(n)}-\psi)/r^{-\delta_{\rm K}(n)/2}$, we have
\begin{equation}\label{eq:annular-pointwise-decay}
\begin{aligned}
|A_M|^n\le \psi^{n/\delta_{\rm K}(n)} 
\le&2^{n/\delta_{\rm K}(n)-1}
 \Bigl(\|r|A_M|\|_{L^\infty(\{r=s\})}^n s^{n/2}r^{-3n/2} +\|r|A_M|\|_{L^\infty(\{r=R\})}^n R^{-n/2}r^{-n/2}\Bigr).
\end{aligned}
\end{equation}

Integrating $r^{-3n/2}$ over
$\{2^kTs\le r<2^{k+1}Ts,\ r<R/T\}$ and $r^{-n/2}$ over
$\{2^{-k-1}R/T<r\le2^{-k}R/T,\ r>Ts\}$, summing over $k\ge0$,
and using \eqref{eq:annular-volume-assumption}, we obtain
\begin{equation}\label{eq:annular-weighted-integrals}
\begin{aligned}
\int_{\mathcal A_{Ts,R/T}(z)}r^{-3n/2}\,d\mu
&\le2^n\Lambda(Ts)^{-n/2}\sum_{k=0}^\infty2^{-kn/2}\le C_n\Lambda(Ts)^{-n/2},\\
\int_{\mathcal A_{Ts,R/T}(z)}r^{-n/2}\,d\mu
&\le2^{n/2}\Lambda(R/T)^{n/2}\sum_{k=0}^\infty2^{-kn/2} \le C_n\Lambda(R/T)^{n/2}.
\end{aligned}
\end{equation}
We conclude \eqref{eq:annular-trimmed-energy} by integrating
\eqref{eq:annular-pointwise-decay} over $\mathcal A_{Ts,R/T}(z)$
and using \eqref{eq:annular-weighted-integrals}.

For $2s\le u<v\le R/2$, by \eqref{eq:index-stable-ball} and
\eqref{eq:annular-volume-assumption} with $\Lambda=L$,
\begin{equation}\label{eq:index-local-annular-energy}
\begin{aligned}
\int_{\mathcal A_{u,v}(z)}|A_M|^n\,d\mu
\le&C\int_{\mathcal A_{u,v}(z)}r^{-n}\,d\mu\\
=&Cv^{-n}\mu(\mathcal A_{u,v}(z))
 +Cn\int_u^v t^{-n-1}\mu(\mathcal A_{u,t}(z))\,dt\\
\le&CL\bigl(1+n\log(v/u)\bigr).
\end{aligned}
\end{equation}
Now take $\varepsilon=\varepsilon_n$. If $R>8K^2s$, choose regular
values $s'\in(Ks,2Ks)$ and $R'\in(R/(2K),R/K)$. Then
$\int_{\mathcal A_{s',R'}(z)}|A_M|^n\,d\mu\le C_nL\varepsilon_n^n$ by
\eqref{eq:annular-trimmed-energy} with $\Lambda=L$ and $T=1$.
Applying \eqref{eq:index-local-annular-energy} with $(u,v)=(2s,s')$
and $(u,v)=(R',R/2)$, we obtain
\[
\begin{aligned}
\int_{\mathcal A_{2s,R/2}(z)}|A_M|^n\,d\mu
\le&C_nL\varepsilon_n^n
 +CL\Bigl(2+n\log\frac{s'}{2s}+n\log\frac R{2R'}\Bigr)\\
\le&C_nL\varepsilon_n^n+CL(1+\log K).
\end{aligned}
\]
For $R\le8K^2s$, by \eqref{eq:index-local-annular-energy} with
$(u,v)=(2s,R/2)$,
\[
\int_{\mathcal A_{2s,R/2}(z)}|A_M|^n\,d\mu
\le CL\Bigl(1+n\log\frac R{4s}\Bigr)
\le CL(1+\log K).
\]
This proves \eqref{eq:index-stable-annulus}.

We first prove \eqref{eq:index-energy-bound} for $\theta=1$ and all
$\Lambda<\infty$ by induction on $I$.
Applying \eqref{eq:stationary-monotonicity} with center $z\in B_{3/2}$
and outer radius $1/4$, using $B_{1/4}(z)\subset B_2$, we obtain
\eqref{eq:index-local-mass} for $0<t<1/4$; for $t\ge1/4$, use
$\mu(F^{-1}(B_t(z)))\le\Lambda\le4^n\Lambda t^n$. Thus
\begin{equation}\label{eq:index-local-mass}
\mu(F^{-1}(B_t(z)))\le4^n\Lambda t^n
\quad\forall\,z\in B_{3/2},\ \forall\,t>0.
\end{equation}
If $\operatorname{ind}_{\delta_+}(F^{-1}(B_t(z)))=0$ for
$z\in B_{3/2}$ and $0<t<1/256$, choose $z_i\in B_1$ with
\[
B_1\subset\bigcup_{i=1}^N B_{1/1024}(z_i),\qquad N\le C(n).
\]
By \eqref{eq:index-stable-ball} and
\eqref{eq:index-local-mass}, we have
\[
\int_{F^{-1}(B_1)}|A_M|^n\,d\mu\le C.
\]
In particular, \eqref{eq:index-energy-bound} holds for $I=0$.
Otherwise, assume \eqref{eq:index-energy-bound}
for indices at most $I-1$ and let
\[
d\coloneqq \inf\bigl\{t\in(0,1/256):
\exists\,z\in B_{3/2},\quad
\operatorname{ind}_{\delta_+}(F^{-1}(B_t(z)))>0\bigr\}.
\]
By properness, $\sup_{F^{-1}(\overline B_{7/4})}\Pcal<\infty$.
For $z\in B_{3/2}$, $0<t<1/4$, and
$\phi\in C_c^1(F^{-1}(B_t(z)))$, by H\"older's inequality,
\eqref{eq:MS}, and \eqref{eq:index-local-mass},
\[
\begin{aligned}
\delta_+\int_{F^{-1}(B_t(z))}\Pcal\phi^2\,d\mu
\le&\delta_+\Bigl(\sup_{F^{-1}(\overline B_{7/4})}\Pcal\Bigr)
 \int_{F^{-1}(B_t(z))}\phi^2\,d\mu\\
\le&\delta_+\Bigl(\sup_{F^{-1}(\overline B_{7/4})}\Pcal\Bigr)
 \mu(F^{-1}(B_t(z)))^{2/n}\\
&\qquad\times\Bigl(\int_M|\phi|^{2n/(n-2)}\,d\mu\Bigr)^{(n-2)/n}\\
\le&C_n\delta_+\Bigl(\sup_{F^{-1}(\overline B_{7/4})}\Pcal\Bigr)
 (4^n\Lambda)^{2/n}t^2\int_M|\nabla\phi|^2\,d\mu.
\end{aligned}
\]
For sufficiently small $t>0$, uniformly in $z\in B_{3/2}$,
\eqref{eq:delta-stability} with $\delta=\delta_+$ holds on $F^{-1}(B_t(z))$.
By the definition of $d$, $0<d<1/256$. Choose $B_\rho(p)$ and
$\psi\in C_c^1(F^{-1}(B_\rho(p)))$ with
$d\le\rho<2d<1/128$ and $Q_{\delta_+,M}(\psi)<0$.
If $V$ is a negative subspace supported in
$F^{-1}(B_2\setminus\overline B_\rho(p))$, disjoint supports imply
\[
Q_{\delta_+,M}(a\psi+\phi)
=a^2Q_{\delta_+,M}(\psi)+Q_{\delta_+,M}(\phi)
\qquad\forall\,a\in\R,\ \forall\,\phi\in V.
\]
Thus $\R\psi\oplus V$ is negative and
$1+\dim V\le\operatorname{ind}_{\delta_+}(M)\le I$, so
\begin{equation}\label{eq:index-exterior-drop}
\operatorname{ind}_{\delta_+}
\bigl(F^{-1}(B_2\setminus\overline B_\rho(p))\bigr)\le I-1.
\end{equation}
If $|p|\ge5/4$, choose centers $z_i$ with
\[
\begin{gathered}
B_1\subset\bigcup_{i=1}^N B_{1/32}(z_i),\quad N\le C(n), \quad
B_{1/16}(z_i)\subset B_{19/16}\setminus\overline B_\rho(p).
\end{gathered}
\]
Applying the induction hypothesis on $B_{1/16}(z_i)$, using
\eqref{eq:index-exterior-drop} and \eqref{eq:index-local-mass}, and
summing over $i$, we obtain \eqref{eq:index-energy-bound}.

If instead $|p|<5/4$, then each $x\in F^{-1}(B_{4\rho}(p))$ satisfies
$F(x)\in B_{3/2}$ and, by the definition of $d$,
\[
\operatorname{ind}_{\delta_+}\bigl(F^{-1}(B_{d/2}(F(x)))\bigr)=0.
\]
Applying \eqref{eq:index-stable-ball} with $z=F(x)$ and $t=d/4$,
we obtain $|A_M|(x)\le C/d\le2C/\rho$. Hence, by
\eqref{eq:index-local-mass},
\begin{equation}\label{eq:index-core-energy}
\int_{F^{-1}(B_{4\rho}(p))}|A_M|^n\,d\mu
\le C\rho^{-n}\mu(F^{-1}(B_{4\rho}(p)))\le C,
\end{equation}
where $C=C(n,\delta_+,\Lambda)$.
For $2\rho\le s<t\le1/4$ with $t/s\le4$, choose centers $z_i$ with
\[
\begin{gathered}
\{s\le|y-p|\le t\}\subset\bigcup_{i=1}^N B_{s/16}(z_i),
\quad N\le C(n),\quad
B_{s/8}(z_i)\subset B_2\setminus\overline B_\rho(p).
\end{gathered}
\]
By the induction hypothesis on $B_{s/8}(z_i)$,
\eqref{eq:index-exterior-drop}, and \eqref{eq:index-local-mass},
\begin{equation}\label{eq:index-collar-energy}
\int_{\{s<|F-p|<t\}}|A_M|^n\,d\mu\le C,
\quad C=C(n,\delta_+,\Lambda,I-1).
\end{equation}
Fix $R_0\in(1/8,1/4)$ and use the radii $2\rho,4\rho,8\rho,\ldots$
below $R_0$, followed by $R_0$.
Starting with $u=2\rho$, if
$\operatorname{ind}_{\delta_+}(\mathcal A_{u,R_0}(p))=0$, take
$v=w=R_0$ and stop. Otherwise, choose the first listed radius $v>u$ with
$\operatorname{ind}_{\delta_+}(\mathcal A_{u,v}(p))>0$, and let $w$ be the
preceding radius. Then
$\operatorname{ind}_{\delta_+}(\mathcal A_{u,w}(p))=0$ and $v\le2w$.
Begin the next annulus at $v$ and repeat.

For each chosen $\mathcal A_{u,v}(p)$ with
$\operatorname{ind}_{\delta_+}(\mathcal A_{u,v}(p))>0$, choose
$\phi\in C_c^1(\mathcal A_{u,v}(p))$ with $Q_{\delta_+,M}(\phi)<0$.
By the disjoint supports and \eqref{eq:index-exterior-drop}, at most
$I-1$ chosen annuli satisfy
$\operatorname{ind}_{\delta_+}(\mathcal A_{u,v}(p))>0$.
The only possible additional annulus is $\mathcal A_{u,R_0}(p)$ with
$\operatorname{ind}_{\delta_+}(\mathcal A_{u,R_0}(p))=0$.

By analyticity and $\Delta|F-p|^2=2n$,
$\mu\{|F-p|=t\}=0$ for $t>0$.
For each chosen $\mathcal A_{u,v}(p)$, applying
\eqref{eq:index-stable-annulus} to $\mathcal A_{u,w}(p)$ when $w>4u$,
and \eqref{eq:index-collar-energy} to
$\mathcal A_{u,\min\{2u,v\}}(p)$, where $\min\{2u,v\}/u\le2$,
and to $\mathcal A_{\max\{u,w/2\},v}(p)$, where
$v/\max\{u,w/2\}\le4$, we obtain
\[
\begin{aligned}
\int_{\mathcal A_{u,v}(p)}|A_M|^n\,d\mu
&\le\int_{\mathcal A_{2u,w/2}(p)}|A_M|^n\,d\mu
+\int_{\mathcal A_{u,\min\{2u,v\}}(p)}|A_M|^n\,d\mu +\int_{\mathcal A_{\max\{u,w/2\},v}(p)}|A_M|^n\,d\mu\le C,
\end{aligned}
\]
where $\mathcal A_{2u,w/2}(p)=\varnothing$ if $w\le4u$.
Hence, after summation over at most $I$ annuli,
\[
\int_{\mathcal A_{2\rho,R_0}(p)}|A_M|^n\,d\mu\le C,
\quad C=C(n,\delta_+,\Lambda,I).
\]
Choose centers $z_i$ with
\[
\begin{gathered}
B_1\setminus B_{R_0}(p)\subset\bigcup_{i=1}^N B_{1/64}(z_i),
\quad N\le C(n),\quad
B_{1/32}(z_i)\subset B_2\setminus\overline B_\rho(p).
\end{gathered}
\]
Applying the induction hypothesis on $B_{1/32}(z_i)$ and summing over $i$,
we also have
\[
\int_{F^{-1}(B_1\setminus B_{R_0}(p))}|A_M|^n\,d\mu\le C.
\]
Adding \eqref{eq:index-core-energy} and the integral estimates on
$\mathcal A_{2\rho,R_0}(p)$ and
$F^{-1}(B_1\setminus B_{R_0}(p))$, we obtain
\eqref{eq:index-energy-bound} at index $I$.

Finally, for $0<\theta<2$, choose $z_i\in\overline B_\theta$ with
\[
\begin{gathered}
\overline B_\theta\subset\bigcup_{i=1}^{N_\theta}B_{(2-\theta)/4}(z_i),
\quad N_\theta\le C(n)(2-\theta)^{-2n},\\
B_{(2-\theta)/2}(z_i)\subset B_2.
\end{gathered}
\]
Applying \eqref{eq:index-energy-bound} with $\theta=1$ to
$4(F-z_i)/(2-\theta)$ and summing,
we obtain \eqref{eq:index-energy-bound}.
\end{proof}

To apply Proposition~\ref{thm:index-energy} on expanding balls, we first
prove properness and Euclidean volume growth from finite index.

\begin{proposition}\label{prop:properness-volume-growth}
Let $F:M^n\to\C^n$, $n\ge3$, be a complete connected special
Lagrangian immersion without boundary. Assume
$\operatorname{ind}_{\delta_+}(M)<\infty$ and that one of conditions
\textup{(i)}--\textup{(iii)} of Theorem~\ref{thm:intro-finite-index} holds.
Then $F$ is proper and, for some $\Lambda>0$,
\begin{equation}\label{eq:all-dimensional-growth}
\mu\bigl(F^{-1}(B_R(0))\bigr)\le\Lambda R^n\quad\forall\,R\ge1.
\end{equation}
\end{proposition}

\begin{proof}
Suppose first that $\delta_+>D_n$.
Set
\begin{equation}\label{eq:all-dimensional-basic-constants}
a_n^2\coloneqq \min\Bigl\{\frac{n+2}{n-1},
\frac{2n(n-2)}{3(n+1)}\Bigr\},
\quad c_n^2\coloneqq \frac{n-1}{n+2}.
\end{equation}
By Proposition~\ref{prop:stable-intrinsic-ball}, choose
$a_0=a_0(n,\delta_+)\in(0,a_n)$ and fix $a_0<a_1<a_n$. If
$F$ is $\delta_+$-Simons stable on $B_R^g(p)$, then
\begin{equation}\label{eq:large-coefficient-local-curvature}
R|A_M|(p)\le a_0<a_1<a_n.
\end{equation}

By finite index, choose a compact domain $K$ outside which
$F$ is $\delta_+$-Simons stable.
Enlarge $K$ if necessary, fix $p\in K$, translate $F(p)$ to zero, and let
$\rho\coloneqq \dist_g(p,\cdot)$ and $r\coloneqq |F|\le\rho$.
Here $C=C(n,\delta_+,a_1,K,p)$.
Choose $R_K<R_1$ such that $K\subset B_{R_K}^g(p)$ and
$a_0\rho/(\rho-R_K)\le a_1$ for $\rho\ge R_1$. Applying
\eqref{eq:large-coefficient-local-curvature} on
$B_{\rho(x)-R_K}^g(x)$ for $x\notin B_{R_1}^g(p)$, we derive
\begin{equation}\label{eq:all-dimensional-initial-decay}
\rho|A_M|\le a_1.
\end{equation}
Since $r\le\rho$, we also have $r|A_M|\le a_1$.
For $|X|=1$, \eqref{eq:directional-second-form} and
\eqref{eq:all-dimensional-initial-decay}, together with
$\nabla^2(r^2/2)=g+\langle F,A_M\rangle$, imply
$\nabla^2(r^2/2)\ge(1-a_1c_n)g>0$ on $\{\rho\ge R_1\}$.
Fix $\varepsilon=(1-a_1c_n)/2\in(0,1/2)$.
For $x\in M$ with $\rho(x)>R_1$, choose a unit speed minimizing
geodesic $\gamma:[0,\rho(x)]\to M$ from $p$ to $x$.
Since $\rho(\gamma(t))=t$ and $\nabla_{\dot\gamma}\dot\gamma=0$,
\[
\Bigl(\frac{r^2\circ\gamma}{2}\Bigr)''(t)
=\nabla^2(r^2/2)(\dot\gamma(t),\dot\gamma(t))
\ge1-a_1c_n
\qquad\forall\,t\in[R_1,\rho(x)].
\]
At $t=R_1$, using $F(p)=0$ and $|(F\circ\gamma)'|=1$, we have
\[
\begin{aligned}
\Bigl|\Bigl(\frac{r^2\circ\gamma}{2}\Bigr)'(R_1)\Bigr|
=&\bigl|\langle F(\gamma(R_1)),(F\circ\gamma)'(R_1)\rangle\bigr|\\
\le&|F(\gamma(R_1))-F(p)|
\le \int_0^{R_1}|(F\circ\gamma)'(s)|\,ds=R_1.
\end{aligned}
\]
Integrating $\bigl(r^2\circ\gamma/2\bigr)''$ twice over $[R_1,\rho(x)]$
and using $\bigl|\bigl(r^2\circ\gamma/2\bigr)'(R_1)\bigr|\le R_1$,
we obtain
\begin{equation}\label{eq:all-dimensional-properness}
\begin{aligned}
r(x)^2
=&r(\gamma(R_1))^2
+2(\rho(x)-R_1)\Bigl(\frac{r^2\circ\gamma}{2}\Bigr)'(R_1) +2\int_{R_1}^{\rho(x)}(\rho(x)-t)
 \Bigl(\frac{r^2\circ\gamma}{2}\Bigr)''(t)\,dt\\
\ge&\varepsilon(\rho(x)-R_1)^2-2R_1(\rho(x)-R_1) \\
\ge&\varepsilon\rho(x)^2-C\rho(x)-C.
\end{aligned}
\end{equation}
The constant $C$ is independent of $x$ and $\gamma$.
Increase $R_1$, keeping $C$ fixed, so that $C\rho+C\le\varepsilon(1-\varepsilon)\rho^2$ for $\rho\ge R_1$.
For $\rho(x)\ge R_1$, \eqref{eq:all-dimensional-properness} then implies
$r(x)^2\ge\varepsilon^2\rho(x)^2$,
and hence $\varepsilon\rho(x)\le r(x)\le\rho(x)$.
If $r(x)\le R$ and $\rho(x)>R_1$, then
$\rho(x)\le r(x)/\varepsilon\le R/\varepsilon$. Thus, for $R>0$,
\[
F^{-1}(\overline B_R(0))
\subset\overline B_{\max\{R_1,R/\varepsilon\}}^g(p).
\]
By completeness and Hopf--Rinow, $F^{-1}(\overline B_R(0))$ is a closed
subset of a compact intrinsic ball, and hence is compact. Therefore $F$ is proper.

By \eqref{eq:all-dimensional-initial-decay}, $r\le\rho$, and
\eqref{eq:P-estimates},
\begin{equation}\label{eq:all-dimensional-initial-potential}
r^2\Pcal\le\tfrac32a_1^2.
\end{equation}
Since $a_1<a_n$,
\[
\tfrac32\delta_{\rm K}(n)a_1^2
<\delta_{\rm K}(n)(n-2-\delta_{\rm K}(n))
\le\tfrac{3(n-2)^2}{16}.
\]
We may therefore choose $0<c_0<\delta_{\rm K}(n)$ close enough to
$\delta_{\rm K}(n)$ that
$c_0(n-2-c_0)>3\delta_{\rm K}(n)a_1^2/2$.
By \eqref{eq:finite-total-radial-barrier} and
\eqref{eq:all-dimensional-initial-potential}, we have
$(\Delta+\delta_{\rm K}(n)\Pcal)r^{-c_0}<0$ and
$(\Delta+\delta_{\rm K}(n)\Pcal)r^{-3(n-2)/4}<0$.
With $R_0$ and $C>0$ chosen so that
$|A_M|^{\delta_{\rm K}(n)}\le Cr^{-3(n-2)/4}$ on $\{r=R_0\}$,
\[
\frac{|A_M|^{\delta_{\rm K}(n)}-Cr^{-3(n-2)/4}}{r^{-c_0}}
\longrightarrow0\quad\text{as }r\to\infty
\]
by \eqref{eq:all-dimensional-initial-decay}.
Repeating \eqref{eq:radial-comparison-maximum} with denominator
$r^{-c_0}$, we obtain
\[
|A_M|\le Cr^{-3(n+1)/4}\quad(r\ge R_0).
\]
By \eqref{eq:all-dimensional-properness},
$\varepsilon t\le r(\gamma(t))\le t$ for $R_1\le t\le\rho(x)$ along a
unit speed minimizing geodesic from $p$ to $x$.
With $|A_M|\le Cr^{-3(n+1)/4}$, repeating the integration leading to
\eqref{eq:finite-total-radial-gradient}, we obtain
$|\nabla r|^2\ge1-C/r$.
Hence \eqref{eq:all-dimensional-growth} holds by the divergence
and coarea calculation in \eqref{eq:finite-total-derived-volume},
after translating the origin back.

It remains to consider \eqref{eq:low-dimensional-index-range} or
condition~\textup{(ii)} of Theorem~\ref{thm:intro-finite-index}, for which
$3\le n\le7$.
Because $\operatorname{ind}_{\delta_+}(M)<\infty$, there exists a compact
smooth domain $K\subset M$ such that $F$ is $\delta_+$-Simons stable on
$M\setminus K$.

Suppose first that \eqref{eq:low-dimensional-index-range} holds. Let
\[
\chi_n\coloneqq
\frac{(n-1)(n+\sqrt{n+2})}{2(n+1)\sqrt{n+2}}.
\]
At any point of $\{|A_M|>0\}$, fix an arbitrary unit tangent vector $X$.
By the Gauss equation, $0\le-\Ric_g(X,X)\le|A_M|^2$, and
\begin{equation}\label{eq:low-dimensional-ricci}
\begin{aligned}
\frac{-\Ric_g(X,X)}{\Pcal}
\le&\frac{-\Ric_g(X,X)|A_M|^2}
 {\frac{n+1}{n(n-1)}|A_M|^4+\frac{n+2}{n-2}|\mathring{\Ric}_g|^2}\\
\le&\frac{-n(n-1)\Ric_g(X,X)|A_M|^2}
 {(n+1)|A_M|^4+\frac{n+2}{n-2}(n\Ric_g(X,X)+|A_M|^2)^2}\\
=&\chi_n-
\frac{n^2\chi_n\bigl(|A_M|^2+\sqrt{n+2}\,\Ric_g(X,X)\bigr)^2}
 {(n-2)(n+1)|A_M|^4+(n+2)(n\Ric_g(X,X)+|A_M|^2)^2}
\le\chi_n.
\end{aligned}
\end{equation}
In particular, $\Ric_g(X,X)\ge-\chi_n\Pcal$.
In \eqref{eq:low-dimensional-ricci}, we used
\eqref{eq:sharp-curvature-decomposition} in the first inequality and
Cauchy--Schwarz with $\tr_g\mathring{\Ric}_g=0$ in the second.
At $A_M=0$, $\Ric_g=\Pcal=0$.

Let $\omega$ be an $L^2$ harmonic one form. On $\{|\omega|>0\}$,
applying Bochner's formula in a local orthonormal frame $e_i$, we obtain
\begin{equation}\label{eq:harmonic-one-form-bochner}
\begin{aligned}
|\omega|\Delta|\omega|
=&|\nabla\omega|^2-\bigl|\nabla|\omega|\bigr|^2
 +\sum_{i,j}(\Ric_g)_{ij}\omega(e_i)\omega(e_j)\\
\ge&\frac1{n-1}\bigl|\nabla|\omega|\bigr|^2
 -\chi_n\Pcal|\omega|^2,
\end{aligned}
\end{equation}
where we used $|\nabla\omega|^2\ge\frac n{n-1}\bigl|\nabla|\omega|\bigr|^2$
and \eqref{eq:low-dimensional-ricci} in the inequality.
For $\eta\in C_c^1(M\setminus K)$, multiplying
\eqref{eq:harmonic-one-form-bochner} by $\eta^2$, integrating by parts
over $M$, and applying \eqref{eq:delta-stability} to $\eta|\omega|$,
we obtain, by Young's inequality,
\[
\begin{aligned}
&\Bigl(\frac n{n-1}-\frac{\chi_n}{\delta_+}\Bigr)
 \int_M\eta^2\bigl|\nabla|\omega|\bigr|^2\,d\mu\\
\le&\frac{\chi_n}{\delta_+}\int_M|\omega|^2|\nabla\eta|^2\,d\mu\\
&-2\Bigl(1-\frac{\chi_n}{\delta_+}\Bigr)
 \int_M\eta|\omega|\langle\nabla\eta,\nabla|\omega|\rangle\,d\mu\\
\le&\frac12\Bigl(\frac n{n-1}-\frac{\chi_n}{\delta_+}\Bigr)
 \int_M\eta^2\bigl|\nabla|\omega|\bigr|^2\,d\mu
 +C\int_M|\omega|^2|\nabla\eta|^2\,d\mu.
\end{aligned}
\]
Since $\delta_+>(n-1)\chi_n/n$ by \eqref{eq:low-dimensional-index-range},
we obtain from \eqref{eq:MS}
\[
\begin{aligned}
\|\eta|\omega|\|_{L^{2n/(n-2)}(M)}^2
\le&C\int_M|\nabla(\eta|\omega|)|^2\,d\mu\\
\le&C\int_M\eta^2\bigl|\nabla|\omega|\bigr|^2
 +|\omega|^2|\nabla\eta|^2\,d\mu\\
\le&C\int_M|\omega|^2|\nabla\eta|^2\,d\mu,
\end{aligned}
\]
where $C=C(n,\delta_+)$.
At $\{\omega=0\}$, we use the same regularization argument as in
\eqref{eq:regularized-power}--\eqref{eq:zero-level-energy}.

Choose smooth compact domains $K\Subset K'$ and
$\xi\in C^1(M)$ with $\xi=0$ near $K$ and
$\xi|_{M\setminus K'}=1$. Fix $p\in M$ and let
$\zeta_j\in C_c^1(M)$ satisfy $\zeta_j|_{B_j^g(p)}=1$,
$0\le\zeta_j\le1$, and $|\nabla\zeta_j|\le C/j$.
Substituting $\xi\zeta_j$ for $\eta$,
\[
\|\xi\zeta_j|\omega|\|_{L^{2n/(n-2)}(M)}^2
\le C\int_{K'\setminus K}|\omega|^2|\nabla\xi|^2\,d\mu
+\frac C{j^2}\int_M|\omega|^2\,d\mu.
\]
Letting $j\to\infty$ and using $\omega\in L^2(M)$, we find
\[
\|\omega\|_{L^{2n/(n-2)}(M\setminus K')}
\le C\|\omega\|_{L^2(K')}.
\]
On a smooth compact domain $K''$ with $K'\Subset K''$, interior
estimates for $d\omega=d^*\omega=0$ imply
\begin{equation}\label{eq:harmonic-local-estimate}
\|\omega\|_{L^{2n/(n-2)}(K')}
\le C\|\omega\|_{W^{1,2}(K')}
\le C\|\omega\|_{L^2(K'')}.
\end{equation}
The constant in \eqref{eq:harmonic-local-estimate} is $C=C(K',K'')$;
hence, after enlarging $K$ to $K''$,
\begin{equation}\label{eq:harmonic-global-estimate}
\|\omega\|_{L^{2n/(n-2)}(M)}
\le C\|\omega\|_{L^2(K)},
\quad C=C(n,\delta_+,K).
\end{equation}

Set
\[
\mathcal H^1_{(2)}(M)
\coloneqq \{\omega\in L^2(T^*M):d\omega=d^*\omega=0\}.
\]
By unique continuation, the $L^2(K)$ inner product is nondegenerate on
$\mathcal H^1_{(2)}(M)$. If $\dim\mathcal H^1_{(2)}(M)=\infty$, choose an
$L^2(K)$-orthonormal sequence $\omega_j$.
Since $\|\omega_j\|_{L^2(K)}=1$, interior estimates on a fixed
neighborhood of $K$, H\"older's inequality, and
\eqref{eq:harmonic-global-estimate} imply
\[
\|\omega_j\|_{W^{1,2}(K)}
\le C_K\|\omega_j\|_{L^{2n/(n-2)}(M)}\le C_K.
\]
Compactness of $W^{1,2}(K)\hookrightarrow L^2(K)$ then supplies a
subsequence converging in $L^2(K)$, contradicting
$\|\omega_i-\omega_j\|_{L^2(K)}=\sqrt2$.
Thus $\dim\mathcal H^1_{(2)}(M)<\infty$.

Let $\mathcal D^{1,2}(M)$ be the completion of $C_c^\infty(M)$
with respect to the norm $\bigl(\int_M|dv|^2\,d\mu\bigr)^{1/2}$.
By \eqref{eq:MS},
\[
\|v\|_{L^{2n/(n-2)}(M)}\le C\|dv\|_{L^2(M)}
\quad\forall\,v\in\mathcal D^{1,2}(M),\quad C=C(n).
\]
Since $v\mapsto dv$ is an isometry and $\mathcal D^{1,2}(M)$ is complete,
$\{dv:v\in\mathcal D^{1,2}(M)\}$ is closed in $L^2(T^*M)$.
For a smooth compactly supported closed one form $\omega_0$, choose
$u\in\mathcal D^{1,2}(M)$ so that $du$ is the orthogonal projection of
$\omega_0$ onto $\{dv:v\in\mathcal D^{1,2}(M)\}$. Then
\[
\int_M\langle du,dv\rangle\,d\mu
=\int_M\langle\omega_0,dv\rangle\,d\mu
\quad\forall\,v\in\mathcal D^{1,2}(M).
\]
The form $\omega_0-du$ is closed, coclosed, and square integrable, hence a
smooth $L^2$ harmonic one form, and is unchanged under
$\omega_0\mapsto\omega_0+dv$ and $u\mapsto u+v$ for
$v\in C_c^\infty(M)$, so $[\omega_0]\mapsto\omega_0-du$ is a well defined
linear map $H_c^1(M;\R)\to\mathcal H^1_{(2)}(M)$.

To check injectivity, suppose $\omega_0=du$. Choose a compact smooth domain
$K_0\supset\operatorname{supp}\omega_0$ such that every component
of $M\setminus K_0$ has noncompact closure.
For each component $E$ of $M\setminus K_0$, by
\cite[Lemma~2.1 and Remark~2.2]{FloritSimon2026},
\[
\mu(E)=\infty,\quad du|_E=0,\quad
u\in L^{2n/(n-2)}(E)
\quad\Longrightarrow\quad u|_E=0.
\]
Since $du=\omega_0\in C_c^\infty(T^*M)$ and $u=0$ on $M\setminus K_0$,
$u\in C_c^\infty(M)$, so $[\omega_0]=0$ and $[\omega_0]\mapsto\omega_0-du$
is injective. Together with $\dim\mathcal H^1_{(2)}(M)<\infty$, this implies
\begin{equation}\label{eq:low-dimensional-cohomology}
\dim H_c^1(M;\R)<\infty.
\end{equation}
Under condition \textup{(ii)} of Theorem~\ref{thm:intro-finite-index},
\eqref{eq:low-dimensional-cohomology} is assumed.

To prove \eqref{eq:all-dimensional-growth}, choose
$z\in\C^n\setminus F(M)$, which is possible since
$\operatorname{vol}_{2n}(F(M))=0$, and let $r\coloneqq |F-z|>0$.
For $n=3$ and $5\le n\le7$, take $a\coloneqq 1$; for $n=4$, take $a\coloneqq 9/10$.
With $\widetilde g\coloneqq r^{-2a}g$ and $\psi\in C_c^1(M\setminus K)$,
expanding $|\nabla(r^{-a(n-2)/2}\psi)|^2$ and integrating by parts over
$M\setminus K$, we obtain, by \eqref{eq:delta-stability},
\begin{equation}\label{eq:conformal-stability}
\begin{aligned}
&Q_{\delta_+,M}\bigl(r^{-a(n-2)/2}\psi\bigr)\\
=&\int_{M\setminus K}r^{-a(n-2)}|\nabla\psi|^2  - \Bigl(
 r^{-a(n-2)/2}\Delta r^{-a(n-2)/2}
 +\delta_+r^{-a(n-2)}\Pcal\Bigr)\psi^2\,d\mu\\
=&\int_{M\setminus K}
 \bigl(|\widetilde\nabla\psi|^2-V\psi^2\bigr)\,d\widetilde\mu\ge0,
\end{aligned}
\end{equation}
where $d\widetilde\mu=r^{-an}d\mu$,
$|\widetilde\nabla\psi|^2=r^{2a}|\nabla\psi|^2$, and, by $\Delta r^2=2n$,
\[
V \coloneqq \delta_+r^{2a}\Pcal-\frac{an(n-2)}2r^{2a-2} +\Bigl(a(n-2)+\frac{a^2(n-2)^2}{4}\Bigr)
 r^{2a-2}|\nabla r|^2.
\]

When $n=3$, take $3/(4\delta_+)<k<3$. Then
\[
\begin{aligned}
kV+\frac12\operatorname{Scal}_{\widetilde g}
=&r^2\Bigl(k\delta_+\Pcal-\frac12|A_M|^2\Bigr)
 +6-\frac{3k}{2}+\Bigl(\frac{5k}{4}-5\Bigr)|\nabla r|^2\\
\ge&r^2\Bigl(\frac{2k\delta_+}{3}-\frac12\Bigr)|A_M|^2
 +1-\frac k4\ge1-\frac k4>0,
\end{aligned}
\]
where we used \eqref{eq:sharp-curvature-decomposition} and
$|\nabla r|^2\le1$ in the first inequality.

For $g$-orthonormal $X,Y$, define
\begin{equation}\label{eq:weighted-biric-definition}
\operatorname{BiRic}_{\alpha,g}(X,Y)
\coloneqq \Ric_g(X,X)+\alpha\Ric_g(Y,Y)-\alpha K_g(X,Y),
\quad \operatorname{BiRic}_g\coloneqq \operatorname{BiRic}_{1,g}.
\end{equation}
Using $\nabla^2(r^2/2)=g+\langle F-z,A_M\rangle$ and
$\Delta r^2=2n$, we compute
\begin{equation}\label{eq:weighted-conformal-biric}
\begin{aligned}
&\operatorname{BiRic}_{\alpha,\widetilde g}(r^aX,r^aY)\\
=&r^{2a}\operatorname{BiRic}_{\alpha,g}(X,Y)\\
&+ar^{2a-2}\bigl\langle F-z,
(n-2-\alpha)A_M(X,X)+\alpha(n-3)A_M(Y,Y)\bigr\rangle\\
&+r^{2a-2}\Bigl[2a\bigl(n-1+\alpha(n-2)\bigr)\\
&\qquad\qquad-\bigl(2a(1+\alpha)+a^2(n-2+\alpha(n-3))\bigr)|\nabla r|^2\\
&\qquad\qquad+a(a-2)\bigl((n-2-\alpha)|dr(X)|^2
+\alpha(n-3)|dr(Y)|^2\bigr)\Bigr].
\end{aligned}
\end{equation}
When $n=4$, choose $\alpha=1$ and $k,c$ from
\eqref{eq:weighted-dimensional-parameters}. Then
\[
\begin{aligned}
&kV+\operatorname{BiRic}_{\widetilde g}(r^aX,r^aY)\\
\ge&\frac{r^{2a}}{4c}|A_M(X,X)+A_M(Y,Y)|^2
 +ar^{2a-2}\langle F-z,A_M(X,X)+A_M(Y,Y)\rangle\\
&+r^{2a-2}\Bigl(10a-4ak
 +[-6a-2a^2+k(2a+a^2)]|\nabla r|^2\Bigr)\\
=&\frac{r^{2a}}{4c}
 \Bigl|A_M(X,X)+A_M(Y,Y)+\frac{2ac}{r^2}(F-z)^\perp\Bigr|^2
 +a(2-a)(2-k)r^{2a-2}\\
\ge&a(2-a)(2-k)r^{2a-2},
\end{aligned}
\]
where we used \eqref{eq:weighted-conformal-biric},
Lemma~\ref{lem:low-dimensional-cubic-coefficients}, and
$|dr(X)|^2+|dr(Y)|^2\le|\nabla r|^2$ in the first inequality,
$|(F-z)^\perp|^2=r^2(1-|\nabla r|^2)$ and
$c=(6-2k)/a+2-k$ in the equality.
Moreover,
\[
a(2-a)(2-k)-\frac{(1-a)^2}{2}=\frac{59}{8000}>0.
\]

In dimensions $5\le n\le7$, choose $\alpha,k,c$ from
\eqref{eq:weighted-dimensional-parameters}, and let $X,Y$ be
$g$-orthonormal tangent vectors. Since $a=1$,
\begin{equation}\label{eq:weighted-conformal-cross-term}
\begin{aligned}
&\frac{r^2}{4c}\bigl|(n-2-\alpha)A_M(X,X)+\alpha(n-3)A_M(Y,Y)\bigr|^2\\
&+\bigl\langle F-z,(n-2-\alpha)A_M(X,X)+\alpha(n-3)A_M(Y,Y)\bigr\rangle\\
=&\frac{r^2}{4c}\Bigl|(n-2-\alpha)A_M(X,X)+\alpha(n-3)A_M(Y,Y)
 +\frac{2c}{r^2}(F-z)^\perp\Bigr|^2\\
&-c(1-|\nabla r|^2)
\ge-c(1-|\nabla r|^2).
\end{aligned}
\end{equation}
By \eqref{eq:weighted-dimensional-parameters},
\[
\begin{aligned}
&kV+\operatorname{BiRic}_{\alpha,\widetilde g}(rX,rY)\\
\ge&2(n-1)+2\alpha(n-2)-\frac{kn(n-2)}2-c\\
&+\Bigl[c-2(n-1)-\alpha(n-2)
 +k\Bigl(n-2+\frac{(n-2)^2}{4}\Bigr)\Bigr]|\nabla r|^2\\
\ge&\frac3{20},
\end{aligned}
\]
where we used \eqref{eq:weighted-conformal-biric},
\eqref{eq:weighted-conformal-cross-term}, and
Lemma~\ref{lem:low-dimensional-cubic-coefficients} in the first
inequality, and \eqref{eq:weighted-dimensional-parameters} with
$0\le|\nabla r|^2\le1$ in the last.

Choose $0<q_0<1-k/4$ when $n=3$,
$(1-a)^2/2<q_0<a(2-a)(2-k)$ when $n=4$, and
$0<q_0<3/20$ when $5\le n\le7$.
Take $\varepsilon,\vartheta>0$ such that
$4q_0\vartheta>(1+\varepsilon)^2(1-a)^2$, with
$\vartheta<(3-k)/(4-k)$ when $n=3$ and
\begin{equation}\label{eq:mu-bubble-square-parameters}
\begin{gathered}
1-k\vartheta>0,\\
\frac{4-(n-2)\alpha^2}{4(n-1-(n-2)\alpha)}
>\vartheta+\frac{k(1-2\vartheta)^2}{4(1-k\vartheta)}
\qquad(4\le n\le7).
\end{gathered}
\end{equation}
For $n=4$, we take $\vartheta<1/2$ sufficiently close to $1/2$.
For $5\le n\le7$, we take $\vartheta>0$ sufficiently small, using
\[
\frac{4-(n-2)\alpha^2}{4(n-1-(n-2)\alpha)}>\frac k4
\]
from \eqref{eq:weighted-dimensional-parameters}.
By \cite[Proposition~3.4]{BianchiniMariRigoli2016} applied to
\eqref{eq:conformal-stability}, we choose a positive
$v\in C^2(M\setminus K)$ satisfying $-\widetilde\Delta v=Vv$.
Approximating $v$ by a positive smooth function sufficiently closely
in $C^2$ near each point, we obtain, for $g$-orthonormal $X,Y$,
\begin{equation}\label{eq:smooth-separator-supersolution}
\begin{aligned}
-k\frac{\widetilde\Delta v}{v}
+\frac12\operatorname{Scal}_{\widetilde g}
&\ge\frac{1-k/4+q_0}{2} &&(n=3),\\
-k\frac{\widetilde\Delta v}{v}
+\operatorname{BiRic}_{\alpha,\widetilde g}(r^aX,r^aY)
&\ge\frac{a(2-a)(2-k)+q_0}{2}r^{2a-2} &&(n=4),\\
-k\frac{\widetilde\Delta v}{v}
+\operatorname{BiRic}_{\alpha,\widetilde g}(rX,rY)
&\ge\frac{3/20+q_0}{2} &&(5\le n\le7).
\end{aligned}
\end{equation}

Fix $p\in M$. By
\cite[Theorem~1]{AzagraFerreraLopezMesasRangel2007}, there exists a smooth
function $\tau$ such that
\[
r\le\dist_g(p,\cdot)+r(p)+1\le\tau
\le\dist_g(p,\cdot)+r(p)+2,
\quad |\nabla\tau|\le1+\varepsilon.
\]
For $0<t<L$, choose $f$ solving
\begin{equation}\label{eq:mu-bubble-profile-equation}
\begin{gathered}
f'=(1-a)f-\frac{q_0+\vartheta f^2}{1+\varepsilon},
\quad f(0+)=+\infty,\quad f(L-)=-\infty,
\end{gathered}
\end{equation}
where
\[
L\coloneqq \frac{2\pi(1+\varepsilon)}
{\sqrt{4q_0\vartheta-(1+\varepsilon)^2(1-a)^2}}.\]
Since $\tau\ge\dist_g(p,\cdot)+r(p)+1$, by Hopf--Rinow,
$\{R\le\tau\le e^LR\}$ is compact and, for
sufficiently large $R$ with $R$ and $e^LR$ regular,
$\{R\le\tau\le e^LR\}\Subset M\setminus K$ has smooth boundary.
On $\{R<\tau<e^LR\}$, set
\[
h_0\coloneqq \tau^{a-1}f(\log(\tau/R)).
\]
By $f(0+)=+\infty$ and $f(L-)=-\infty$,
\[
\lim_{\tau\downarrow R}h_0=+\infty,\quad
\lim_{\tau\uparrow e^LR}h_0=-\infty.
\]

Let $d\widetilde\sigma$ denote hypersurface measure induced by
$\widetilde g$, and let $\widehat g\coloneqq v^{2k/(n-1)}\widetilde g$.
Since $R$ and $e^LR$ are regular values and $h_0$ diverges at both levels,
choose regular values $R<s_-<s_+<e^LR$ such that, with the
$\widehat g$-unit normal toward increasing $\tau$,
\[
\begin{aligned}
H_{\widehat g}(\{\tau=s_-\})&<h_0v^{-k/(n-1)},\quad 
H_{\widehat g}(\{\tau=s_+\}) >h_0v^{-k/(n-1)}.
\end{aligned}
\]
On each component of $\{s_-\le\tau\le s_+\}$ meeting both boundary
levels, minimize $\mathcal J$ in \eqref{eq:weighted-separator-functional}
over finite perimeter sets $E$ with $\mathbf1_E=1$
near $\{\tau=s_-\}$ and $\mathbf1_E=0$ near $\{\tau=s_+\}$:
\begin{equation}\label{eq:weighted-separator-functional}
\begin{aligned}
\mathcal J(E)\coloneqq&\int_{\partial^*E\cap\{s_-<\tau<s_+\}}v^k\,d\widetilde\sigma
-\int_{E\cap\{s_-<\tau<s_+\}}h_0v^k\,d\widetilde\mu\\
=&\operatorname{Area}_{\widehat g}
 \bigl(\partial^*E\cap\{s_-<\tau<s_+\}\bigr)
-\int_{E\cap\{s_-<\tau<s_+\}}h_0v^{-k/(n-1)}\,d\mu_{\widehat g}.
\end{aligned}
\end{equation}
By \cite[Lemma~2.2]{CecchiniWangXieZhu2026}, $\mathcal J$ admits a
$\mu$-bubble minimizer $E$.
For a connected component $U$ of $\{s_-\le\tau\le s_+\}$ meeting
exactly one boundary level, define
\[
E\cap U\coloneqq
\begin{cases}
U,&U\cap\{\tau=s_-\}\ne\varnothing,\\
\varnothing,&U\cap\{\tau=s_+\}\ne\varnothing.
\end{cases}
\]
Outside $\{s_-<\tau<s_+\}$, take
$\mathbf1_E=\mathbf1_{\{\tau\le s_-\}}$.
The strict mean curvature inequalities at $\{\tau=s_-\}$ and
$\{\tau=s_+\}$ imply $\partial E\subset\{s_-<\tau<s_+\}$.
Since $n\le7$, applying \cite[Lemma~2.2]{CecchiniWangXieZhu2026}
on $\{s_-<\tau<s_+\}$, we conclude that $\partial E$ is a smooth
hypersurface in $(M,\widehat g)$.

Equip each component $\Sigma$ of $\partial E$ with the metric induced by $\widetilde g$.
Let $\nu,B_\Sigma$ be the outer unit normal and second fundamental form of
$\Sigma\subset(M,\widetilde g)$, and let
$H_{\Sigma,\widetilde g}\coloneqq\operatorname{tr}_{\widetilde g}B_\Sigma$. By the first variation of
$\mathcal J$ at $E$,
\[
H_{\Sigma,\widehat g}
=v^{-k/(n-1)}\bigl(H_{\Sigma,\widetilde g}+k\nu(\log v)\bigr)
 =h_0v^{-k/(n-1)},
\]
and hence $H_{\Sigma,\widetilde g}=h_0-k\nu(\log v)$.
The second variation of $\mathcal J$ at $E$, as in
\cite[Lemma~5.2]{HongLiWang2024}, is nonnegative, so for $\psi\in C^1(\Sigma)$,
\[ 
 \int_\Sigma v^k|\nabla_\Sigma\psi|^2\,d\widetilde\sigma 
\ge \int_\Sigma v^k\Big(|B_\Sigma|^2+\Ric_{\widetilde g}(\nu,\nu)
-k\widetilde\nabla^2(\log v)(\nu,\nu)+\nu(h_0)\Big)\psi^2\,d\widetilde\sigma. 
\]
Substituting $\psi=v^{-k/2}\phi$ with $\phi\in C^1(\Sigma)$ and
integrating by parts over $\Sigma$, we obtain, by Cauchy--Schwarz,
\begin{equation}\label{eq:separator-second-variation}
\begin{aligned}
&\frac4{4-k}\int_\Sigma|\nabla_\Sigma\phi|^2\,d\widetilde\sigma\\
\ge&\int_\Sigma |\nabla_\Sigma\phi|^2
+k\phi\langle\nabla_\Sigma\phi,\nabla_\Sigma\log v\rangle -\frac{k(4-k)}4\phi^2|\nabla_\Sigma\log v|^2 \,d\widetilde\sigma\\
\ge&\int_\Sigma\Bigl(-k\frac{\widetilde\Delta v}{v}
+|B_\Sigma|^2+\Ric_{\widetilde g}(\nu,\nu) +k\bigl(\nu(\log v)\bigr)^2
+kH_{\Sigma,\widetilde g}\nu(\log v)+\nu(h_0)\Bigr)\phi^2\,d\widetilde\sigma.
\end{aligned}
\end{equation}
If $n=3$,
\begin{equation}\label{eq:mu-bubble-three-dimensional-square}
\begin{aligned}
&|B_\Sigma|^2+\Ric_{\widetilde g}(\nu,\nu)
 +k\bigl(\nu(\log v)\bigr)^2+kH_{\Sigma,\widetilde g}\nu(\log v)
 -\frac{3-k}{4-k}\bigl(H_{\Sigma,\widetilde g}+k\nu(\log v)\bigr)^2\\
\ge&\tfrac12\operatorname{Scal}_{\widetilde g}-K_\Sigma
 +\tfrac34H_{\Sigma,\widetilde g}^2+k\bigl(\nu(\log v)\bigr)^2+kH_{\Sigma,\widetilde g}\nu(\log v)
 -\frac{3-k}{4-k}\bigl(H_{\Sigma,\widetilde g}+k\nu(\log v)\bigr)^2\\
=&\tfrac12\operatorname{Scal}_{\widetilde g}-K_\Sigma
 +\frac{k}{4(4-k)}\bigl(H_{\Sigma,\widetilde g}+2(k-2)\nu(\log v)\bigr)^2\\
\ge&\tfrac12\operatorname{Scal}_{\widetilde g}-K_\Sigma,
\end{aligned}
\end{equation}
where we used the Gauss equation and $|B_\Sigma|^2\ge H_{\Sigma,\widetilde g}^2/2$ in the first
inequality.

For $n\ge4$, set
\[
\lambda_\Sigma(z)\coloneqq \min_{\substack{e\in T_z\Sigma\\|e|_{\widetilde g}=1}}
 \Ric_\Sigma(e,e),\quad z\in\Sigma.
\]
Choosing a unit eigenvector $e\in T\Sigma$ for
$\lambda_\Sigma$ and applying the Gauss equation, we obtain
\[
\Ric_{\widetilde g}(\nu,\nu)
=\operatorname{BiRic}_{\alpha,\widetilde g}(\nu,e)
-\alpha\lambda_\Sigma+\alpha H_{\Sigma,\widetilde g}B_\Sigma(e,e)-\alpha|B_\Sigma(e,\cdot)|^2.
\]
For $0<\alpha\le1$ as in \eqref{eq:weighted-dimensional-parameters},
by Cauchy--Schwarz on $e^\perp\subset T\Sigma$,
\[
\begin{aligned}
&|B_\Sigma|^2-\alpha|B_\Sigma(e,\cdot)|^2+\alpha H_{\Sigma,\widetilde g}B_\Sigma(e,e)\\
\ge&(1-\alpha)B_\Sigma(e,e)^2+
 \frac{(H_{\Sigma,\widetilde g}-B_\Sigma(e,e))^2}{n-2}+\alpha H_{\Sigma,\widetilde g}B_\Sigma(e,e)\\
=&\frac{n-1-(n-2)\alpha}{n-2}
 \Bigl(B_\Sigma(e,e)-\frac{2-(n-2)\alpha}{2(n-1-(n-2)\alpha)}H_{\Sigma,\widetilde g}\Bigr)^2 +\frac{4-(n-2)\alpha^2}{4(n-1-(n-2)\alpha)}H_{\Sigma,\widetilde g}^2\\
\ge&\frac{4-(n-2)\alpha^2}{4(n-1-(n-2)\alpha)}H_{\Sigma,\widetilde g}^2.
\end{aligned}
\]
By \eqref{eq:mu-bubble-square-parameters}, for $4\le n\le7$,
\[
\begin{aligned}
&\frac{4-(n-2)\alpha^2}{4(n-1-(n-2)\alpha)}H_{\Sigma,\widetilde g}^2
 +k\bigl(\nu(\log v)\bigr)^2+kH_{\Sigma,\widetilde g}\nu(\log v)
 -\vartheta\bigl(H_{\Sigma,\widetilde g}+k\nu(\log v)\bigr)^2\\
\ge&k(1-k\vartheta)\Bigl(\nu(\log v)
 +\frac{1-2\vartheta}{2(1-k\vartheta)}H_{\Sigma,\widetilde g}\Bigr)^2\ge0.
\end{aligned}
\]
Hence, for $4\le n\le7$,
\begin{equation}\label{eq:mu-bubble-higher-dimensional-square}
\begin{aligned}
&|B_\Sigma|^2-\alpha|B_\Sigma(e,\cdot)|^2+\alpha H_{\Sigma,\widetilde g}B_\Sigma(e,e)
+k\bigl(\nu(\log v)\bigr)^2+kH_{\Sigma,\widetilde g}\nu(\log v)
\ge\vartheta\bigl(H_{\Sigma,\widetilde g}+k\nu(\log v)\bigr)^2.
\end{aligned}
\end{equation}
By \eqref{eq:mu-bubble-profile-equation},
\begin{equation}\label{eq:low-dimensional-unified-barrier}
\begin{aligned}
|\widetilde\nabla h_0|
=&r^a\tau^{a-2}
 \bigl|(a-1)f(\log(\tau/R))+f'(\log(\tau/R))\bigr|\,|\nabla\tau|\\
=&\frac{r^a\tau^{a-2}}{1+\varepsilon}
 \bigl(q_0+\vartheta f(\log(\tau/R))^2\bigr)|\nabla\tau|\\
\le&q_0r^a\tau^{a-2}+\vartheta(r/\tau)^a h_0^2\\
\le&q_0r^{2a-2}+\vartheta h_0^2,
\end{aligned}
\end{equation}
where we used $|\nabla\tau|\le1+\varepsilon$ in the first inequality,
and $r\le\tau$, $0<a\le1$ in the last.
For $n=3$, substituting $\phi\equiv1$ and
$h_0=H_{\Sigma,\widetilde g}+k\nu(\log v)$ into
\eqref{eq:separator-second-variation}, and using
\eqref{eq:low-dimensional-unified-barrier},
\eqref{eq:smooth-separator-supersolution},
\eqref{eq:mu-bubble-three-dimensional-square}, $\vartheta<(3-k)/(4-k)$, and
Gauss--Bonnet, we obtain
\[
\frac{1-k/4-q_0}{2}\operatorname{Area}_{\widetilde g}(\Sigma)
\le\int_\Sigma K_\Sigma\,d\widetilde\sigma
=2\pi\chi(\Sigma)\le4\pi.
\]

For $n=4$, by \eqref{eq:separator-second-variation},
\eqref{eq:low-dimensional-unified-barrier},
\eqref{eq:smooth-separator-supersolution}, and
\eqref{eq:mu-bubble-higher-dimensional-square},
\[
\begin{aligned}
\frac4{4-k}\int_\Sigma|\nabla_\Sigma\phi|^2\,d\widetilde\sigma
&\ge\int_\Sigma
\bigl(\frac{a(2-a)(2-k)-q_0}{2}r^{2a-2}-\lambda_\Sigma\bigr)
\phi^2\,d\widetilde\sigma 
 \ge\int_\Sigma
\bigl(\frac{a(2-a)(2-k)-q_0}{2}(e^LR)^{2a-2}-\lambda_\Sigma\bigr)
\phi^2\,d\widetilde\sigma.
\end{aligned}
\]
Hence
\[
\lambda_1\!\Bigl(-\frac4{4-k}\Delta_\Sigma+\lambda_\Sigma\Bigr)
\ge\frac{a(2-a)(2-k)-q_0}{2}(e^LR)^{2a-2}.
\]
Approximating $\lambda_\Sigma$ uniformly from below by smooth functions
and applying \cite[Theorem~5.1]{ChodoshLiMinterStryker2025} after rescaling,
with $4/(4-k)<2$, we obtain
\[
\operatorname{vol}_{3,\widetilde g}(\Sigma)
\le C\bigl(\frac{a(2-a)(2-k)-q_0}{2}(e^LR)^{2a-2}\bigr)^{-3/2}\le C R^{3(1-a)}.
\]

For $5\le n\le7$, after dividing
\eqref{eq:separator-second-variation} by $\alpha$ and using
\eqref{eq:low-dimensional-unified-barrier},
\eqref{eq:smooth-separator-supersolution}, and
\eqref{eq:mu-bubble-higher-dimensional-square},
\[
\int_\Sigma\Bigl(\frac4{\alpha(4-k)}|\nabla_\Sigma\phi|^2
+\lambda_\Sigma\phi^2\Bigr)d\widetilde\sigma
\ge\frac{3/20-q_0}{2\alpha}
\int_\Sigma\phi^2\,d\widetilde\sigma.
\]
Here, by \eqref{eq:weighted-dimensional-parameters},
\[
\frac4{\alpha(4-k)}<\frac{n-2}{n-3}.
\]
Applying \cite[Theorem~1(2) and Remark~1]{AntonelliXu2026} to the
$(n-1)$-dimensional manifold $\Sigma$, we obtain
\[
\operatorname{vol}_{n-1,\widetilde g}(\Sigma)\le C.
\]
Since $g=r^{2a}\widetilde g$ and $r\le e^LR$ on $\Sigma$,
for $3\le n\le7$ we conclude that
\begin{equation}\label{eq:mu-bubble-component-area}
\operatorname{vol}_{n-1,g}(\Sigma)\le CR^{n-1},
\quad C=C(n,\delta_+).
\end{equation}

For sufficiently large $R$, by Sard's theorem choose
$R'\in(R+r(p)+2,R+r(p)+3)$ so that $R'$ and $e^LR'$
are regular values of $\tau$.
Let $E$ minimize $\mathcal J$ with $R'$ in place of $R$.

Since $R'>R+r(p)+2$ and
$\tau\le\dist_g(p,\cdot)+r(p)+2$, we have
$B_R^g(p)\subset\{\tau<R'\}\subset E$.
Let $\Omega$ be the component of $E$ containing $B_R^g(p)$, and
enlarge $\Omega$ by adjoining the components of $M\setminus\Omega$
with compact closure. Then $\Omega\Subset M$ is connected and smooth,
contains $B_R^g(p)$, and has
boundary equal to a finite union of components of $\partial E$;
the components of $M\setminus\Omega$ have noncompact closure.
By \eqref{eq:low-dimensional-cohomology} and
\cite[Lemma~3.2]{Ranieri2026},
\[
\#\pi_0(\partial\Omega)\le1+\dim H_c^1(M;\R)<\infty.
\]
By \eqref{eq:mu-bubble-component-area} with $R'$ in place of $R$,
and since $R'\le R+r(p)+3$,
\[
\operatorname{vol}_{n-1,g}(\partial\Omega)
=\sum_{\Sigma\subset\partial\Omega}
\operatorname{vol}_{n-1,g}(\Sigma)\le CR^{n-1}.
\]
Here $C=C\bigl(n,\delta_+,\dim H_c^1(M;\R)\bigr)$ for all sufficiently
large $R$; since $B_R^g(p)\subset\Omega$, applying the isoperimetric
inequality in \cite[Theorem~2.1]{MichaelSimon1973}, we find
\[
\mu(B_R^g(p))
\le\mu(\Omega)
\le C\operatorname{vol}_{n-1,g}(\partial\Omega)^{n/(n-1)}
\le CR^n.
\]
By
\cite[Theorem~1.1]{FloritSimon2026}, $F$ is proper and
$\mu(F^{-1}(B_R(0)))\le\Lambda R^n$ for some $\Lambda>0$ and all $R\ge1$.
\end{proof}

We now prove Theorem~\ref{thm:intro-finite-index}.

\begin{proof}[Proof of Theorem~\ref{thm:intro-finite-index}]
By Proposition~\ref{prop:properness-volume-growth}, $F$ is proper and
satisfies \eqref{eq:all-dimensional-growth}.
The immersion $R^{-1}F|_{F^{-1}(B_{2R})}$ satisfies
\[
\operatorname{ind}_{\delta_+}
\bigl(F^{-1}(B_{2R}),R^{-2}g\bigr)
\le\operatorname{ind}_{\delta_+}(M),
\quad R^{-n}\mu(F^{-1}(B_{2R}))\le2^n\Lambda.
\]
Since \eqref{eq:index-energy-parameter} holds in cases~\textup{(i)}--\textup{(iii)}
of Theorem~\ref{thm:intro-finite-index}, applying
Proposition~\ref{thm:index-energy} to $R^{-1}F|_{F^{-1}(B_{2R})}$ with
$\theta=1$ and letting $R\to\infty$, we obtain
\[
\int_M|A_M|^n\,d\mu\le C<\infty.
\]
Here $C=C\bigl(n,\delta_+,\Lambda,\operatorname{ind}_{\delta_+}(M)\bigr)$.
Finally, repeating \eqref{eq:finite-total-initial-decay}--
\eqref{eq:finite-total-improved-decay},
\eqref{eq:finite-total-radial-gradient}, and
\eqref{eq:finite-total-sharp-barrier}--\eqref{eq:finite-total-sharp-decay},
we obtain, for some $C,R_0>0$,
\[
|A_M|(x)\le C|F(x)|^{-n-1}\quad\forall\,x\in\{|F|\ge R_0\}.
\]
\end{proof}

We next prove Theorem~\ref{thm:intro-bernstein}.

\begin{proof}[Proof of Theorem~\ref{thm:intro-bernstein}]
Under conditions~\textup{(ii)}, \textup{(iii)}, or \textup{(i)} with
$5\le n\le7$ in Theorem~\ref{thm:intro-bernstein}, stability implies
$\operatorname{ind}_{\delta_+}(M)=0$. Theorem~\ref{thm:intro-finite-index}
then implies that $F$ is proper and $\int_M|A_M|^n\,d\mu<\infty$.
Since $\delta_+>\delta_{\rm K}(n)$, $F$ is flat by
Theorem~\ref{thm:finite-total-rigidity}.
It remains to prove condition~\textup{(i)} for $n=3$ with
$\delta_+\ge(5+3\sqrt5)/40$ and $n=4$ with $\delta_+\ge7/12$.

Let $\pi:\widetilde M\to M$ be the universal covering map. For
$\psi\in C_c^1(\widetilde M)$ and
$\Phi(x)\coloneqq(\sum_{\pi(y)=x}\psi(y)^2)^{1/2}\in\operatorname{Lip}_c(M)$,
by \eqref{eq:lipschitz-stability},
\[
\begin{aligned}
\delta_+\int_{\widetilde M}(\Pcal\circ\pi)\psi^2\,d\widetilde\mu
=&\delta_+\int_M\Pcal\Phi^2\,d\mu 
\le \int_M|\nabla\Phi|^2\,d\mu
\le\int_{\widetilde M}|\nabla\psi|^2\,d\widetilde\mu.
\end{aligned}
\]
Since $F\circ\pi$ is $\delta_+$-Simons stable, we may replace $F$ by
$F\circ\pi$ and assume that $M$ is simply connected.

For $n=3$, suppose that $M$ has at least two ends. By
\eqref{eq:low-dimensional-ricci},
\[
\Ric_g\ge-\frac{5+3\sqrt5}{20}\Pcal g.
\]
Moreover, stability implies
\[
\int_M\Bigl(|\nabla\phi|^2-\frac{5+3\sqrt5}{40}\Pcal\phi^2\Bigr)\,d\mu
\ge0\quad\forall\,\phi\in C_c^1(M).
\]
Applying \cite[Theorem~3.3]{CatinoMariMastroliaRoncoroni2026} with
$V=(5+3\sqrt5)\Pcal/40\in C^{0,1}_{\mathrm{loc}}(M)$ and $\beta=2$,
we obtain
\[
\lambda_{\max}(-\Ric_g)=\frac{5+3\sqrt5}{20}\Pcal
\quad\text{on }M.
\]
At $|A_M|>0$, applying \eqref{eq:low-dimensional-ricci} with $n=3$
to an eigenvector of $-\Ric_g$ corresponding to
$\lambda_{\max}(-\Ric_g)$, we have
\[
\lambda_{\max}(-\Ric_g)=\frac{|A_M|^2}{\sqrt5}.
\]
For any orthonormal tangent frame $X,Y,Z$,
\[
K_g(X,Y)=\tfrac12\operatorname{Scal}_g-\Ric_g(Z,Z)
\le\Bigl(\frac1{\sqrt5}-\frac12\Bigr)|A_M|^2\le0,
\]
where $\operatorname{Scal}_g=-|A_M|^2$ by \eqref{eq:ric-matrices}.
At $A_M=0$, $K_g(X,Y)=0$. By the Cartan--Hadamard theorem,
$M$ is diffeomorphic to $\R^3$, contradicting the assumption of at
least two ends, so $M$ has one end.

For $n=4$, first assume $\sup_M|A_M|<\infty$. Fix $p\in M$ and let
$r\coloneqq|F-F(p)|$ and $\rho\coloneqq(1+r^2)^{-1}$. Using
$\Delta r^2=8$ and $\nabla(r^2/2)=(F-F(p))^\top$, we compute,
for $\phi\in C_c^1(M)$,
\begin{equation}\label{eq:bernstein-rho-poincare}
\begin{aligned}
\int_M\rho\phi^2\,d\mu
\le&\int_M\Bigl(\frac4{1+r^2}
-\frac{3|(F-F(p))^\top|^2}{(1+r^2)^2}\Bigr)\phi^2\,d\mu 
= \int_M|\nabla\phi|^2 
- \Bigl|\nabla\phi
+\frac{(F-F(p))^\top}{1+r^2}\phi\Bigr|^2\,d\mu\\
\le&\int_M|\nabla\phi|^2\,d\mu.
\end{aligned}
\end{equation}
Since $\delta_+\ge7/12>9/16$, by \eqref{eq:delta-stability} and
\eqref{eq:bernstein-rho-poincare},
\[
\begin{aligned}
&\int_M\Bigl(\frac9{16}\Pcal
 +\Bigl(1-\frac9{16\delta_+}\Bigr)\rho\Bigr)\phi^2\,d\mu\\
=&\frac9{16\delta_+}\Bigl(\delta_+\int_M\Pcal\phi^2\,d\mu\Bigr)
 +\Bigl(1-\frac9{16\delta_+}\Bigr)\int_M\rho\phi^2\,d\mu\\
\le&\Bigl(\frac9{16\delta_+}+1-\frac9{16\delta_+}\Bigr)
 \int_M|\nabla\phi|^2\,d\mu
 =\int_M|\nabla\phi|^2\,d\mu.
\end{aligned}
\]
Also, by \eqref{eq:low-dimensional-ricci} and $\chi_4<27/32$,
\begin{equation}\label{eq:bernstein-weighted-ricci}
\Ric_g\ge-\chi_4\Pcal g
>-\Bigl(\frac{27}{32}\Pcal
+\Bigl(\frac32-\frac{27}{32\delta_+}\Bigr)\rho\Bigr)g.
\end{equation}
Since $r\le\dist_g(p,\cdot)$, the metric $\rho g$ is complete;
by $\sup_M|A_M|<\infty$ and \eqref{eq:P-estimates},
$\rho,\Pcal\in L^\infty(M)$. If $M$ has two nonparabolic ends, by \cite[proof of Theorem~1.2]{ChengZhou2009} there exists a
nonconstant bounded harmonic function $u$ with
$\int_M|\nabla u|^2\,d\mu<\infty$ such that, on $\{|\nabla u|>0\}$,
\[
\begin{aligned}
0=&\Delta|\nabla u|^{2/3}
+\Bigl(\frac9{16}\Pcal
+\Bigl(1-\frac9{16\delta_+}\Bigr)\rho\Bigr)|\nabla u|^{2/3}\\
=&\frac23|\nabla u|^{-4/3}
\Bigl(|\nabla^2u|^2-\frac43\bigl|\nabla|\nabla u|\bigr|^2
+\Ric_g(\nabla u,\nabla u) 
 +\Bigl(\frac{27}{32}\Pcal
+\Bigl(\frac32-\frac{27}{32\delta_+}\Bigr)\rho\Bigr)|\nabla u|^2\Bigr)>0,
\end{aligned}
\]
where we used Bochner's formula in the second equality, and
$|\nabla^2u|^2\ge\frac43\bigl|\nabla|\nabla u|\bigr|^2$ with
\eqref{eq:bernstein-weighted-ricci} in the strict inequality.
Hence $M$ has at most one nonparabolic end.
Every end of a complete minimal immersion of
dimension at least three is nonparabolic
\cite[Theorem~3]{LiWang2002}, so $M$ has one end.

For $n=3$, or for $n=4$ with $\sup_M|A_M|<\infty$, $M$ is simply
connected and has one end. Given $\omega\in C_c^\infty(T^*M)$ with
$d\omega=0$, choose $u\in C^\infty(M)$ with $\omega=du$.
Since $du=0$ outside $\operatorname{supp}\omega$ and $M$ has one end,
$u-c\in C_c^\infty(M)$ for some $c\in\R$.
Therefore $\omega=d(u-c)$, and hence $H_c^1(M;\R)=0$.
Applying condition~\textup{(ii)} of
Theorem~\ref{thm:intro-finite-index} to $F\circ\pi:\widetilde M\to\C^n$,
we obtain properness and $\int_{\widetilde M}|A_{\widetilde M}|^n\,d\widetilde\mu<\infty$.
By Theorem~\ref{thm:finite-total-rigidity},
$\pi^*A_M=A_{\widetilde M}=0$, so $A_M=0$.

For $n=4$ without the assumption $\sup_M|A_M|<\infty$, suppose that
$|A_M|(p)>0$ at some $p\in M$. Choose $y_j$ maximizing
$(j-\dist_g(p,\cdot))|A_M|$ on $\overline{B_j^g(p)}$, and let
$\lambda_j\coloneqq|A_M|(y_j)$ and
$s_j\coloneqq(j-\dist_g(p,y_j))/2$.
For $\lambda_j(F-F(y_j))$ with induced metric $\lambda_j^2g$, we have
\begin{equation}\label{eq:bernstein-rescaling-data}
\begin{gathered}
\lambda_j^{-1}|A_M|(y_j)=1,\quad
\sup_{B_{\lambda_js_j}^{\lambda_j^2g}(y_j)}\lambda_j^{-1}|A_M|\le2, \quad 
\lambda_js_j\ge\tfrac j2|A_M|(p)\longrightarrow\infty
\quad\text{as }j\to\infty.
\end{gathered}
\end{equation}
For fixed $R>0$, $R<\lambda_js_j$ for all large $j$ by
\eqref{eq:bernstein-rescaling-data}. By \cite[Lemmas~2.3--2.4]{Cooper2011}, Schauder estimates for
\eqref{eq:minimal-graph-equations} applied to the graph functions of
$\lambda_j(F-F(y_j))$, and the Gauss equation, for each $R>0$, $k\ge0$,
and all large $j$,
\[
\sup_{B_R^{\lambda_j^2g}(y_j)}
 \lambda_j^{-k-1}|\nabla^k A_M|\le C(R,k).
\]
Also, by \cite[Lemmas~2.3--2.4]{Cooper2011} and Bishop volume
comparison, for some $\rho_0>0$ independent of $j$ and all large $j$,
\[
B_R^{\lambda_j^2g}(y_j)
\subset\bigcup_{\nu=1}^{N_j}
 B_{\rho_0}^{\lambda_j^2g}(z_{j,\nu}),\qquad N_j\le C(R),
\]
where $\lambda_j(F-F(y_j))$ is a graph over
$dF(T_{z_{j,\nu}}M)$ on $B_{\rho_0}^{\lambda_j^2g}(z_{j,\nu})$.
Gluing the graphs of $\lambda_j(F-F(y_j))$ on increasing intrinsic
balls as in \cite[proof of Theorem~3.1]{Cooper2011}, we pass to a
smooth complete connected
special Lagrangian immersion $F_\infty:M_\infty^4\to\C^4$.
By \eqref{eq:bernstein-rescaling-data} and smooth convergence,
\[
|A_{M_\infty}|(y_\infty)=1,\qquad \sup_{M_\infty}|A_{M_\infty}|\le2.
\]
For $\phi\in C_c^1(M_\infty)$, let $\phi_j\in C_c^1(M)$ be the
corresponding test functions under the smooth convergence of
$\lambda_j(F-F(y_j))$, extended by zero. Since $\dim M=4$, by
\eqref{eq:delta-stability} and \eqref{eq:P-estimates},
\[
Q_{\delta_+,M_\infty}(\phi)
=\lim_{j\to\infty}\lambda_j^2 Q_{\delta_+,M}(\phi_j)\ge0.
\]
By the $n=4$ case of Theorem~\ref{thm:intro-bernstein} already proved
under $\sup_M|A_M|<\infty$, applied to
$F_\infty$, $A_{M_\infty}=0$, contradicting $|A_{M_\infty}|(y_\infty)=1$, so $F$ is flat.
\end{proof}

\appendix

\section{Coefficient calculations}\label{app:coefficients}

For condition~\textup{(iii)} of Theorem~\ref{thm:intro-bernstein}, we now
construct $D_n$ and prove the local curvature estimate
\eqref{eq:large-coefficient-local-curvature} used in
Proposition~\ref{prop:properness-volume-growth}.

\begin{proposition}\label{prop:stable-intrinsic-ball}
For $n\ge3$, the constant $D_n$ in \eqref{eq:all-dimensional-constant}
satisfies $D_n<14n^3$.
Let $F:M^n\to\C^n$ be a special Lagrangian immersion, and fix
$p\in M$, $R>0$ with $\overline{B_R^g(p)}\Subset M$.
If $F$ is $\delta$-Simons stable on $B_R^g(p)$ and $\delta>D_n$, then
\[
R|A_M|(p)<a_n.
\]
Here $a_n$ is defined in \eqref{eq:all-dimensional-basic-constants}.
\end{proposition}

\begin{proof}
With $a_n,c_n$ from \eqref{eq:all-dimensional-basic-constants} and
$b_n^2\coloneqq n/(n+2)$, choose $\theta$ such that
\[
0<\theta<\frac{\pi}{2b_n},\quad
\frac{\sin(c_n\theta)}{c_n\theta}>\sin(b_n\theta).
\]

Let $y\in\overline{B_R^g(p)}$ satisfy
\[
(R-\dist_g(p,y))|A_M|(y)
=\max_{z\in\overline{B_R^g(p)}}(R-\dist_g(p,z))|A_M|(z).
\]
If $(R-\dist_g(p,y))|A_M|(y)=0$, then $A_M=0$ on $B_R^g(p)$ and
$R|A_M|(p)<a_n$. Otherwise, take
\[
s\coloneqq \frac{\theta}{a_n+\theta}\min\bigl\{(R-\dist_g(p,y))|A_M|(y),a_n\bigr\}.
\]
Replacing $F$ by $|A_M|(y)(F-F(y))$ and $(g,R)$ by
$(|A_M|^2(y)g,|A_M|(y)R)$, we have
$F(y)=0$ and $|A_M|(y)=1$ and, by the maximal choice of $y$,
\[
|A_M|\le1+\theta/a_n\quad\text{on }B_s^g(y).
\]
By the Gauss equation and $|A_M|\le1+\theta/a_n$,
$|K_g|\le|A_M|^2/2\le\bigl(1+\theta/a_n\bigr)^2/2$.
Since $s<\pi/(\sqrt2\bigl(1+\theta/a_n\bigr))$, $\exp_y$ has no conjugate points for
$|v|\le s$.
On $U_s\coloneqq \{v\in T_yM:|v|<s\}$, we use the metric $\exp_y^*g$
and the immersion $F_s\coloneqq F\circ\exp_y$, with second fundamental form $A_{U_s}$
and tensor $h_{U_s}$ defined by \eqref{eq:lagrangian-cubic}. On a smooth domain $\Omega$ with
$\overline{B_s^g(y)}\Subset\Omega\Subset B_R^g(p)$, let $u>0$ be a first
Dirichlet eigenfunction of $-\Delta-\delta\Pcal$. By
\eqref{eq:delta-stability},
\[
-\Delta u-\delta\Pcal u
=\frac{u}{\int_\Omega u^2\,d\mu}
 \int_\Omega\bigl(|\nabla u|^2-\delta\Pcal u^2\bigr)\,d\mu\ge0.
\]
Pulling $-\Delta u\ge\delta\Pcal u$ back by the local isometry
$\exp_y:(U_s,\exp_y^*g)\to M$ and integrating by parts on $U_s$, we obtain
\[
\begin{aligned}
\delta\int_{U_s}\Pcal\phi^2\,d\mu
\le&-\int_{U_s}\frac{\Delta(u\circ\exp_y)}{u\circ\exp_y}\phi^2\,d\mu\\
=&\int_{U_s}|\nabla\phi|^2\,d\mu
 -\int_{U_s}(u\circ\exp_y)^2
 \Bigl|\nabla\!\Bigl(\frac{\phi}{u\circ\exp_y}\Bigr)\Bigr|^2\,d\mu\\
\le&\int_{U_s}|\nabla\phi|^2\,d\mu,
\qquad\forall\,\phi\in C_c^\infty(U_s).
\end{aligned}
\]

Let $G$ be the oriented Gauss map of $F_s$. For any unit tangent
vector $X$ of $U_s$, applying Lemma~\ref{lem:kato-equality-algebra}
to $h_{U_s}$ defined in \eqref{eq:lagrangian-cubic}, we obtain
\[
|dG(X)|^2=-\Ric_{\exp_y^*g}(X,X)\le b_n^2|A_{U_s}|^2.
\]
For the orthogonal projection $\pi$ onto $T_yM$, integrating
$|(G\circ\gamma)'|\le b_n\bigl(1+\theta/a_n\bigr)$ along unit speed radial geodesics of
length at most $s\le\frac{a_n\theta}{a_n+\theta}$, we obtain
\begin{equation}\label{eq:normal-projection-jacobian}
\operatorname{Jac}(\pi|_{T_xU_s})=\langle G(x),G(0)\rangle
\ge\cos(b_n\theta)>0.
\end{equation}
Thus the largest principal angle between
$d(F_s)_x(T_xU_s)$ and $dF_y(T_yM)$ is at most $b_n\theta$.
By Hessian comparison, $\partial U_s$ is strictly convex, so any two
distinct interior points can be joined in $U_s$ by a minimizing geodesic
$\gamma$ of length $0<L\le2s$.
Fix $x\in U_s$ and an arbitrary unit vector $X\in T_xU_s$, and complete
$X$ to an orthonormal basis $e_1=X,e_2,\ldots,e_n$ of $T_xU_s$.
By \eqref{eq:contraction-decomposition},
\begin{equation}\label{eq:directional-second-form}
\begin{aligned}
|A_{U_s}(X,X)|^2
=&(h_{U_s})_{111}^2+\sum_{\alpha=2}^n(h_{U_s})_{11\alpha}^2\\
\le&\frac{n-1}{n+2}|A_{U_s}|^2=c_n^2|A_{U_s}|^2.
\end{aligned}
\end{equation}
Along $\gamma$,
\[
|(F_s\circ\gamma)''|
=|A_{U_s}(\dot\gamma,\dot\gamma)|\le c_n|A_{U_s}|\le c_n\bigl(1+\theta/a_n\bigr).
\]
Since $|(F_s\circ\gamma)'|=1$,
\begin{equation}\label{eq:normal-ball-tangent-angle}
\arccos\bigl\langle(F_s\circ\gamma)'(t),(F_s\circ\gamma)'(L/2)\bigr\rangle
\le c_n\bigl(1+\theta/a_n\bigr)|t-L/2|\qquad\forall\,t\in[0,L].
\end{equation}
If $\pi(F_s(\gamma(0)))=\pi(F_s(\gamma(L)))$, then
$F_s(\gamma(L))-F_s(\gamma(0))$ is normal to $T_yM$.
By \eqref{eq:normal-ball-tangent-angle} and
$c_n\bigl(1+\theta/a_n\bigr) L/2\le c_n\theta<\pi/2$,
\[
\begin{aligned}
L\frac{\sin(c_n\theta)}{c_n\theta}
\le&\frac{2}{c_n\bigl(1+\theta/a_n\bigr)}\sin\frac{c_n\bigl(1+\theta/a_n\bigr) L}{2} 
= \int_0^L\cos(c_n\bigl(1+\theta/a_n\bigr)|t-L/2|)\,dt\\
\le&\bigl\langle F_s(\gamma(L))-F_s(\gamma(0)),
(F_s\circ\gamma)'(L/2)\bigr\rangle\\
\le&|F_s(\gamma(L))-F_s(\gamma(0))|\sin(b_n\theta)
\le L\sin(b_n\theta),
\end{aligned}
\]
contradicting $\sin(c_n\theta)/(c_n\theta)>\sin(b_n\theta)$.
Hence $\pi\circ F_s$ is injective.

For a unit speed radial geodesic $\gamma:[0,s]\to U_s$ with
$\gamma(0)=0$, applying \eqref{eq:normal-ball-tangent-angle} on $[0,t]$
and using $F_s(0)=0$, we obtain
\begin{equation}\label{eq:explicit-normal-ball-chord}
\begin{aligned}
|F_s(\gamma(t))|
\ge&\bigl\langle F_s(\gamma(t)),(F_s\circ\gamma)'(t/2)\bigr\rangle\\
\ge&\int_0^t\cos(c_n\bigl(1+\theta/a_n\bigr)|u-t/2|)\,du
=\frac{2}{c_n\bigl(1+\theta/a_n\bigr)}\sin\frac{c_n\bigl(1+\theta/a_n\bigr) t}{2}.
\end{aligned}
\end{equation}
Let $r\coloneqq |F_s|$, $q\coloneqq \frac{2}{c_n(1+\theta/a_n)}\sin\frac{c_n(1+\theta/a_n)s}{2}$, and let $\omega_n$
be the volume of the Euclidean unit $n$-ball. Since $\{r<t\}\Subset U_s$ for $0<t<q$, integration of $\Delta r^2=2n$
over $\{r<t\}$, together with the injectivity of $\pi\circ F_s$ and
\eqref{eq:normal-projection-jacobian}, yields
\begin{equation}\label{eq:coupled-local-volume}
\begin{aligned}
\mu\{r<t\}&\le\frac{\omega_nt^n}{\cos(b_n\theta)}
&&\forall\,t\in(0,q),\\
\int_{\{r=t\}}|\nabla r|\,d\sigma
&=\frac n t\mu\{r<t\}
&&\text{for a.e. }t\in(0,q).
\end{aligned}
\end{equation}
From $|A_{U_s}|^2(0)=1$,
$\Delta|A_{U_s}|^2\ge-2\bigl(1+\theta/a_n\bigr)^2\Pcal$, and
\cite[Proposition~2.3, equation~(2.7)]{ColdingMinicozzi2005},
\begin{equation}\label{eq:coupled-curvature-mean}
t^{-n}\int_{\{r<t\}}|A_{U_s}|^2\,d\mu
\ge\omega_n-\bigl(1+\theta/a_n\bigr)^2\int_{\{r<t\}}\Pcal K_t(r)\,d\mu,
\end{equation}
where $K_t(r)\coloneqq \frac{2r^{2-n}}{n(n-2)}
-\frac{t^{2-n}}{n-2}+\frac{r^2t^{-n}}n$.
For $0<\xi<1$, applying \eqref{eq:lipschitz-stability} to
$\sqrt{\min\{K_t(r),K_t(\xi t)\}}$ on $\{r<t\}$, extended by zero
for $r\ge t$, and using \eqref{eq:coupled-local-volume}, we obtain
\begin{equation}\label{eq:coupled-truncated-kernel}
\begin{aligned}
 \delta\int_{\{r<t\}}\Pcal\min\{K_t(r),K_t(\xi t)\}\,d\mu 
\le&\frac14\int_{\{\xi t<r<t\}}
 \frac{K_t'(r)^2}{K_t(r)}|\nabla r|^2\,d\mu \\
=&\frac14\int_{\xi t}^t\frac{K_t'(u)^2}{K_t(u)}
 \Bigl(\int_{\{r=u\}}|\nabla r|\,d\sigma\Bigr)\,du\\
\le&\frac{n\omega_n}{4\cos(b_n\theta)}
 \int_{\xi t}^t\frac{K_t'(u)^2}{K_t(u)}u^{n-1}\,du
=\frac{\omega_n I_n(\xi)}{\cos(b_n\theta)},
\end{aligned}
\end{equation}
where
\[
I_n(\xi)\coloneqq
\int_\xi^1\frac{(n-2)(1-z^n)^2}
{z(2-nz^{n-2}+(n-2)z^n)}\,dz.
\]
Integration by parts, using \eqref{eq:coupled-local-volume}, gives
\begin{equation}\label{eq:coupled-kernel-remainder}
\begin{aligned}
\int_{\{r<\xi t\}}(K_t(r)-K_t(\xi t))\,d\mu
&=-\int_0^{\xi t}K_t'(u)\mu\{r<u\}\,du \le\frac{\omega_nt^2\xi^2}{n\cos(b_n\theta)}
\Bigl(1-\frac{2\xi^n}{n+2}\Bigr).
\end{aligned}
\end{equation}
Substituting
$K_t(r)=\min\{K_t(r),K_t(\xi t)\}
+(K_t(r)-K_t(\xi t))\mathbf1_{\{r<\xi t\}}$, valid for $0<r<t$,
into \eqref{eq:coupled-curvature-mean} and using
\eqref{eq:coupled-truncated-kernel}, \eqref{eq:coupled-kernel-remainder},
and $\Pcal\le3\bigl(1+\theta/a_n\bigr)^2/2$, we obtain
\begin{equation}\label{eq:explicit-curvature-mass}
\begin{aligned}
t^{-n}\int_{\{r<t\}}|A_{U_s}|^2\,d\mu
&\ge\omega_n\Bigl(1-\frac{\bigl(1+\theta/a_n\bigr)^2I_n(\xi)}{\delta\cos(b_n\theta)} -\frac{3\bigl(1+\theta/a_n\bigr)^4t^2\xi^2}{2n\cos(b_n\theta)}
 \bigl(1-\frac{2\xi^n}{n+2}\bigr)\Bigr).
\end{aligned}
\end{equation}
Let
$\lambda_n$ be the first Dirichlet eigenvalue of the Euclidean unit
$n$-ball and $\eta$ the decreasing radial first Dirichlet eigenfunction
of the Euclidean $q$-ball. Extending $\eta(r)$ by zero for $r\ge q$,
we have $\eta(r)\in W_0^{1,2}(U_s)$ by
\eqref{eq:explicit-normal-ball-chord}.
Since $\eta$ is decreasing,
\begin{equation}\label{eq:radial-eigenfunction-moment}
(n+2)\int_0^q\eta(t)^2t^{n+1}\,dt
\le nq^2\int_0^q\eta(t)^2t^{n-1}\,dt.
\end{equation}
Multiplying \eqref{eq:explicit-curvature-mass} by
$-t^n(\eta^2)'(t)$, integrating by parts over $(0,q)$, and using
\eqref{eq:radial-eigenfunction-moment} and $\bigl(1+\theta/a_n\bigr) q\le B_\theta\coloneqq 2c_n^{-1}\sin(c_n\theta/2)$, we have
\begin{equation}\label{eq:coupled-weighted-curvature}
\begin{aligned}
&\int_{\{r<q\}}|A_{U_s}|^2\eta(r)^2\,d\mu 
\ge \frac{n\omega_n}{\cos(b_n\theta)}
 \Bigl(Q_n(\theta,\xi)-\frac{\bigl(1+\theta/a_n\bigr)^2I_n(\xi)}{\delta}\Bigr)
 \int_0^q\eta(t)^2t^{n-1}\,dt,
\end{aligned}
\end{equation}
where
\[
Q_n(\theta,\xi)\coloneqq \cos(b_n\theta)
-\frac{3\bigl(1+\theta/a_n\bigr)^2B_\theta^2\xi^2}{2n}
\Bigl(1-\frac{2\xi^n}{n+2}\Bigr).
\]
On the other hand, the coarea formula, \eqref{eq:coupled-local-volume},
and $-(t^{n-1}\eta')'=(\lambda_n/q^2)t^{n-1}\eta$ on $(0,q)$ imply
\begin{equation}\label{eq:coupled-weighted-gradient}
\begin{aligned}
\int_{\{r<q\}}|\nabla\eta(r)|^2\,d\mu
&\le\frac{n\omega_n}{\cos(b_n\theta)}
 \int_0^q\eta'(t)^2t^{n-1}\,dt =\frac{n\omega_n\lambda_n}{q^2\cos(b_n\theta)}
 \int_0^q\eta(t)^2t^{n-1}\,dt.
\end{aligned}
\end{equation}
Substituting $\eta(r(x))$ into
\eqref{eq:delta-stability} and combining
\eqref{eq:coupled-weighted-curvature}--\eqref{eq:coupled-weighted-gradient}
with \eqref{eq:P-estimates}, we derive
\begin{equation}\label{eq:stable-ball-stability-estimate}
(\delta Q_n(\theta,\xi)-\bigl(1+\theta/a_n\bigr)^2I_n(\xi))q^2
\le\frac{n(n-1)\lambda_n}{n+1}.
\end{equation}
Assume
\begin{equation}\label{eq:stable-ball-parameter-conditions}
0<\theta<\frac{\pi}{2b_n},\quad 0<\xi<1,\quad
\frac{\sin(c_n\theta)}{c_n\theta}>\sin(b_n\theta),\quad
Q_n(\theta,\xi)>0
\end{equation}
and let
\[
E_n(\theta,\xi)\coloneqq \frac{\bigl(1+\theta/a_n\bigr)^2}{Q_n(\theta,\xi)}
 \Bigl(\frac{n(n-1)\lambda_n}{(n+1)B_\theta^2}+I_n(\xi)\Bigr).
\]
If $\delta>E_n(\theta,\xi)$ and $s=\frac{a_n\theta}{a_n+\theta}$, then
$q=\frac{a_nB_\theta}{a_n+\theta}$ and
\[
(\delta Q_n(\theta,\xi)-\bigl(1+\theta/a_n\bigr)^2I_n(\xi))q^2
>\frac{n(n-1)\lambda_n}{n+1},
\]
contradicting \eqref{eq:stable-ball-stability-estimate}. Hence
\[
s=\frac{\theta}{a_n+\theta}(R-\dist_g(p,y))|A_M|(y).
\]
By \eqref{eq:stable-ball-stability-estimate}, with
$q/s\ge B_\theta/\theta$ and
$R|A_M|(p)\le(R-\dist_g(p,y))|A_M|(y)=(a_n+\theta)s/\theta$,
\begin{equation}\label{eq:stable-ball-explicit-bound}
\begin{aligned}
R|A_M|(p)
&\le\frac{a_n+\theta}{B_\theta}
 \sqrt{\frac{n(n-1)\lambda_n}
 {(n+1)(\delta Q_n(\theta,\xi)-\bigl(1+\theta/a_n\bigr)^2I_n(\xi))}} <a_n.
\end{aligned}
\end{equation}
Finally, define
\begin{equation}\label{eq:all-dimensional-constant}
D_n\coloneqq \max\Bigl\{\delta_{\rm K}(n),\delta_{\rm cone}(n),
 \inf_{(\theta,\xi)} E_n(\theta,\xi)\Bigr\},
\end{equation}
where the infimum is over pairs $(\theta,\xi)$ satisfying
\eqref{eq:stable-ball-parameter-conditions}.
For $\theta=17/20$ and $\xi=1/(5\sqrt n)$, using $\sin x\le x$,
$\sin x\ge x-x^3/6$, and $\cos x\ge1-x^2/2$, we derive
\[
1+\theta/a_n<\frac94,\quad Q_n(\theta,\xi)>\frac35,
\quad B_\theta>\frac45.
\]
Using $a_n^2\ge1/2$ and $b_n,c_n<1$, we also verify
\eqref{eq:stable-ball-parameter-conditions} for
$(\theta,\xi)=(17/20,1/(5\sqrt n))$.
Moreover,
\[
2-nz^{n-2}+(n-2)z^n
\ge\frac{n-2}{n}(1-z^n)^2\quad\forall\,z\in(0,1),
\]
so $I_n(\xi)\le n\log(1/\xi)\le n^2$.
Using $1-|x|^2$ on $B_1(0)\subset\R^n$, we have
\[
\lambda_n\le
\frac{4\int_{B_1(0)}|x|^2\,dx}{\int_{B_1(0)}(1-|x|^2)^2\,dx}
=\frac{n(n+4)}2.
\]
Consequently,
\[
\begin{aligned}
E_n(\theta,\xi)
&<\frac{135}{16}\Bigl(
\frac{25}{32}n^2\frac{(n-1)(n+4)}{n+1}+n^2\Bigr) &<14n^3.
\end{aligned}
\]
Combining $E_n(\theta,\xi)<14n^3$ with
$\delta_{\rm K}(n)<1$ and $\delta_{\rm cone}(n)<n/4$, we derive
$D_n<14n^3$.

Since $\delta>D_n$, choose $(\theta,\xi)$ satisfying
\eqref{eq:stable-ball-parameter-conditions} and $E_n(\theta,\xi)<\delta$.
Then \eqref{eq:stable-ball-explicit-bound} implies $R|A_M|(p)<a_n$.
\end{proof}

\begin{lemma}\label{lem:low-dimensional-cubic-coefficients}
Let $F:M^n\to\C^n$ be a special Lagrangian immersion, $4\le n\le7$,
and take
\begin{equation}\label{eq:weighted-dimensional-parameters}
\begin{aligned}
n=4:\quad &(\alpha,k,c,\delta_0)
\coloneqq \bigl(1,\tfrac{159}{80},\tfrac{181}{80},\tfrac7{12}\bigr),\\
n=5:\quad &(\alpha,k,c,\delta_0)
\coloneqq \bigl(\tfrac{19}{20},\tfrac{11}{10},5,\tfrac32\bigr),\\
n=6:\quad &(\alpha,k,c,\delta_0)
\coloneqq \bigl(\tfrac78,\tfrac{11}{20},10,\tfrac72\bigr),\\
n=7:\quad &(\alpha,k,c,\delta_0)
\coloneqq \bigl(\tfrac{17}{20},\tfrac{11}{50},\tfrac{33}{2},10\bigr).
\end{aligned}
\end{equation}
Then, for $g$-orthonormal tangent vectors $X,Y$,
\[
\begin{aligned}
&-\operatorname{BiRic}_{\alpha,g}(X,Y)
+\frac{\bigl|(n-2-\alpha)A_M(X,X)+\alpha(n-3)A_M(Y,Y)\bigr|^2}{4c} 
\le k\delta_0\Pcal.
\end{aligned}
\]
\end{lemma}

\begin{proof}
Fix $p\in M$ and arbitrary $g$-orthonormal vectors $X,Y\in T_pM$.
Extending $X,Y$ to an orthonormal basis $e_1=X,e_2=Y,e_3,\ldots,e_n$
of $T_pM$, we first prove
\begin{equation}\label{eq:weighted-cubic-majorant}
\begin{aligned}
&\bigl|(n-2-\alpha)A_M(X,X)+\alpha(n-3)A_M(Y,Y)\bigr|^2 
\le -\beta_1(\Ric_g)_{11}-\beta_2(\Ric_g)_{22}+2\beta_3K_g(X,Y),
\end{aligned}
\end{equation}
where
\[
\begin{aligned}
n=4:\quad &(\beta_1,\beta_2,\beta_3)\coloneqq \bigl(\tfrac23,\tfrac23,\tfrac23\bigr),\\
n=5:\quad &(\beta_1,\beta_2,\beta_3)\coloneqq \bigl(\tfrac{13}{4},3,3\bigr),\\
n=6:\quad &(\beta_1,\beta_2,\beta_3)\coloneqq \bigl(8,6,\tfrac{13}{2}\bigr),\\
n=7:\quad &(\beta_1,\beta_2,\beta_3)\coloneqq (15,11,12).
\end{aligned}
\]
Indeed, by $\sum_{j=1}^n (h_M)_{ijj}=0$ and Cauchy--Schwarz,
\begin{equation}\label{eq:weighted-trace-estimates}
\sum_{j=3}^n (h_M)_{ijj}^2
\ge\frac{((h_M)_{i11}+(h_M)_{i22})^2}{n-2}\quad\forall\,i\in\{1,2\}.
\end{equation}
Combining \eqref{eq:weighted-trace-estimates} with the Gauss equation
and $\beta_1,\beta_2,\beta_1+\beta_2-\beta_3>0$, we obtain
\[
\begin{aligned}
&-\beta_1(\Ric_g)_{11}-\beta_2(\Ric_g)_{22}+2\beta_3K_g(X,Y) -\bigl|(n-2-\alpha)A_M(X,X)+\alpha(n-3)A_M(Y,Y)\bigr|^2\\
\ge&((h_M)_{111},(h_M)_{122})M_1\binom{(h_M)_{111}}{(h_M)_{122}}
 +((h_M)_{222},(h_M)_{112})M_2\binom{(h_M)_{222}}{(h_M)_{112}}\\
&+\sum_{j=3}^n((h_M)_{11j},(h_M)_{22j})M_3\binom{(h_M)_{11j}}{(h_M)_{22j}}.
\end{aligned}
\]
Here the symmetric $2\times2$ matrices $M_1,M_3$ have entries
\[
\begin{aligned}
(M_1)_{11}&\coloneqq \frac{n-1}{n-2}\beta_1-(n-2-\alpha)^2,\\
(M_1)_{12}&\coloneqq \beta_3-\alpha(n-3)(n-2-\alpha)+\frac{\beta_1}{n-2},\\
(M_1)_{22}&\coloneqq \frac{n-1}{n-2}\beta_1+2\beta_2-2\beta_3-\alpha^2(n-3)^2,\\
(M_3)_{11}&\coloneqq 2\beta_1-(n-2-\alpha)^2,\\
(M_3)_{12}&\coloneqq \beta_3-\alpha(n-3)(n-2-\alpha),\\
(M_3)_{22}&\coloneqq 2\beta_2-\alpha^2(n-3)^2.
\end{aligned}
\]
For $M_2$, interchange $\beta_1$ and $\beta_2$, and $n-2-\alpha$ and
$\alpha(n-3)$, in the entries of $M_1$.
Substituting \eqref{eq:weighted-dimensional-parameters} and the values
of $\beta_1,\beta_2,\beta_3$, we have
\[
(M_i)_{11}\ge0,\quad (M_i)_{22}\ge0,\quad
\det M_i\ge0\quad\forall\,i\in\{1,2,3\}.
\]
Hence $M_i\ge0$ for $i=1,2,3$, which proves
\eqref{eq:weighted-cubic-majorant}.

Using
\[
|W_g(X,Y,X,Y)|^2\le\frac{n-3}{4(n-1)}|W_g|^2,
\]
we combine \eqref{eq:weighted-cubic-majorant} with
\eqref{eq:weighted-biric-definition} and apply Cauchy--Schwarz in
\eqref{eq:sharp-curvature-decomposition} to obtain
\[
\begin{aligned}
&-\operatorname{BiRic}_{\alpha,g}(X,Y)
+\frac{\bigl|(n-2-\alpha)A_M(X,X)+\alpha(n-3)A_M(Y,Y)\bigr|^2}{4c}\\
\le&-\Bigl(1+\frac{\beta_1}{4c}\Bigr)(\Ric_g)_{11}
-\Bigl(\alpha+\frac{\beta_2}{4c}\Bigr)(\Ric_g)_{22} +\Bigl(\alpha+\frac{\beta_3}{2c}\Bigr)K_g(X,Y)\\
\le&d_0|A_M|^2+d_1\sqrt{\Pcal|A_M|^2
-\frac{n+1}{n(n-1)}|A_M|^4},
\end{aligned}
\]
where $d_1>0$ and
\[
\begin{aligned}
n=4:\quad &(d_0,d_1^2)\coloneqq \bigl(\tfrac{64}{145},\tfrac{17}{88}\bigr),\\
n=5:\quad &(d_0,d_1^2)\coloneqq \bigl(\tfrac25,\tfrac12\bigr),\\
n=6:\quad &(d_0,d_1^2)\coloneqq \bigl(\tfrac13,\tfrac23\bigr),\\
n=7:\quad &(d_0,d_1^2)\coloneqq \bigl(\tfrac7{24},\tfrac78\bigr).
\end{aligned}
\]
Finally, for $n=4,5,6,7$,
\[
\frac{n+1}{n(n-1)}k^2\delta_0^2-d_0k\delta_0-\frac{d_1^2}{4}>0.
\]
On $\{|A_M|>0\}$, by \eqref{eq:sharp-curvature-decomposition},
\[
\begin{aligned}
&d_0|A_M|^2+d_1\sqrt{\Pcal|A_M|^2
-\frac{n+1}{n(n-1)}|A_M|^4}\\
=&k\delta_0\Pcal
-\Bigl(\frac{n+1}{n(n-1)}k\delta_0-d_0
-\frac{d_1^2}{4k\delta_0}\Bigr)|A_M|^2\\
&-\Bigl(\sqrt{k\delta_0}
 \sqrt{\Pcal-\frac{n+1}{n(n-1)}|A_M|^2}
 -\frac{d_1}{2\sqrt{k\delta_0}}|A_M|\Bigr)^2\\
\le&k\delta_0\Pcal.
\end{aligned}
\]
At $A_M=0$, $\operatorname{BiRic}_{\alpha,g}(X,Y)=\Pcal=0$.
\end{proof}

\section*{On the use of AI.} 

The research problem and all mathematical ideas in this work are the author's own. ChatGPT was used only to assist with calculations, check details, and polish the exposition. All statements, proofs, references, and conclusions were independently verified by the author, who takes full responsibility for the entire content of the manuscript.

\end{document}